\documentclass[11pt]{amsart}
\pdfoutput=1
\usepackage{amsmath,amssymb,amsthm,array,color,multirow,pdflscape,graphicx,pigpen,stmaryrd,comment}
\usepackage[all]{xypic}
\usepackage[normalem]{ulem}

\usepackage{lmodern}

\usepackage{tikz}
\usepackage{tikz-cd}
\tikzset{shorten <>/.style={shorten >=#1,shorten <=#1}}

\usepackage{float}

\usepackage{todonotes}

\usepackage[margin=1.5in]{geometry}

\usepackage[all,arc,2cell]{xy}
\UseAllTwocells

\usepackage{microtype}

\usepackage{xcolor}
\definecolor{darkgreen}{rgb}{0,0.50,0} 
\definecolor{darkred}{rgb}{0.75,0,0}
\definecolor{darkblue}{rgb}{0,0,0.6} 
\usepackage[pdfborder=0,pagebackref,colorlinks,citecolor=darkgreen,linkcolor=darkgreen,urlcolor=darkblue]{hyperref}
\renewcommand*{\backref}[1]{}
\renewcommand*{\backrefalt}[4]{({%
    \ifcase #1 Not cited.%
          \or On p.~#2%
          \else On pp.~#2%
    \fi%
    })}

\def\makeautorefname#1#2{\expandafter\def\csname#1autorefname\endcsname{#2}}
\makeautorefname{equation}{Equation}%
\makeautorefname{footnote}{footnote}%
\makeautorefname{item}{item}%
\makeautorefname{figure}{Figure}%
\makeautorefname{table}{Table}%
\makeautorefname{part}{Part}%
\makeautorefname{appendix}{Appendix}%
\makeautorefname{chapter}{Chapter}%
\makeautorefname{section}{Section}%
\makeautorefname{subsection}{Section}%
\makeautorefname{subsubsection}{Section}%
\makeautorefname{paragraph}{Paragraph}%
\makeautorefname{subparagraph}{Paragraph}%
\makeautorefname{theorem}{Theorem}%
\makeautorefname{thm}{Theorem}%
\makeautorefname{addm}{Addendum}%
\makeautorefname{mainthm}{Main theorem}%
\makeautorefname{corollary}{Corollary}%
\makeautorefname{cor}{Corollary}%
\makeautorefname{lemma}{Lemma}%
\makeautorefname{lem}{Lemma}%
\makeautorefname{sublemma}{Sublemma}%
\makeautorefname{sublem}{Sublemma}%
\makeautorefname{subl}{Sublemma}%
\makeautorefname{prop}{Proposition}%
\makeautorefname{property}{Property}
\makeautorefname{pro}{Property}
\makeautorefname{sch}{Scholium}%
\makeautorefname{step}{Step}%
\makeautorefname{conject}{Conjecture}%
\makeautorefname{conj}{Conjecture}%
\makeautorefname{questn}{Question}
\makeautorefname{quest}{Question}
\makeautorefname{qn}{Question}
\makeautorefname{definition}{Definition}%
\makeautorefname{defn}{Definition}%
\makeautorefname{defi}{Definition}%
\makeautorefname{def}{Definition}%
\makeautorefname{dfn}{Definition}%
\makeautorefname{notation}{Notation}
\makeautorefname{notn}{Notation}
\makeautorefname{rem}{Remark}%
\makeautorefname{rems}{Remarks}%
\makeautorefname{rmk}{Remark}%
\makeautorefname{rk}{Remark}%
\makeautorefname{remarks}{Remarks}%
\makeautorefname{rems}{Remarks}%
\makeautorefname{rmks}{Remarks}%
\makeautorefname{rks}{Remarks}%
\makeautorefname{example}{Example}%
\makeautorefname{examp}{Example}%
\makeautorefname{exmp}{Example}%
\makeautorefname{exam}{Example}%
\makeautorefname{exa}{Example}%
\makeautorefname{axiom}{Axiom}%
\makeautorefname{axi}{Axiom}%
\makeautorefname{ax}{Axiom}%
\makeautorefname{case}{Case}%
\makeautorefname{claim}{Claim}%
\makeautorefname{clm}{Claim}%
\makeautorefname{assumpt}{Assumption}%
\makeautorefname{asses}{Assumptions}%
\makeautorefname{conclusion}{Conclusion}%
\makeautorefname{concl}{Conclusion}%
\makeautorefname{conc}{Conclusion}%
\makeautorefname{cond}{Condition}%
\makeautorefname{const}{Construction}%
\makeautorefname{con}{Construction}%
\makeautorefname{criterion}{Criterion}%
\makeautorefname{criter}{Criterion}%
\makeautorefname{crit}{Criterion}%
\makeautorefname{exercise}{Exercise}%
\makeautorefname{exer}{Exercise}%
\makeautorefname{exe}{Exercise}%
\makeautorefname{problem}{Problem}%
\makeautorefname{problm}{Problem}%
\makeautorefname{prob}{Problem}%
\makeautorefname{prob}{Problem}%
\makeautorefname{soln}{Solution}%
\makeautorefname{sol}{Solution}%
\makeautorefname{sum}{Summary}%
\makeautorefname{oper}{Operation}%
\makeautorefname{obs}{Observation}%
\makeautorefname{ob}{Observation}%
\makeautorefname{conv}{Convention}%
\makeautorefname{cvn}{Convention}%
\makeautorefname{warn}{Warning}%
\makeautorefname{note}{Note}%
\makeautorefname{fact}{Fact}%
\makeautorefname{ouch0}{Counterexample}%
\newtheorem{thm}{Theorem}[section]
\newtheorem{cor}{Corollary}[section]
\newtheorem{prop}{Proposition}[section]
\newtheorem{lem}{Lemma}[section]

\theoremstyle{definition}
\newtheorem{defn}{Definition}[section]
\newtheorem{con}{Construction}[section]

\newtheorem{exmp}{Example}[section]
\newtheorem{notn}{Notation}[section]

\newtheorem{rem}{Remark}[section]

\numberwithin{equation}{section}
\numberwithin{figure}{section}

\makeatletter
\let\c@cor=\c@thm
\let\c@con=\c@thm
\let\c@prop=\c@thm
\let\c@lem=\c@thm
\let\c@conj=\c@thm
\let\c@defn=\c@thm
\let\c@df=\c@thm
\let\c@exmp=\c@thm
\let\c@rem=\c@thm
\let\c@notn=\c@thm
\let\c@sch=\c@thm
\let\c@equation\c@thm
\makeatother

\newcommand{\R}{\mathbb R}
\newcommand{\C}{\mathbb C}
\newcommand{\N}{\mathbb N}
\newcommand{\Z}{\mathbb Z}

\newcommand{\Map}{\textup{Map}}

\newcommand{\Emb}{\textup{Emb}}

\newcommand{\NTub}{\textup{NTub}}

\newcommand{\id}{\textup{id}}
\newcommand{\colim}{\textup{colim}\,}
\newcommand{\hocolim}{\textup{hocolim}\,}

\newcommand{\ra}{\longrightarrow}

\newcommand{\ti}{\widetilde}
\newcommand{\simar}{\overset\sim\longrightarrow}
\newcommand{\congar}{\overset\cong\longrightarrow}

\newcommand{\mc}{\mathcal}
\newcommand{\op}{\textup{op}}

\newcommand{\Diff}{\textup{Diff}\,}
\newcommand{\tipartial}{\widetilde\partial}

\newcommand{\Fun}{\textup{Fun}}

\newcommand{\Aut}{\textup{Aut}}

\newcommand{\sSet}{\textup{sSet}}

\newcommand{\mirror}{{mirror}}

\newcommand{\Hiso}{\mathcal{H}_{\iso}}

\newcommand{\Man}{\mc M^\mathrm{Emb}}

\newcommand{\Manst}{\mc M^\mathrm{Stab}}
\newcommand{\Manstr}{\mc M^\mathrm{Strat}}
\newcommand{\Hanst}{\mc H^\mathrm{Stab}}
\newcommand{\Ex}{\mathrm{Ex}}

\newcommand{\iso}{\mathrm{iso}}
\newcommand{\eq}{\mathrm{eq}}
\newcommand{\dashiso}{\textrm{-}\mathrm{iso}}
\newcommand{\dasheq}{\textrm{-}\mathrm{eq}}

\newcommand{\St}{\mathrm{St}}

\newcommand{\Wh}{\mathrm{Wh}}
\newcommand{\GL}{\mathrm{GL}}

\newcommand{\A}{\textbf{A}}

\newcommand{\po}{\ar@{}[dr]|{\text{\pigpenfont R}}}
\newcommand{\pb}{\ar@{}[dr]|{\text{\pigpenfont J}}}

\newcommand{\sbt}{\,\begin{picture}(-1,1)(0.5,-1)\circle*{1.8}\end{picture}\hspace{.05cm}}

\newcommand{\ourdefn}[1]{{\it #1}}

\newcommand{\rh}{encasing function }
\newcommand{\therh}{the encasing function }
\newcommand{\rhperiod}{encasing function. }

\newcommand{\arhperiod}{an encasing function. }

\newcommand{\rhs}{encasing functions }

\newlength{\storeparskip}
\title{The splitting of the stable space of equivariant $h$-cobordisms}
\author{Thomas Goodwillie}
\address{Department of Mathematics, Brown University}
\email{thomas\_goodwillie@brown.edu}
\author{Kiyoshi Igusa}
\address{Department of Mathematics, Brandeis University}
\email{igusa@brandeis.edu}
\author{Cary Malkiewich}
\address{Department of Mathematics, Binghamton University}
\email{cmalkiew@binghamton.edu}
\author{Mona Merling}
\address{Department of Mathematics, The University of Pennsylvania}
\email{mmerling@math.upenn.edu}

\begin{document}

\maketitle

%\setcounter{tocdepth}{2}
%\tableofcontents

\begin{abstract}
We prove that the space of equivariant $h$-cobordisms on a smooth $G$-manifold $M$ has a tom Dieck style splitting, after stabilization. As an intermediate step, we introduce a space of isovariant $h$-cobordisms on $M$ and prove that this space also has a splitting, even before stabilizing. This completes the proof of the equivariant version of Waldhausen's stable parametrized $h$-cobordism theorem. Along the way, we establish several foundational results on equivariant smooth manifolds with corners.
\end{abstract}

\begingroup%
\setlength{\parskip}{\storeparskip}% Restore \parskip within this scope
%code added in to indent subsections in table of contents
\makeatletter
\def\l@subsection{\@tocline{2}{0pt}{2.5pc}{5pc}{}}
\def\l@subsubsection{\@tocline{2}{0pt}{5pc}{7.5pc}{}}
\makeatother
% end indentation code
\tableofcontents
\endgroup%

\section{Introduction}

%\cmar{Cary's comments}
%\mmar{Mona's comments}
%\tmar{Tom's comments}
%\kmar{Kiyoshi's comments}

A common technique for constructing diffeomorphisms between smooth compact manifolds $M$ and $N$ is to show first that they are $h$-cobordant. When the dimension is at least $5$, the only obstruction to making the $h$-cobordism trivial (and therefore giving a diffeomorphism $M \cong N$) lives in the Whitehead group $\Wh(\pi_1M)$, the quotient of $K_1(\Z[\pi_1 M]) = \GL_\infty(\Z[\pi_1 M])^{ab}$ by the subgroup $\{\pm 1\} \times \pi_1(M)^{ab}$. In fact, Whitehead torsion gives a bijection between the set of diffeomorphism classes of $h$-cobordisms and $\Wh(\pi_1M)$. This is known as the $s$-cobordism theorem \cite{barden, mazur, stallings, BMS}, or in the case of a simply connected manifold the $h$-cobordism theorem \cite{smale}.

One can make a space $\mathcal{H}(M)$ of $h$-cobordisms on $M$ such that $\pi_0\mathcal H(M)$ is in bijection with the set of diffeomorphism classes of $h$-cobordisms. More generally, homotopy classes of maps of $X$ into $\mathcal{H}(M)$ correspond to isomorphism classes of families of $h$-cobordisms parametrized by $X$ (bundles of $h$-cobordisms on $X$). This space plays a crucial role in determining the homotopy type of the diffeomorphim group $\Diff(M)$ and the homotopy type of the moduli space of smooth manifolds $\coprod_{[M]} B\Diff(M)$. 

The $s$-cobordism theorem computes the set $\pi_0\mathcal{H}(M)$ in terms of $K$-theory when the dimension of $m$ is at least $5$. The stable parametrized $h$-cobordism theorem \cite{wjr} uses Waldhausen's algebraic $K$-theory of spaces to extend this relationship to the higher homotopy groups of $\mathcal{H}(M)$, in a stable range \cite{ww}. To be specific, there is a cofiber sequence of spectra
\[ \Sigma^\infty_+ M \to \A(M) \to \Wh(M), \]
and $\Omega^{\infty+1}  \Wh(M)$ is equivalent to a stabilized space of $h$-cobordisms $\mathcal{H}^\infty(M)$, obtained from $\mathcal{H}(M)$ by multiplying the manifold $M$ by closed disks of arbitrarily high dimension. In a stable range, first computed in \cite{igusa}, $\mathcal{H}^\infty(M)$ agrees with the original space $\mathcal{H}(M)$, and therefore in that range one has a calculation of $\mathcal{H}(M)$ in terms of algebraic $K$-theory. From $\pi_1$ of this cofiber sequence, one recovers the statement that $\pi_0\mathcal{H}(M) \cong \Wh(\pi_1M)$.

Suppose now that a finite group $G$ acts on the smooth compact manifold $M$. Then there is a space of equivariant $h$-cobordisms $\mathcal{H}_G(M)$, constructed in \cite{GIMM1}, with a stabilization $\mathcal{H}_G^\mathcal{U}(M)$. The stabilization involves multiplying the manifold $M$ by representation disks $D(V)$; the superscript $\mathcal U$ refers to a complete universe for $G$. In the case of the trivial group $G = \{e\}$, this is a specific construction of  $\mathcal{H}^\infty(M)$.

The functor $\mathcal{H}_G^\mathcal{U}(M)$ extends in a universal way to a functor on all $G$-spaces, which we denote by $\widetilde{\mc H}_G^\mathcal{U}(X)$ \cite[Thm 7.14]{GIMM1}. In the case where $G$ is trivial, we denote this extension of $\mathcal{H}^\infty(M)$ to all spaces by $\widetilde{\mc H}^\infty(X)$.

In this paper we give the last piece of the proof that this equivariant stabilized $h$-cobordism space is again expressible in terms of algebraic $K$-theory. This leads to an equivariant version of the $s$-cobordism theorem and its generalization to the higher homotopy groups. The main result is that $\mathcal{H}_G^\mathcal{U}(M)$ has a splitting over conjugacy classes of subgroups of $G$, formally analogous to the tom Dieck splitting for fixed points of equivariant suspension spectra \cite{tD}.

 \begin{thm}[\autoref{stable_splitting_thm}]\label{mainthm}
 Let $G$ be a finite group, and let $M$ be a compact smooth $G$-manifold.  Then the space of stable equivariant $h$-cobordisms (stabilized with respect to representation disks) has a splitting
 $$\mathcal{H}_G^\mathcal{U}(M)\simeq \prod_{(H) \leq G} \widetilde{\mc H}^\infty(M^H_{hWH} ),$$ where $M^H_{hWH}$ is the homotopy orbit space for the action of  the Weyl group $WH$ on the fixed point space $M^H$.
 \end{thm}
\noindent
We also prove that the splitting extends to all $G$-spaces.
%If we instead use the functor $\widetilde{\mc H}_G^\mathcal{U}$ that is defined on all $G$-spaces, then the above splitting becomes natural as well: %With this language it becomes possible to talk about naturality, and in fact we prove that the above splitting is natural, in other words it is a zig-zag of weak equivalences of functors.
\begin{thm}[\autoref{natural_stable_splitting_thm}]\label{mainthmspaces}
	There is a (natural) equivalence of functors of $G$-spaces
	$$\widetilde{\mc H}_G^\mathcal{U}(X)\simeq \prod_{(H) \leq G} \widetilde{\mc H}^\infty(X^H_{hWH} ).$$
\end{thm}

When combined with the results of \cite{CaryMona,CaryMona3}, \autoref{mainthm} finishes the proof of the equivariant version of the stable parametrized $h$-cobordism theorem, providing the geometric core of the result. To see why, recall that the equivariant version of algebraic $K$-theory of spaces, denoted $A_G(X)$ has a tom Dieck style splitting for the fixed point spectrum $A_G(X)^G$  \cite{CaryMona}\footnote{It is important to recall that the $G$-spectrum $\A_G(X)$ is not the ``$K$-theory of group actions" for the $G$-action on the category of retractive spaces over $X$. That $G$-spectrum, called $\A_G^\mathrm{coarse}(X)$ in \cite{malkiewich_merling_coassembly}, lacks such a splitting.}. In \cite{CaryMona3} the third and fourth authors construct a map of $G$-spectra $\A_G(M)\to \Sigma^\infty_G M$, for which the $G$-fixed point spectrum of the fiber is equivalent at the infinite-loopspace level to 
$\prod_{(H) \leq G} \widetilde{\mc H}^\infty(M^H_{hWH})$, a product of stable \emph{nonequivariant} $h$-cobordism spaces. Therefore by \autoref{mainthmspaces}, these spaces re-assemble into the stabilized space of \emph{equivariant} $h$-cobordisms, which proves:

\begin{thm}[Equivariant stable parametrized $h$-cobordism theorem]\label{esphcthm} For every compact smooth $G$-manifold $M$, there is a cofiber sequence of $G$-spectra
$$\Sigma^\infty_{G} M_+ \to \A_G(M)\to \Wh_G(M)$$
such that $\Omega^{\infty+1}  \Wh_G(M)^G \simeq \mathcal{H}_G^\mathcal{U}(M)$.
\end{thm}

\begin{rem}
	It appears that \autoref{esphcthm} may be combined with the splitting $\A_G(M) \to \Sigma^\infty_G M_+$ constructed in \cite{hilman_ramzi}, and in unpublished work of Calle, Chan, Chedalavada, and Mejia, to show the stronger result that there is a decomposition $\A_G(M)\simeq \Sigma^\infty_G M_+ \times \Wh_G(M)$, as in the nonequivariant case \cite{wjr}. This requires showing that the form of equivariant $A$-theory used in \cite{hilman_ramzi} is indeed equivalent to $\A_G(M)$. (The geometric core of this stronger result would also be provided by \autoref{mainthmspaces}.) 
	We remark that prior to \cite{hilman_ramzi} it was not clear whether such a splitting should exist. Indeed, the existence of this splitting implies that the transfers in $A$-theory and stable homotopy commute with the projection from $A$-theory to stable homotopy, as in \cite[Thm 4.2]{raptis}, for every covering space of the form $X \to X_{hG}$ with $G$ a finite group:
	\[ \begin{tikzcd}
		A(X_{hG}) \dar \rar & A(X) \dar \\
		\Sigma^\infty_+ X_{hG} \rar & \Sigma^\infty_+ X.
	\end{tikzcd} \]
	But this was not known: \cite{raptis} relies on the functoriality of the Becker--Gottlieb transfer, which was still open \cite{klein_malkiewich_corrigendum}. In fact, it has recently been established that the above commutativity statement fails to hold if $X \to X_{hG}$ is replaced by an arbitrary perfect fibration, thereby proving that the Becker-Gottlieb transfer is not functorial in general \cite{lehner_malkiewich_naef_bg_not,naef_stoll_bg_not}; see also \cite{lind2019transfer,ccry}.
 \end{rem}

Another application of the results of this paper and \cite{GIMM1} occurs in work of the first two authors on equivariant higher Reidemeister torsion. Higher torsion classes are certain real characteristic classes for bundles of smooth compact manifolds, in particular for bundles of smooth $h$-cobordisms. Igusa \cite[Theorem 6.4.2]{igusa2002higher} uses higher torsion to detect the non-triviality of some bundles of $h$-cobordisms over spheres which had been constructed by Hatcher. In terms of spaces of $h$-cobordisms, Hatcher’s construction yields a map $G/O\to \mathcal H^\infty(\ast)$, and one can prove the non-triviality of the construction by showing that universal torsion classes in $H^{4j}(\mathcal H^\infty(\ast))$ pull back nontrivially to $H^{4j}(G/O)\cong H^{4j}(BO)$.
Goodwillie and Igusa (unpublished, extending work of \cite{goodwillie_igusa_ohrt}) generalize these ideas by defining equivariant torsion classes for bundles of smooth compact $G$-manifolds. The non-triviality of an equivariant version of Hatcher’s construction can be detected by pulling back universal equivariant torsion classes from $\mathcal H^{\mathcal U}_G(\ast)$ to the fixed points in an equivariant model for $G/O$. The splitting of $\mathcal H^{\mathcal U}_G(\ast)$ plays a key role in this computation.

In order to prove \autoref{mainthm}, we need to use the space of \emph{isovariant} $h$-cobordisms $\mathcal{H}_{G\dashiso}(M)$. An equivariant map is isovariant if it preserves fixed point strata, and it is an isovariant homotopy equivalence if it has an inverse up to isovariant homotopy. An equivariant $h$-cobordism is said to be an isovariant $h$-cobordism if the inclusions of the top and bottom are isovariant homotopy equivalences. As a part of the technical work in this paper, we extend our treatment of equivariant $h$-cobordisms in \cite{GIMM1} to isovariant $h$-cobordisms.
 
Then our strategy is to prove that the isovariant $h$-cobordism space $\mathcal{H}_{G\dashiso}(M)$ has a splitting over conjugacy classes of subgroups, even before stabilizing. Let us state the result here. For a compact smooth $G$-manifold $M$, let $M_H$ be the submanifold consisting of points with isotropy group $H$. We use a compactification of $M_H$ that is obtained by successively removing open tubular neighborhoods of fixed-point submanifolds with larger isotropy, or equivalently by iterated spherical blowup along fixed-point submanifolds for those subgroups. This produces a $G$-manifold with corners $M_{[H]}$ whose interior can be identified with $M_H$.
 
 \begin{thm}[\autoref{unstable_splitting_thm}]\label{unstablemainthm}
 Let $G$ be a finite group, and let $M$ be a compact smooth $G$-manifold in which $\dim M^H \geq 5$ for all $H \leq G$. Then the unstable space of isovariant $h$-cobordisms satisfies
 $$\mathcal{H}_{G\dashiso}(M)\simeq \prod_{(H) \leq G} \mathcal{H}(M_{[H]}/WH ).$$
 \end{thm}

An equivariant $h$-cobordism on $M$ is not always an isovariant $h$-cobordism \cite{luck}, but is if the $G$-manifold satisfies the \emph{weak gap hypothesis}, meaning that the difference between dimensions of successive fixed point strata is at least 3. It follows that after stabilization with respect to representations of $G$ the spaces of isovariant and equivariant $h$-cobordisms coincide, so that we can deduce \autoref{mainthm} from \autoref{unstablemainthm}.

\autoref{unstablemainthm} is a lift of several known results at the level of $\pi_0$: the isovariant $s$-cobordism theorem of Browder-Quinn \cite{equivsurgery} and Rothenberg \cite{rothenberg} (see also \cite{luck}), as well as Araki and Kawakubo's result that the equivariant Whitehead torsion of an equivariant $h$-cobordism $M\hookrightarrow W$ is trivial if and only if there exists a $G$-representation $V$ such that the equivariant $h$-cobordism $(W\times D(V); M\times D(V), N\times D(V))$ is trivial \cite{AK}. %Our work lifts these earlier results from $\pi_0$ to the entire space of isovariant, respectively equivariant, $h$-cobordisms.

Though this paper takes place entirely in the smooth or DIFF category, we take a moment to comment on the PL and TOP categories. Nonequivariantly, the $s$-cobordism theorem is known in these categories \cite{pl_h_cobordism,steinberger}. Equivariantly, it seems likely that \autoref{unstablemainthm} and \autoref{mainthm} hold in the PL category, and therefore there should be not only an equivariant $s$-cobordism theorem, but also a parametrized theorem for the higher homotopy groups. Proving this would require technical results that go beyond the scope of this paper.

In the TOP category, the equivariant $s$-cobordism theorem is known to hold \cite{steinberger}, but the resulting $\pi_0$ calculations tell us that a tom Dieck style splitting is impossible. Indeed, let us take a prime $p$, let $C_p $ be the cyclic group of order $p$, and take $M$ to be the unit sphere in the $C_p$-representation $\R^k \oplus \C^n$, with $C_p$ acting trivially on the copies of $\R$ and by multiplication by $e^{2\pi i/p}$ on the copies of $\C$. Then it follows from \cite{quinn_jams} that $\pi_0\mc H_G^{\mathrm{top}}(M) = 0$, whereas a splitting as above would give the group $\pi_0\mc H(S^{k-1}) \oplus \pi_0\mc H(M_{hC_p}) \cong 0 \oplus \Wh(C_p) \cong \Z^{\frac{p-3}{2}}$. Despite this failure, it is still possible that the corresponding theorem relating the higher homotopy groups to algebraic $K$-theory still holds.

 \subsection{Overview} In \autoref{sec:preliminaries} we review technical results about $G$-manifolds with corners which will be needed throughout the paper. As in \cite{GIMM1}, a $G$-action on a manifolds with corners is required to satisfy a condition that means that the $G$-manifold can be modeled locally by $G \times_H V \times [0,\infty)^k$, for varying $H \leq G$ and $H$-representations $V$. In \autoref{Gmanifolds} and \autoref{collarsetc} we review and slightly extend definitions and results from \cite{GIMM1} about manifolds with corners and collars, and in \autoref{fiberbundles} we extend some standard results about smooth equivariant vector bundles to the setting of $G$-manifolds with corners. In  \autoref{tub} we prove a tubular neighborhood theorem for neat embeddings of $G$-manifolds with corners, since the tubular neighborhood theorem for non-neat embeddings given in \cite{GIMM1} is not sufficient for some later results in this paper. In \autoref{blowups} we recall the notion of blowing up a manifold along a neat submanifold, and we prove existence and uniqueness of blowups even in the presence of corners, extending results from \cite{arone_kankaanrinta}. Lastly, in \autoref{strata}, we establish notation and some properties of fixed points, strata and compactified strata for $G$-manifolds with corners, which allow us to prove statements about $G$-manifolds by inductive arguments.
  
In \autoref{isov} we discuss isovariance, a stronger condition than equivariance, requiring not only fixed point sets but  strata to be preserved. After observing that isovariant smooth extensions from the boundary do not necessarily exist as they do for equivariant maps, we introduce the notion of strict isovariance, which also requires isovariance of the map on tangent bundles. We prove a smooth extension lemma for strictly isovariant maps.  In \autoref{fiberbundles1} we prove some results about the interaction between equivariant fiber bundles and isovariance. All of this is used in the following section to construct isovariant maps by an inductive procedure.
 
 In \autoref{isovhcob} and \autoref{isovhcobspaces} we adapt the earlier definitions and constructions from \cite{GIMM1} from the equivariant to the isovariant setting. %We introduce the gap hypotheses that will ensure that an equivariant $h$-cobordism is isovariant.
 We define the notion of isovariant $h$-cobordism and the space of isovariant $h$-cobordisms, along with the extra structures on cobordisms that will be needed to stabilize with respect to embeddings. To show that these extra structures do not change the homotopy type of the $h$-cobordism space, we prove that the space of strictly isovariant retractions of an isovariant $h$-cobordism to its bottom forms a contractible Kan complex.
 %, which ensures that the extra structure does not change the homotopy type of the space of isovariant $h$-cobordisms. 
After that, \autoref{isovhcobstab} is mostly a straightforward generalization of the stabilization construction from \cite{GIMM1} to the isovariant setting. Building on the main result in \cite{GIMM1}, we obtain a functoriality result for the space of isovariant $h$-cobordisms, which we use later in the proof of the splitting result.

 In \autoref{sec:torsion} we review the geometric contruction of equivariant and isovariant Whitehead torsion and Whitehead groups, mostly following \cite{luck}. We also introduce the notion of a manifold with facets, often called a manifold with faces in the literature, and extend these notions of torsion to that case. In \autoref{isovhcobrel} we introduce the notion of relative $h$-cobordisms, which are essentially cobordisms that are nontrivial on a single side. We map the space of relative $h$-cobordisms to the space of absolute $h$-cobordisms by ``stacking with a basic $V$-block''. In \autoref{HequalsHV} we use isovariant torsion to show that this map is a homotopy equivalence. Lastly, in \autoref{fiberseq}, we use these ideas to produce the fiber sequence that plays the central role in the splitting result.
  
 In \autoref{unstablesplitting} we prove the isovariant splitting result (\autoref{unstablemainthm}). The proof is by induction over isotropy types, by comparing a split fiber sequence with the fiber sequence constructed in \autoref{fiberseq}. In \autoref{stablesplitting} we stabilize the splitting and prove the stable equivariant splitting result (\autoref{mainthm}), and in \autoref{sec:naturality} we prove the naturality statement in \autoref{mainthmspaces}. A key feature used in this section is the functoriality of the isovariant $h$-cobordism space, which is worked out in \cite{GIMM1} and extended to the isovariant setting in \autoref{isovhcobstab}. This was one of our original motivations for providing a rigorous treatment of the functoriality of the $h$-cobordism space.

\subsection{Acknowledgements}
 It is a pleasure to acknowledge contributions  to this project arising from conversations with Daniel \' Alvarez-Gavela, Andrew Blumberg, Jim Davis, Ian Hambleton, Manuel Krannich, Sander Kupers, Wolfgang L{\"u}ck, Oscar Randal-Williams,
and Shmuel Weinberger. We especially thank Weinberger for sharing with us the example which shows that the analogous splitting result does not hold in the topological category.

Igusa was partially supported by Simons Foundation Grant 686616. Malkiewich was partially supported by the NSF grants DMS-2005524, DMS-2052923, and DMS-2506430, the Simons Fellowship [SFI-MPS-SFM-00011827] and the travel support grant from Simons Foundation International [SFI-MPS-TSM-00013425]. Merling was partially supported by NSF grants DMS-1943925, DMS-2052988, and DMS-2506429, and the Simons Fellowship [SFI-MPS-SFM-00011841]. Malkiewich and Merling thank the Max Planck Institute in Bonn, where they were in residence in the fall of 2018, and where the origins of this project can be traced back to.

\noindent
\textbf{AI Declaration.}
	The mathematics of this paper was carried out, and the paper itself was written, illustrated, and proofread, without the use of generative AI or large language models.

% !TEX root = splitting.tex

\section{Preliminaries on $G$-manifolds with corners}\label{sec:preliminaries}

\subsection{$G$-manifolds with corners}\label{Gmanifolds} We begin by briefly recalling some definitions and facts from \cite{GIMM1}. Throughout, $G$ is a fixed finite group. A \ourdefn{manifold} %of dimension $n$ 
is a smooth manifold with corners, and a \ourdefn{$G$-manifold} is a smooth manifold with a smooth $G$-action. The action is required to be \ourdefn{trivial along corners} in the sense of \autoref{Gcornerdef} below. 

The \ourdefn{depth} of a point $x \in M$ is the unique number $k$ such that some chart sends $x$ to the origin in $[0,\infty)^k \times \R^{n-k}$. The \ourdefn{boundary subspace} $\partial M \subset M$ is the boundary of $M$ as a topological manifold, i.e. the set of all points of depth $\geq 1$. The set of \ourdefn{corner points}, points of depth $\geq 2$, is denoted by $\partial^c M$.

Following \cite{joyce}, for any point $x\in M$, the set of \ourdefn{local boundary components} at $x$ is the inverse limit of $\pi_0(U\cap (\partial M\backslash \partial^c M))$ over neighborhoods $U$ of $x$, so that the depth of $x$ is the number of local boundary components at that point. The \ourdefn{smooth boundary} $\tipartial M$ is the set of pairs $(x,b)$ where $x \in M$ and $b$ is a local boundary component at $x$. This inherits a smooth atlas from $M$, so that $\tipartial M$ becomes a smooth manifold with corners \cite[Definition 2.7.]{joyce}. There is a canonical smooth map $i: \tipartial M\to M$ whose image is $\partial M$. 
For instance, the smooth boundary of an $n$-simplex is the disjoint union of the $(n-1)$-simplices that are its faces.

\begin{figure}[h]
	\centering
	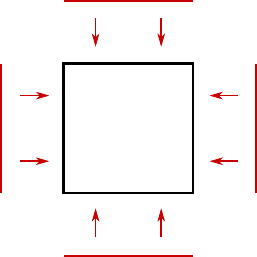
	\vspace{-.5em}
	\caption{The smooth boundary of a square consists of four disjoint line segments.}\label{fig:faces}
\end{figure}

As mentioned above, we consider only $G$-manifolds in which the action of $G$ is trivial along corners.

\begin{defn}\label{Gcornerdef}
A smooth $G$-action on a manifold $M$ with corners is \ourdefn{trivial along corners} (or $M$ is a \emph{$G$-manifold with $G$-trivial corners}) if any of the following equivalent conditions hold:
	\begin{itemize}	
		\item For each boundary point $x \in \partial M$, the action of the stabilizer $G_x$ on the set of local boundary components at $x$ is trivial.
		\item $M$ is locally modeled by $G \times_H V \times [0,\infty)^k$, for varying subgroups $H \leq G$ and $H$-representations $V$.
		\item $M$ is locally modeled by a finite product of smooth $G$-manifolds without corners.
	\end{itemize}
\end{defn}

%\begin{defn}
	Every $G$-manifold without corners has this property. The property is preserved by products.
	
	Throughout this paper, as in \cite{GIMM1}, a \ourdefn{$G$-manifold with corners} is always a manifold $M$ with corners, with a smooth $G$ action that is trivial along corners. 
%\end{defn}

\subsection{Smooth extensions, embeddings, and collars}\label{collarsetc}

In this subsection we recall some technical results from \cite{GIMM1}.

A map defined on an arbitrary subset of $ \R^{n}$ (for example, on an open subset of $[0,\infty)^k \times \R^{n-k}$) is called smooth if it extends to a smooth map defined on an open subset of $\R^n$. The smooth extension exists locally iff it exists globally \cite[Lem 2.12]{GIMM1}. A map defined on a manifold with corners is called smooth if with respect to local coordinates it appears as a smooth map in this sense defined on an open subset of $[0,\infty)^k \times \R^{n-k}$.

An \ourdefn{embedding} is a smooth map that is a topological embedding and whose derivative has full rank (rank equal to the dimension of the domain) at every point. An \ourdefn{open (respectively closed) embedding} is an embedding whose image is open (respectively closed).

\begin{lem}[{\cite[Lem 2.19]{GIMM1}}]\label{depth_preserving_means_open_embedding}
	Suppose that $M$ and $N$ are manifolds (not necessarily compact) of the same dimension, with corners. If $f\colon M \to N$ is a smooth and injective map, and if $f$ has full rank and is depth-preserving, then $f$ is an open embedding.
\end{lem}

\begin{defn}\label{neat_embedding}
An embedding $M\to N$ of manifolds with corners is \ourdefn{neat} if around each depth $k$ point of $M$ it is locally modeled on the inclusion map
\[  \R^{m-k} \times \{0\}^{n-m} \times [0,\infty)^k \to  \R^{m-k} \times \R^{n-m} \times [0,\infty)^k \] 
A \ourdefn{neat submanifold} is a subset of $M$ that is the image of such an embedding.
\end{defn}

In particular, neat embeddings are depth-preserving, and in codimension 0 a neat embedding is the same thing as an open embedding, or an injective map that is locally a diffeomorphism.

\begin{lem}\label{local_model}
If a neat embedding is also equivariant, then it is locally modeled on
\begin{equation*}
 \xymatrix{ G \times_H V \times [0,\infty)^k \ar[r] & G \times_H W \times [0,\infty)^k } 
 \end{equation*}
where $H$ is a subgroup of $G$ and where $V\to W$ is an inclusion of $H$-representations.
\end{lem}

\begin{proof}
Let $x\in M$ be a point of depth $k$ and let $H=G_x$. Since $f$ is an equivariant embedding, $G_{f(x)}=H$. Let $U$ be a small enough $H$-invariant neighborhood of  $f(x)$ such that the $G/H$ coset translates of $U$ are disjoint. Thus it suffices to construct the local model at $x$ $H$-equivariantly.

In local coordinate charts for $M$ and $N$, the embedding is now given near $x$ by an $H$-equivariant embedding $f$ near the origin %, which when composed with the non-equivariant diffeomorphism $\phi$ gives the inclusion of $V$ into $V \oplus V'$:
\[ \begin{tikzcd}
	V\times [0, \infty)^k \rar{f} & W\times [0, \infty)^k % \rar[<->]{\phi}[swap]{\cong} & V \oplus V' \times [0,\infty)^k
\end{tikzcd} \]
of the form $f(v, z)=(f_1(v,z), f_2(v,z))$. 
Since this is non-equivariantly neat, there is a non-equivariant diffeomorphism of the target under which this is the inclusion of $V$ into a linear subspace of $W$. Taking first derivatives, we can deduce that as an $H$-representation we have a decomposition $W \cong V \oplus V'$ along which the matrix $\frac{\partial{f_1}}{{\partial v}}$ is the inclusion of $V$ into $W$. We also deduce that the matrix $\frac{\partial{f_2}}{{\partial z}}$ is diagonal with positive diagonal entries, and the matrix $\frac{\partial{f_2}}{{\partial v}}$ vanishes because any column is a vector in $[0,\infty)^k$ whose negative is also in $[0,\infty)^k$. Therefore the equivariant map
\[ \begin{tikzcd}
	V \times V' \times [0,\infty)^k \rar{F}[swap]{\cong} & W\times [0, \infty)^k
\end{tikzcd} \]
defined by $F(v,t,z) = (f_1(v,z),t,f_2(v,z))$ has first derivative an invertible matrix at the origin. It also preserves depth, and so by the inverse function theorem its inverse is another equivariant smooth function that preserves depth, and is therefore by \autoref{depth_preserving_means_open_embedding} is a local diffeomorphism of manifolds with corners. Composing the second coordinate chart with the inverse of this diffeomorphism gives a new equivariant coordinate chart in which the embedding is just the inclusion
\[ \begin{tikzcd}
	V \times \{0\} \times [0,\infty)^k \subseteq V \times V' \times [0,\infty)^k.
\end{tikzcd} \]
%
%of the form $f(v, z)=(f_1(v,z), f_2(v,z))$. Since the embedding is non-equivariantly neat, the first derivative in the $V$ direction must identify $V$ with a sub-representation of $W$, and so without loss of generality $W \cong V \oplus V'$ and the first derivative of $f_1$ at the origin in the $v$ direction just gives the inclusion $V \subseteq W$. Since 
%
% and up to some non-equivariant diffeomorphism of $W\times [0, \infty)^k$ the embedding is identified with the inclusion $V \subseteq W$. In particular, 
%Without loss of generality $V \subseteq W$ and the first derivative of the embedding at the origin is given by 
\end{proof}

\begin{figure}[h]
	\centering
	\def\svgwidth{4.5in}
	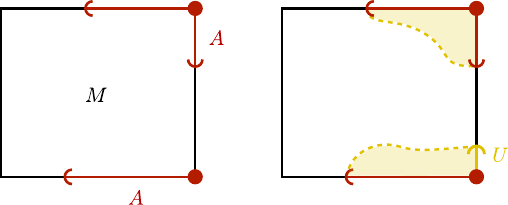
	\vspace{-.5em}
	\caption{A partial boundary $A \subseteq \partial M$ and an open neighborhood of $A$.}\label{fig:partial_boundary}
\end{figure}

%The boundary of a smooth manifold with corners is not a smooth manifold if the corner set is non-empty. Nevertheless, as already discussed we say that a map $\partial M\to N$ is smooth if  (for example, in stating the next result) we may speak of a smooth map $\partial M\to N$, meaning one that with respect to local coordinates appears as a smooth map defined on an open subset of $\partial ([0,\infty)^k \times \R^{n-k})$.
\begin{defn}\label{tame_def}
	A subset $A \subseteq \partial M$ is a \ourdefn{partial boundary} if near each point of $A$, $A$ is a union of local boundary components. Stated more precisely, this means that %$A$ is the image of an open subset $\ti A \subseteq \ti \partial M$ of the smooth boundary.
	for every point $x \in A$ of depth $k$ there is some number $1\le a\le k$ such that in some coordinate chart $U 		\subseteq [0,\infty)^k \times \R^{m-k}$ the set $A\cap U$ corresponds to 
	\[ U \cap \big(\partial([0,\infty)^a) \times [0,\infty)^{k-a} \times \R^{m-k}\big). \]
	Let $\ti A\subseteq \ti \partial M$ be the preimage of $A$. See \autoref{fig:partial_boundary}.
\end{defn}

%%%PICTURE

%\begin{lem}\label{smooth_on_boundary}
%	A map $f\colon \partial M \to N$ is smooth iff $f \circ i\colon \ti\partial M \to N$ is smooth.
%\end{lem}

%\tmar{I prefer to reserve the term ``inclusion'' for the canonical injection from a subset to a set.}
%\begin{lem}\label{smooth_extension}
%	An equivariant map $f\colon \partial M \to \textup{int } N$ is smooth iff it extends to an equivariant smooth map $ U \to N$ for some open subset $U \subseteq M$ containing $\partial M$.
%\end{lem}
%\cmar{Add results involving the partial boundary}

%\begin{prop}\label{smooth_extension}
%	If $A \subseteq \partial M$ is a $G$-invariant partial boundary, $f\colon A \to \textup{int }N$ is an equivariant map, and $f \circ i\colon \ti A \to \textup{int }N$ is smooth, then $f$ extends to an equivariant smooth map $U \to \textup{int }N$ for some $G$-invariant open set $U$ containing $A$.
%\end{prop}

\begin{prop}[{\cite[Cor 2.32]{GIMM1}}]\label{smooth_extension}
	If $S \subseteq M$ is the union of a $G$-invariant partial boundary $A \subseteq \partial M$ and a $G$-invariant neat submanifold $P \subseteq M$, any equivariant map $f\colon S \to \textup{int }N$ whose restriction to $P$ is smooth and whose restriction to each component of $\ti A$ is smooth, must extend to a smooth equivariant map $U \to \textup{int }N$ for some $G$-invariant open set $U$ containing $S$.
\end{prop}

In particular, when $S = \partial M$, a map that is smooth on each component of $\ti\partial M$ extends to a smooth map on a neighborhood of $\partial M$. As observed in \cite[Rem 2.10]{GIMM1}, this does not always hold for maps $f\colon \partial M \to N$ that are allowed to hit the boundary of $N$; it can happen that every smooth extension of $f$ ``shoots outside'' of $N$. This phenomenon only occurs in the presence of corners. %It occurs even non-equivariantly.
Fortunately, for any $N$ there is always a slightly larger manifold containing $N$ in its interior:

%\begin{prop}[{\cite[Prop 2.46]{GIMM1}}]\label{deform_embedding_to_interior}
%	If $N$ is a compact smooth $G$-manifold with corners, there is an isotopy of equivariant embeddings from $\id_N$ to an embedding $N \to N$ whose image is contained in the interior of $N$.
%\end{prop}

\begin{lem}[{\cite[Cor 2.47]{GIMM1}}]\label{extend_to_open_manifold}
	Every compact smooth $G$-manifold with corners can be smoothly equivariantly embedded into the interior of another $G$-manifold of the same dimension.
\end{lem}

Even in the presence of corners, equivariant continuous maps of $G$-manifolds can be approximated by equivariant smooth maps, assuming again that the map to be approximated lands in the interior:

\begin{lem}[{\cite[Cor 1.12]{wasserman},\cite[Lem 2.28]{GIMM1}}]\label{smooth_approximation}
	Any continuous equivariant map $f\colon M \to \textup{int } N$ can be approximated by a smooth equivariant map. If $f$ is smooth on a neighborhood of a closed subset $C \subseteq M$ then the approximation can be rel $C$ (that is, leaving $f$ unchanged in some neighborhood of $C$).
\end{lem}

\begin{defn}\label{df:face}
A \ourdefn{face} $F$ of $M$ is a $G$-invariant subspace of the smooth boundary $\tipartial M$ such that
\begin{itemize}
	\item $F$ is a union of components of $\tipartial M$ and
	\item the map $i\colon \tipartial M \to M$ is injective when restricted to $F$.
\end{itemize}
The set $i(F)\subset \partial M$ may also be referred to as a face of $M$.
\end{defn}

To prevent misunderstandings, we observe a few things about this notion of face. First, there may be points of the boundary that do not belong to any face. For example, the boundary of a teardrop-shaped closed region in the plane is not a face  since its inclusion is not injective. Second, even if $G$ is trivial, a face is not required to be connected. Third, in the case of an $n$-simplex (but not a cube) with trivial action, the non-empty faces in this sense are what are usually called the faces.

%The following is a useful consequence of \autoref{smooth_extension} and \autoref{extend_to_open_manifold}. 
%\begin{cor}\label{smooth_on_partial_delta} 
%%\mmar{This is not immediate from \autoref{smooth_on_boundary} since the target there is $\R$. }
%	An equivariant map $M \times \partial \Delta^k \to N$ that is smooth on $M\times F$ for each face $F$ of $\Delta^n$ can always be extended to a smooth equivariant map $M \times U \to \ti N$, for some $G$-invariant open neighborhood $U$ of $\partial \Delta^k$ in $\Delta^k$ and for some $\ti N$ that contains $N$ in its interior.
%\end{cor}

\begin{defn}\label{df:collar}
	Let $F$ be a face of a $G$-manifold $M$. A \ourdefn{collar} on $F$ is an extension of the inclusion map $F\cong F\times 0 \to M$ to an equivariant embedding $c\colon F \times I \to M$ whose restriction to $F \times [0,1)$ is depth-preserving (and so is an open embedding of that subset). A collar is \ourdefn{neat} if it decreases depth by exactly 1 on $F \times \{1\}$. 
\end{defn}

A \ourdefn{partially defined collar} is an extension of $F \to M$ to a neat equivariant embedding $ c:L\to M$ where $L$ is an open subset of $F \times I$ containing the bottom and sides:
	\[ (F \times \{0\}) \cup (\partial F \times I) \subseteq L \subseteq F \times I \]
A \ourdefn{collar germ} is an equivalence class of partially defined collars with respect to the relation of agreeing on some neighborhood of the bottom and sides. In particular, any collar represents some collar germ. The following lemmas provide motivation for our definition of an $h$-cobordism in \autoref{isovhcob}.
	
\begin{lem}[{\cite[Lem 2.55]{GIMM1}}]\label{collars_unique}\label{collar_isotopy_extension}
	Every face $F$ in a smooth compact $G$-manifold with corners $M$ has a neat collar $F \times I \to M$. Any two collars are isotopic through collars, and any two neat collars are related by an ambient isotopy of $M$. 
\end{lem}

\begin{lem}[{\cite[Lem 2.56]{GIMM1}}]\label{partial_collar}
Every collar germ is represented by some neat collar. If two neat collars represent the same collar germ, then they are related by an (equivariant) isotopy that is fixed in a neighborhood of the bottom and sides. 
\end{lem}

\subsection{Smooth extension for sections of bundles}\label{fiberbundles}

We discuss some results about smooth equivariant fiber bundles, which will be used in building strictly isovariant maps by an inductive process. %Results like this are standard in the continuous and smooth-with-boundary setting -- the point is to carefully check them in our setting of $G$-manifolds with corners.
Recall the following from \cite[Def 2.22]{GIMM1}.

\begin{defn}\label{fiber_bundle}
A \ourdefn{smooth $G$-equivariant fiber bundle} with fiber $F$ is a smooth equivariant map of $G$-manifolds with corners, $p\colon E \to B$, such that for each $x \in B$ with isotropy group $G_x$ there is a $G_x$-invariant open neighborhood $U_x \subseteq B$, an action of $G_x$ on $F$, and a $G_x$-equivariant diffeomorphism $p^{-1}(U_x)\cong F\times U_x$ preserving the product projection to $U_x$.
	
	Now suppose that a group $\Gamma$ is acting smoothly and faithfully on $F$. The bundle has \ourdefn{structure group} $\Gamma$ if there is a privileged class of trivializations as above, covering $B$, such that for each point $x$ the action of $G_x$ on $F$ is given by a homomorphism $G_x \to \Gamma$, and such that for any two of the trivializations the transition map $U_x\cap U_y\to \Diff(F)$ takes values in $\Gamma$. As a special case, a \ourdefn{smooth $G$-equivariant vector bundle} of rank $n$ is a smooth equivariant fiber bundle with fiber $\mathbb R^n$ and structure group $GL_n(\R)$.
\end{defn}

%In particular, a map of $G$-manifolds with corners $E \to M$ is a smooth equivariant vector bundle if it is locally of the form $U \times V \to U$, where $V$ is a $G$-representation, not necessarily an orthogonal one, and the transition maps are diffeomorphisms that are linear on each fiber. 

\begin{rem}
	A smooth $G$-equivariant fiber bundle is always equivariantly locally trivial in the sense that every point $x\in M$ has a $G$-invariant open neighborhood of the form $G \times_H U$, with $H$ being $G_x$ and $U$ being an open $H$-invariant neighborhood of $x$, such that the restriction of the bundle to this set is equivariantly isomorphic to the projection 
$$
G\times_H (F\times U)\to G\times_H U.
$$
	This is true because we may pick $U$ small enough that the sets $gU$ are disjoint for different cosets $gH$. Then the above isomorphism follows immediately from the $H$-equivariant isomorphism $p^{-1}(U) \cong F \times U$.
\end{rem}

%Associated to $E$ are three important but distinct bundles, each of which has fiber $GL(V)$, the group of non-equivariant isomorphisms $V \to V$, with $G_x$ acting by conjugation.

\begin{defn}\label{associated_bundles}
	Let $E,E' \to M$ be two smooth equivariant vector bundles, with transition maps $(u,v) \mapsto (u,T_u(v))$ and $(u,v) \mapsto (u,T_u'(v))$, respectively. The \ourdefn{bundle of fiberwise embeddings} $\Emb(E,E')$ is locally of the form $U \times I(\R^m,\R^n) \to U$, where $I(\R^m,\R^n)$ is the space of injective linear maps $\R^m \to \R^n$. The transition functions are of the form
	\[ (u,A) \mapsto (u,T_u' \cdot A \cdot T_u^{-1}). \]
	This is a smooth $G$-fiber bundle with $G$ acting on $I(\R^m,\R^n)$ through $GL_m(\R) \times GL_n(\R)^\op$. An equivariant section is identified with a fiberwise-linear injective map $E \to E'$.
\end{defn}

As a special case, when $E = E'$ this is the bundle of fiberwise linear automorphisms $\Aut(E) \to M$, which is locally of the form $U \times GL(V) \to U$ and has transition map $(u,A) \mapsto (u,T_u \cdot A \cdot T_u^{-1})$. This is not to be confused with the principal bundle $P \to M$ which is also locally of the form $U \times GL(V) \to U$, but whose transition map is $(u,A) \mapsto (u,T_u \cdot A)$, which gives $P$ a right action by $GL(V)$.
%	\begin{itemize}
%		\item The \ourdefn{principal bundle} $P \to M$ is locally of the form $U \times GL(V) \to U$, and each transition map is of the form $(u,A) \mapsto (u,T_u \cdot A)$. This bundle has a right action by $GL(V)$, and $P \times_{GL(V)} V$ is isomorphic to $E$ again. (All of this is done locally, as the isomorphism type of the $G_x$-representation $V$ could change when moving between different components of the fixed points.)
%		\item The \ourdefn{automorphism bundle} $\Aut(E) \to M$ is also locally of the form $U \times GL(V) \to U$, but the transition map is of the form $(u,A) \mapsto (u,T_u \cdot A \cdot T_u^{-1})$. This bundle does not have a right action by $GL(V)$; rather, sections of this bundle are identified with isomorphisms $E \cong E$ over $M$. An equivariant section is identified with an equivariant isomorphism.
%		\item More generally, for smooth equivariant vector bundles $E$ and $E'$ of possibly different dimensions, 
%	\end{itemize}
%\end{defn}
As in \cite{GMM}, these are bundles in which $G$ acts on the fiber and on the structure group, and in which the action of the structure group on the fiber is $G$-equivariant. %In the first two cases $G$ acts on the fiber $GL(V)$ and structure group $GL(V)$ by conjugation. The only difference is that in the principal bundle, the structure group acts on the fiber by multiplication on the left, while in the automorphism bundle the structure group acts on the fiber by conjugation. In the principal bundle there is also a leftover right action of $GL(V)$ on the fibers, and therefore on the entire bundle, that is also a $G$-equivariant action.

\begin{prop}\label{section_approximation}
	If $E \to M$ is a smooth equivariant fiber bundle, then any equivariant continuous section $s\colon M \to E$ can be approximated by an equivariant smooth section $\tilde s$. If $s$ is smooth on an open neighborhood of some closed set $C_0 \subseteq M$ then $\tilde s$ can be taken to agree with $s$ on $C_0$.
\end{prop}

\begin{proof}
	Each orbit $G/H \subseteq M$ has an open neighborhood on which the bundle is trivialized, and in that neighborhood the approximation rel $C_0$ exists by \autoref{smooth_approximation}. Covering $M$ with finitely many closed neighborhoods of such orbits $C_1$, $\ldots$, $C_k$, we can successively make $s$ smooth on an open neighborhood of the union $\bigcup_{i = 0}^j C_i$ by making it smooth near $C_j$ while not changing its value on the closed set $\bigcup_{i = 0}^{j-1} C_i$.
\end{proof}

\begin{cor}\label{embedding_approximation}
	For any two smooth vector bundles $E$ and $E'$ on $M$, any continuous equivariant fiberwise linear embedding $\varphi\colon E \to E'$ can be approximated by a smooth equivariant fiberwise linear embedding $\tilde\varphi$. If $\varphi$ is smooth on an open neighborhood of some closed set $C_0 \subseteq M$ then $\tilde\varphi$ can be taken to agree with $\varphi$ on $C_0$.
\end{cor}

\begin{proof}
	Apply \autoref{section_approximation} to the bundle of fiberwise embeddings $\Emb(E,E')$ from \autoref{associated_bundles}.
\end{proof}

The proof of the next result uses \ref{section_approximation}.

\begin{prop}\label{section_extension}
	Suppose that $A \subseteq \partial M$ is a closed partial boundary, and that the inclusion $A \to M$ is an equivariant homotopy equivalence, and let $E \to M$ be a smooth equivariant fiber bundle. For any smooth equivariant section $s\colon U \to E$ defined on a neighborhood of $A$ in $M$, there is a smooth global equivariant section $M\to E$ that agrees with $s$ on a smaller neighborhood of $A$.
\end{prop}

\begin{proof}
Because the pair $(M,A)$ admits a $G$-CW structure, $A$ has an equivariant mapping cylinder neighborhood in $M$. This may be chosen to be contained in $U$, and we may take it to be $U$. 
%Since $A$ is a subcomplex of a $G$-CW complex, $A$ has a mapping cylinder neighborhood $U$ that admits an equivariant deformation retraction to $A$.  Therefore there is an equivariant homotopy equivalence $U\to M$. 
Because the inclusions $A\to M$ and $A\to U$ are equivariant homotopy equivalences, the same is true of $U\to M$. Pick a continuous equivariant homotopy inverse $g\colon M \to U$ that is the identity near $A$. Pulling back $E$ along a homotopy $h\colon g \sim \id$ rel $A$, we get a continuous equivariant fiber bundle $h^*E \to M \times I$ that is identified with $E \to M$ along one end and $g^*E \to M$ along the other, and is also identified with $(E \times I)|_{A \times I}$ near $A \times I$. By the usual covering homotopy property of fiber bundles, $h^*E$ is continuously, equivariantly isomorphic to either of the two trivial bundles $E \times I$ or $g^*E \times I$. Furthermore since $A \subseteq M$ is a $G$-invariant subcomplex on which the bundle has already been trivialized in this way, the isomorphism can be chosen to respect that trivialization.

The given smooth section $s\colon U \to E$ gives a continuous section $g^*s\colon M \to g^*E$ that agrees with $s$ near $A$. Since the bundle $h^*E$ is trivial in the $I$ direction, this gives a continuous section $(g^*s \times \id)\colon M \times I \to h^*E$, that agrees with $s \times \id$ near $A \times I$. At the other end of $h^*E$, this gives a continuous section $s'\colon M \to E$ that agrees with $s$ near $A$. Finally, by \autoref{section_approximation}, $s'$ can be approximated by a smooth section $s''$ which still agrees with $s$ near $A$.
\end{proof}

%\begin{cor}\label{isomorphism_extension}
%	Suppose $A \subseteq \partial M$ is a closed partial boundary, $A \to M$ is an equivariant homotopy equivalence, and $E, E' \to M$ are two smooth vector bundles. Any isomorphism $\varphi\colon E|_U \cong E'|_U$ defined on an open set $U \subseteq M$ containing $A$ can be extended to an isomorphism on all of $M$, that agrees with $\varphi$ on a smaller neighborhood of $A$.
%\end{cor}
%
%\begin{proof}
%	As in the proof of \autoref{section_extension}, after shrinking $U$, the inclusion $U \to M$ is an equivariant homotopy equivalence. Since $E \cong E'$ on $U$, by \autoref{bundle_he_isom}, there must exist an isomorphism $E \cong E'$ on $M$. Therefore without loss of generality $E = E'$ and we are extending a self-isomorphism of $E$ from $U$ to $M$. Now we just apply \autoref{section_extension} to the automorphism bundle $\Aut(E) \to M$.
%\end{proof}

\begin{cor}\label{embedding_extension}
	Suppose  $A \subseteq \partial M$ is a closed partial boundary, $A \to M$ is an equivariant homotopy equivalence, and $E, E' \to M$ are two smooth vector bundles. Then any fiberwise linear embedding $\varphi\colon E|_U \to E'|_U$ defined on an open set $U \subseteq M$ containing $A$ can be extended to a fiberwise linear embedding on all of $M$, that agrees with $\varphi$ on a smaller neighborhood of $A$.
\end{cor}

\begin{proof}
We  apply \autoref{section_extension} to obtain a section to the bundle of fiberwise embeddings $\Emb(E,E')$ from \autoref{associated_bundles}, and then we  apply  \autoref{section_approximation} to smoothly approximate it without changing it near $A$.
\end{proof}

\subsection{Neat tubular neighborhoods}\label{tub}

In this section we prove a parametrized tubular neighborhood theorem for neat equivariant submanifolds in the presence of corners. 
The tubular neighborhood theorem from \cite[Thm 2.40]{GIMM1} is insufficient for our purposes, because it applies only to embeddings $M \to N$ that land in the interior of $N$ (whereas neat embeddings never land in the interior unless $\partial M$ is empty). The proof of the old result required that limitation, %here requires something different, because when an embedding touches the boundary of $N$ the smooth extension result from \autoref{smooth_extension} breaks down, 
for reasons discussed in \cite[Rem 2.10]{GIMM1}. The key to the proof of the new result here is \autoref{neat_isotopy} below.

Recall the notion of a smooth equivariant vector bundle from \autoref{fiber_bundle}.
%By a smooth equivariant vector bundle, or a $G$ vector bundle, on a $G$-manifold $M$ we mean a smooth vector bundle $\pi:E\to M$ together with a smooth action of $G$ on $E$ such that the map $\pi$ is equivariant and such that the map $\pi^{-1}(x)\to \pi^{-1}(gx)$ given by an element $g\in G$ is always linear. In other words, a $G$ vector bundle is a vector bundle with an action of $G$ by vector bundle morphisms. 
The tangent bundle $TM$ of a smooth $G$-manifold is an example of a smooth equivariant vector bundle on $M$, as is the normal bundle $\nu_i$ of an equivariant embedding $i\colon M\to N$. By ``normal bundle'' we mean the quotient $(i^\ast TN)/TM$; we do not always assume that an inner product is given.

%We will prove an equivariant tubular neighborhood theorem for neat embeddings $i\colon M \to N$ of corner-trivial $G$-manifolds. %The proof given in \cite{GIMM1} does not work when the submanifold has corners.

Suppose that $M$ is neatly embedded in $N$, and let $\nu\to M$ be the normal bundle. Note that $M$ is also neatly embedded in $\nu$ as the zero section, and that the normal bundle of this embedding is canonically identified with $\nu$. Roughly, a tubular neighborhood is an isomorphism between a neighborhood of $M$ in $N$ and a neighborhood of $M$ in $\nu$ such that the induced isomorphism of normal bundles is ``the identity''. We will be concerned with both closed tubular neighborhoods and germs of tubular neighborhoods.

\begin{defn}\label{neat_tubular_nbhd}
A (closed) \ourdefn{tubular neighborhood} of an equivariant embedding $i\colon M \to N$ is given by a $G$-invariant inner product on the normal bundle $\nu$ together with a smooth equivariant embedding of the closed unit disk bundle
\[ \tilde i:D(\nu)\to N \]
such that the restriction of $\ti i$ to the zero-section is $i$, the image of $\ti i$ is a neighborhood of $i(M)$, and $\ti i$ induces the identity  isomorphism of normal bundles, see \cite[Def 2.23]{GIMM1}. It is a \ourdefn{neat tubular neighborhood} if it preserves depth on the open disc bundle and $\tilde{i}^{-1}(\partial N)\subseteq \partial D(\nu)$, and decreases depth by $1$ on the boundary sphere bundle $S(\nu)$.
\end{defn}

\begin{rem}A neat tubular neighborhood is not actually a neat embedding of $D(\nu)$, because on the subspace $S(\nu)$ it does not preserve depth.\end{rem}

%What really matters is the isomorphism class of the tubular neighborhood, where $\tilde i:D(\nu)\to M$ and $\tilde i':D(\nu')\to M$ are considered isomorphic if they are related by composing with a bundle isomorphism $\nu \cong \nu'$ (that preserves inner products). Of course, this isomorphism is unique if it exists. Note that the embedding $\tilde i$ determines an isomorphism between $\nu$ and the normal bundle of the embedding $i$ (defined as a quotient of tangent bundles). Each isomorphism class has a unique representative in which $\nu$ is the normal bundle of $i$ and the bundle isomorphism is the identity map. We will ordinarily use this representative.

%A path of such is allowed to smoothly vary the inner product, but fix the bundle $\nu$.

We also consider tubular neighborhood germs. These are represented by partially defined tubular neighborhoods.

\begin{defn}\label{tubular neighborhood germ}
A \ourdefn{partial tubular neighborhood} of $i$ is an equivariant smooth embedding $\tilde i:U\to N$ where $U\subseteq \nu$ is a $G$-invariant neighborhood of the zero section. In this case no inner product on $\nu$ is specified. Again the normal bundle isomorphism determined by the embedding is required to be the identity.
\end{defn} 

As in \cite[Def 2.25]{GIMM1}, two such partially defined embeddings are said to determine the same germ if they agree on some neighborhood of the zero section, 

A partial tubular neighborhood $\tilde i$ is \ourdefn{neat} if it preserves depth. By \autoref{depth_preserving_means_open_embedding}, this is equivalent to asking that $\tilde i$ is an open embedding. A germ is \ourdefn{neat} if it has a neat representative.

The following lemma will help us to extend neat tubular neighborhoods, simply by extending them as embeddings and then checking a condition on the coordinates.

\begin{lem}\label{neat_isotopy}
	Suppose that U is an open neighborhood of $0$ in $[0,\infty)^k \times \R^{n-k}$ and that 
	\[ \phi_t\colon U  \to \R^k\times \R^{n-k}. \]
is a smooth family of embeddings parametrized by $t\in T$ for some compact connected manifold $T$. Assume that for every $t\in T$, for every $1 \leq i \leq k$, for every $x = (x_1,\ldots,x_n)\in U$, if $x_i=0$ then $\phi_t(x)_i=0$. If there exists some $t_0\in T$ such that $\phi_{t_0}(U)$ is an open subset of $[0,\infty)^k \times \R^{n-k}$, then there is a smaller open neighborhood $V\subset U$ such that for every $t\in T$ the image $\phi_{t}(V)$ is an open subset of $[0,\infty)^k \times \R^{n-k}$. 

 %If we have an open embedding defined near the origin
	%\[ \phi_0\colon [0,\infty)^k \times \R^{n-k} \to [0,\infty)^k \times \R^{n-k}, \]
	%which extends to an isotopy of smooth embeddings defined near the origin
	%\[ \phi_t\colon [0,\infty)^k \times \R^{n-k} \to \R^n. \]
	%Assume that for any $x = (x_1,\ldots,x_n)$ in the domain, and any $1 \leq i \leq k$, if $x_i = 0$ then $\phi_t(x)_i = 0$ for all $t$. Then $\phi_t$ is also a family of open embeddings $[0,\infty)^k \times \R^{n-k} \to [0,\infty)^k \times \R^{n-k}$ on some smaller neighborhood $U$ of the origin.
	
	%The conclusion also holds if instead of an isotopy, $\phi_t$ is a $T$-family of such embeddings parametrized by a connected compact manifold $T$. (A single choice of open set $U$ works across the entire family.)
\end{lem}

\begin{proof}

%	A diffeomorphism of $[0,\infty)^k \times \R^{n-k}$ must preserve depth, so it sends each coordinate axis of $[0,\infty)^k$ into a subspace of the form $[0,\infty) \times \R^{n-k}$. Furthermore, it must also preserve the subspace of the tangent space corresponding to $\R^{n-k}$ by the argument in \autoref{depth_decreasing}. Therefore at any point $x = (0,\ldots,0,x_{k+1},\ldots,x_n)$ where the first $k$ coordinates are zero, the first derivative of $\phi_0$ must have the block form
%	\[ \begin{pmatrix} PD & 0 \\ A & B \end{pmatrix} \]
%	where $P$ is a $k \times k$ permutation matrix, $D$ is a diagonal matrix with positive diagonal entries, $A$ is a $k \times (n-k)$-matrix and $B$ is an invertible $(n-k) \times (n-k)$ matrix. By the assumptions on the isotopy, the permutation matrix $P$ must be the identity.
%	
%	On the other hand, the condition on the first $k$ coordinates forces the first derivative $D\phi_t$ to have its first $k$ rows also be a diagonal $k \times k$ matrix and a zero $k \times (n-k)$ matrix, so that the above block form is maintained throughout the isotopy at all points of the form $x = (0,\ldots,0,x_{k+1},\ldots,x_n)$. Therefore the first derivative, at the very least, preserves the subspace $[0,\infty)^k \times \R^{n-k}$.
	For convenience, extend $\phi_t(x)$ to be a smooth function defined for $(t,x)$ in a neighborhood of $T\times 0$ in $T\times \R^n$. By the compactness of $T$ we may take this domain to be $T\times W$ for some neighborhood $W$ of $0$. Since $\phi_t(x)_i$ vanishes wherever $x_i=0$, there is a smooth function $\psi_t(x)_i$ having the same domain and satisfying $\phi_t(x)_i=x_i\psi_t(x)_i$. Since $\phi_t$ has full rank and $\phi_t(x)_i$ vanishes where $x_i=0$, the partial derivative $\frac{\partial \phi_t(x)_i}{\partial x_i}$ is not zero wherever $x_i=0$. That is, $\psi_t(0)_i\neq 0$. This number is positive when $t=t_0$ and, because $T$ is connected, $\psi_t(0)_i>0$ for all $t$. We can insure that $\psi_t(x)_i$ is positive for all $(t,x)$ by making $W$ still smaller. Now $\phi_t(x)_i$ has the same sign as $x_i$ throughout the domain.
	
	Using the Inverse Function Theorem and the compactness of $T$, we see that by making $W$ smaller we can arrange for the restriction of $\phi_t$ to $W$ to be an open embedding for all $t$. Because of the statement about signs, the restriction of $\phi_t$ to $W\cap ([0,\infty)^k\times \R^{n-k})$ is an open embedding in $([0,\infty)^k\times \R^{n-k})$. We may take $V$ to be $W\cap U$.
\end{proof}

Next we define a space of neat tubular neighborhoods as in \cite{GIMM1}.

\begin{defn}\label{tubular_nbhd_space}
Suppose $M$ is compact and $i\colon M \to N$ is a fixed neat embedding. Then let $\NTub_k(M)$ be the space of neat tubular neighborhoods: a $k$-simplex is a metric on $\Delta^p \times \nu$ and a $\Delta^p$-family of neat tubular neighborhoods of $i$, represented by a single smooth map $D(\Delta^p \times \nu) \to N$. %whose track $\Delta^p \times D(\nu) \to \Delta^p \times N$ is a smooth embedding.
Note that the inner product on the normal bundle $\nu$ is allowed to vary smoothly over $\Delta^p$.

The simplicial set of neat tubular neighborhood germs $\NTub_{\sbt}^g(M)$ is defined similarly. A $k$-simplex is a map $\Delta^p \times U \to N$ for a fixed open set $U \subseteq \nu$ containing the zero section, that for each point of $\Delta^p$ is a partial tubular neighborhood of $i$, and two such simplices are identified if they agree on $\Delta^p \times U'$ for some smaller open subset $U'$ containing the zero section.
\end{defn}

In the above definition, we take a $\Delta^p$-family of neat embeddings $M \to N$ to be a smooth map $\Delta^p \times M \to N$ whose track $\Delta^p \times M \to \Delta^p \times N$ is locally modeled by the inclusion
\[  [0,\infty)^d \times \R^{p-d} \times \R^{m-k} \times \{0\}^{n-m} \times [0,\infty)^k \to [0,\infty)^d \times \R^{p-d} \times \R^{m-k} \times \R^{n-m} \times [0,\infty)^k \]
where $[0,\infty)^d \times \R^{p-d}$ is a local model for $\Delta^p$, $\R^{m-k} \times [0,\infty)^k$ is a local model for $M$, and $\R^{n-k} \times [0,\infty)^k$ is a local model for $N$. By the same argument as in \autoref{local_model}, an equivariant family of neat embeddings is locally modeled by the inclusion
\[ \xymatrix{ [0,\infty)^d \times \R^{p-d} \times G \times_H V \times [0,\infty)^k \ar[r]^-{\subseteq} & [0,\infty)^d \times \R^{p-d} \times G \times_H W \times [0,\infty)^k, } \]
arising from a single inclusion of $H$-representations $V \subseteq W$. We will use this when we prove the following theorem, although we will only use the theorem in cases where the family of neat embeddings is a constant family.
%as the simplicial set in which a $k$-simplex is a smooth map $\Delta^p \times U \to N$, where $U \subseteq \nu$ is an open set containing the zero section, that at every point of $\Delta^p$ is a partial tubular neighborhood of $M$ in $N$. There is no inner product on $\nu$ in this case, and these are considered up to germs, so we can shrink $U$ as much as we like.

\begin{thm}[Neat Tubular Neighborhood Theorem]\label{neat_tubular_nbhd_thm}
	Assume that $M$ and $N$ are compact $G$-manifolds with corners. Then for any $\Delta^p$-family of neat embeddings $M \to N$, any $\partial \Delta^p$-family of neat tubular neighborhoods $D(\partial \Delta^p \times \nu) \to N$ can be extended to a $\Delta^p$-family of neat tubular neighborhoods $D(\Delta^p \times \nu) \to N$.
	% the spaces $\NTub_{\sbt}(M)$ and $\NTub_{\sbt}^g(M)$ are contractible Kan complexes, and the forgetful map between them $\NTub_{\sbt}(M) \to \NTub_{\sbt}^g(M)$ is an acyclic Kan fibration.
\end{thm}

\begin{cor}
	For any fixed neat embedding $i\colon M \to N$, the spaces $\NTub_{\sbt}(M)$ and $\NTub_{\sbt}^g(M)$ are contractible Kan complexes. %, and the forgetful map between them $\NTub_{\sbt}(M) \to \NTub_{\sbt}^g(M)$ is an acyclic Kan fibration.
\end{cor}

In fact, our proof will show that the forgetful map $\NTub_{\sbt}(M) \to \NTub_{\sbt}^g(M)$ is also an acyclic Kan fibration. Since contractible spaces are nonempty, we get as a corollary the existence of a neat tubular neighborhood for any neat embedding $i\colon M \to N$.

The proof of \autoref{neat_tubular_nbhd_thm} occupies the rest of this section. As in \cite[Theorem 2.40]{GIMM1}, it is easiest to prove this first for the space of neat tubular neighborhood germs. The following is the analogue of \cite[Prop 2.42]{GIMM1}.

%\begin{rem}
%	We don't actually need $M$ to be compact in the proof below, we just need it to be covered by finitely many coordinate charts on which the neat embedding is modeled by $G \times_H V \times [0,\infty)^k \to G \times_H W \times [0,\infty)^k$.
%\end{rem}
\begin{prop}\label{ntn_step1}
	For any $\Delta^p$-family of neat embeddings $M \to N$, any $\partial \Delta^p$-family of neat tubular neighborhood germs $\partial \Delta^p \times U \to N$ can be extended to a $\Delta^p$-family of neat tubular neighborhood germs $\Delta^p \times U' \to N$ for some smaller open set $U'$ containing the zero section of $M \times \nu$.
\end{prop}

\begin{proof}
	The idea is to prove the statement on open subsets $U_\beta \subseteq \Delta^p \times \nu$, and to make $U_\beta$ successively larger by adding one open set $U_\alpha$ at a time, where $U_\alpha$ is contained in a single coordinate chart of $\Delta^p \times \nu$.
	
	We are given a map $\phi_0\colon (\partial \Delta^p \times U) \cup (\Delta^p \times M) \to N$ which for each $t \in \partial \Delta^p$ gives a germ of a neat tubular neighborhood, and for each $t \in \Delta^p$ gives a neat embedding $M \to N$. Let $(t,x,0)$ be any point in the zero section $\Delta^p \times M$ of the bundle $\Delta^p \times \nu$. Let $H$ be the stabilizer group of $x$ and let $(G \times_H U_\alpha) \subseteq \Delta^p \times \nu$ be any $G$-invariant open subset about $(t,x,0)$ whose image in $N$ is contained in a single coordinate chart in which the family of neat embeddings $\Delta^p \times M \to \Delta^p \times N$ is modeled by the inclusion of the subspace
\[ \xymatrix{ [0,\infty)^d \times \R^{p-d} \times G \times_H V \times [0,\infty)^k \ar[r]^-{\subseteq} & [0,\infty)^d \times \R^{p-d} \times G \times_H W \times [0,\infty)^k, } \]
arising from a single inclusion of $H$-representations $V \subseteq W$.%Since we are only considering germs of tubular neighborhoods, without loss of generality all of the tubular neighborhoods that we consider are defined on an open set $U_\alpha$ that is contained in $\nu|_O$, and whose image stays in this chart where $N$ is identified with $G \times_H W \times [0,\infty)^k$.
	
	In this coordinate chart, the open set $U_\alpha$ corresponds to an open set about the origin in $[0,\infty)^d \times \R^{p-d} \times W \times [0,\infty)^k$.  A germ of a neat tubular neighborhood defined near $(t,x,0)$ is given by an $H$-equivariant map
	\[ U_\alpha \ra [0,\infty)^d \times \R^{p-d} \times W \times [0,\infty)^k \]
	that
	\begin{enumerate}
		\item is the identity in the $[0,\infty)^d \times \R^{p-d}$ coordinates (and so defines part of a $\Delta^p$-family of maps),
		\item is the identity along the zero section $U_\alpha \cap ([0,\infty)^d \times \R^{p-d} \times V \times [0,\infty)^k)$,
		\item along the zero section has first derivative that also induces the identity on the normal direction $W - V$,
		\item is an embedding, and
		\item preserves depth (and therefore by \autoref{depth_preserving_means_open_embedding} is an open embedding).
	\end{enumerate}
%	Two of these are considered to be the same if they agree after shrinking the neighborhood $D_\alpha$.
	We call a map with these properties a \ourdefn{local neat tubular neighborhood germ}.
	
	In this language, $\phi_0$ defines a local neat tubular neighborhood germ that is only defined on the zero section
	\[ [0,\infty)^d \times \R^{p-d} \times V \times [0,\infty)^k \]
	and the sides of $\Delta^p$
	\[ \partial([0,\infty)^d) \times \R^{p-d} \times W \times [0,\infty)^k. \]
	We extend this to a smooth map
	\[ \phi_\alpha\colon U_\alpha \to [0,\infty)^d \times \R^{p-d} \times W \times [0,\infty)^k \]
	that agrees with $\phi_0$ on the above two subspaces, by the smooth extension from a corner and a neat submanifold as defined in \cite[Cor 2.32]{GIMM1}, referenced above as \autoref{smooth_extension}. On each point $s \in [0,\infty)^d \times \R^p$, $\phi_\alpha$ is obtained by an affine-linear combination of maps $[0,\infty)^d \times \R^{p-d} \times W \times [0,\infty)^k \to [0,\infty)^d \times \R^{p-d} \times W \times [0,\infty)^k$ obtained by projecting one of the $[0,\infty)^d$ directions to 0, or projecting $W$ to $V$, and then applying $\phi_0$. 
	
	These affine-linear combinations satisfy points (1) and (2) in the above list automatically. For (3) we notice that each of the projected maps has first derivative in the $W-V$ direction equal to either the identity or zero, and the projections in which it is the identity also have coefficients summing to 1 (they correspond to the smooth extension just from the corner $[0,\infty)^d$), so the final map is also the identity in the $W-V$ direction. For (4), notice that by points (2) and (3), the first derivative of $\phi_\alpha$ is full-rank along the zero section, so after shrinking the open set $U_\alpha$ to a smaller open set about the origin, we can guarantee that $\phi_\alpha$ is an embedding, by the inverse function theorem.
	
	To prove (5), that depth is preserved, we use \autoref{neat_isotopy}. We know that at the origin in $[0,\infty)^d \times \R^{p-d}$ the resulting map $W \times [0,\infty)^k \to W \times [0,\infty)^k$ preserves depth. Furthermore, for every point in $[0,\infty)^d \times \R^{p-d}$, the resulting map $W \times [0,\infty)^k \to W \times [0,\infty)^k$ preserves the subspace in which the $i$th coordinate in $[0,\infty)^k$ is zero, for each $0 \leq i \leq k$, because this subspace is preserved by the affine-linear combination. Therefore by \autoref{neat_isotopy}, if we shrink $U_\alpha$ to be a compact connected neighborhood of the origin in $[0,\infty)^d \times \R^{p-d}$ times a sufficiently small neighborhood of the origin in $W \times [0,\infty)^k$, $\phi_\alpha$ preserves depth on $U_\alpha$ and is therefore an open embedding there.
	
	In summary, after shrinking $U_\alpha$ to be a smaller open set about the origin, we have shown that $\phi_0$ extends to a local neat tubular neighborhood germ $\phi_\alpha$ on $U_\alpha$.
	
	Now that we have done the local extension, we show how to patch them together. Let $U_\beta \subseteq \Delta^p \times \nu$ be any $G$-invariant open set on which we have defined a local neat tubular neighborhood germ $\phi_\beta$ extending $\phi_0$, and let $C_\beta \subseteq U_\beta \cap (\Delta^p \times M)$ be any $G$-invariant compact subset that is contained in the zero section. Similarly let $U_\alpha$ be another $G$-invariant open set whose image in $\Delta^p \times N$ is contained in a single coordinate chart as above, on which we have defined another local neat tubular neighborhood germ $\phi_\alpha$ extending $\phi_0$, and let $C_\alpha \subseteq U_\alpha \cap (\Delta^p \times M)$ be another $G$-invariant compact subset that is contained in the zero section. Since $C_\alpha$ and $U_\alpha^c$ are disjoint closed subsets of $\Delta^p \times \nu$, we can pick an equivariant smooth bump function $\lambda\colon \Delta^p \times \nu \to [0,1]$ that is 1 on a neighborhood of $C_\alpha$ and 0 on a neighborhood of $U_\alpha^c$.
	
	One might expect that $\phi_\beta$ defines a local neat tubular neighborhood germ on the intersection $U_\alpha \cap U_\beta$, but that is not quite true. Some points in this set may be sent outside the given coordinate chart in $\Delta^p \times N$, so they don't match the definition of local neat tubular neighborhood germ lying inside this coordinate chart. However, the compact set $C_\beta \cap (\textup{supp }\lambda)$ is contained inside this coordinate chart, and so it has some open neighborhood $P \subseteq \Delta^p \times M$ on which $\phi_\beta$ defines a local neat tubular neighborhood germ in the above sense.
	
	Now we define a local neat tubular neighborhood germ $\phi_{\alpha\beta}$ on an open set $U_{\alpha\beta}$ containing $C_\alpha \cup C_\beta$, extending $\phi_0$, as follows. Consider the formula
	\[ \phi_{\alpha\beta} = \lambda \phi_\alpha + (1-\lambda) \phi_\beta, \]
	carried out in the coordinate chart $U_\alpha \subseteq [0,\infty)^d \times \R^{p-d} \times G \times_H W \times [0,\infty)^k$. The formula is well-defined
	\begin{itemize}
		\item in a neighborhood of $C_\alpha$ since $\lambda = 1$ there,
		\item in $(\textup{supp }\lambda)^c \cap U_\beta \supseteq (\textup{supp }\lambda)^c \cap C_\beta$ since $\lambda = 0$ there, and
		\item in the neighborhood $P$ of $(\textup{supp }\lambda) \cap C_\beta$ since both $\phi_\alpha$ and $\phi_\beta$ are defined there.
	\end{itemize}
	It is therefore well-defined on an open subset of $U_\alpha \cup U_\beta$ that contains $C_\alpha \cup C_\beta$. By intersecting the images of this open set under the $G$-action, we can enforce that the domain of definition is $G$-invariant as well. The function is also smooth on this domain, since it is clearly smooth where it is defined in $U_\alpha$, and agrees with $\phi_\beta$ on the open set $(\textup{supp }\lambda)^c \cap U_\beta$. Inside $U_\alpha$, as in the first part of the proof, the interpolation satisfies conditions (1)-(3) for a local neat tubular neighborhood germ that extends $\phi_0$, and also satisfies (5), because the coefficients of the linear combination are nonnegative and so the depth is automatically preserved. Only (4) might fail: $\phi_{\alpha\beta}$ might not be an embedding. However since its first derivative is full rank along the zero section, $\phi_{\alpha\beta}$ is an embedding on some smaller open set $U_{\alpha\beta} \subseteq U_\alpha \cup U_\beta$ that contains $C_\alpha \cup C_\beta$, using the version of the inverse function theorem recalled in \cite[Lem 2.18]{GIMM1}.
	
	Since $\Delta^p \times M$ is compact, we can find finitely many $G$-invariant open sets $U_i \subseteq \Delta^p \times \nu$, each of which lands in a single coordinate chart as above, whose image under the neat embedding is in a single coordinate chart in $\Delta^p \times N$, each containing a $G$-invariant closed ball $C_i$ in the zero section $U_i \cap (\Delta^p \times M)$, such that the interiors of the sets $C_i$ cover the zero section $\Delta^p \times M$.
	
	Then we perform the above argument inductively, finding $G$-invariant open subsets $U_{\{1,2,\ldots,k\}} \subseteq \Delta^p \times \nu$ such that
	\[ C_1 \cup \cdots \cup C_k \subseteq U_{\{1,2,\ldots,k\}} \subseteq U_1 \cup \cdots \cup U_k, \]
	and a local neat tubular neighborhood germ $\phi_k$ defined on $U_{\{1,2,\ldots,k\}}$ that extends $\phi_0$. For the inductive step we define
	\begin{align*}
		U_\beta = U_{\{1,2,\ldots,k\}}, & \qquad U_\alpha = U_{k+1}, \\
		C_\beta = C_1 \cup \cdots \cup C_k, & \qquad C_\alpha = C_{k+1}, \\
	\end{align*}
	and the above argument produces the next open set $U_{\{1,2,\ldots,k,k+1\}}$ and the next extension $\phi_{k+1}$. After finitely many stages of the argument, we have the desired neat tubular neighborhood germ $\phi$, defined on an open set of $\Delta^p \times \nu$ containing all of the $C_i$ and therefore containing the zero section. This must contain an open set of the form $\Delta^p \times U$ by the tube lemma, and the proof is complete.
\end{proof}

The following is the analogue of \cite[Prop 2.43]{GIMM1}.

\begin{prop}\label{ntn_step2}
	For any $\Delta^p$-family of neat embeddings $M \to N$, any $\partial \Delta^p$-family of neat tubular neighborhoods $D(\partial \Delta^p \times \nu) \to N$ and compatible $\Delta^p$-family of neat tubular neighborhood germs $\Delta^p \times U \to N$ can be extended to a $\Delta^p$-family of neat tubular neighborhoods $D(\Delta^p \times \nu) \to N$.
\end{prop}

\begin{proof}
	The idea is to work locally in $\Delta^p \times M$ as before, but now we focus on doing local extensions at points $(t,x,v) \in \Delta^p \times \nu$ where $t \in \partial \Delta^p$, instead of focusing on the zero section where $v = 0$.
	
	We are given a map $\phi_0\colon D(\partial \Delta^p \times \nu) \cup (\Delta^p \times U) \to N$, where $U \subseteq \nu$ is an open set containing the zero section, which for each $t \in \partial \Delta^p$ gives a neat tubular neighborhood, and for each $t \in \Delta^p$ gives a germ of a neat tubular neighborhood. Let $(t,x,v)$ be any point in $D(\partial \Delta^p \times \nu) \setminus (\partial \Delta^p \times U)$; note that $v \neq 0$. Let $H$ be the stabilizer group of $x$ and let $(G \times_H U_\alpha) \subseteq \Delta^p \times \nu$ be any $G$-invariant open subset about $(t,x,v)$ such that the image of $U_\alpha \cap D(\partial \Delta^p \times \nu)$ along $\phi_0$ is contained in a single coordinate chart. Near this point, the embedding $D(\partial\Delta^p \times \nu) \to \Delta^p \times N$ is modeled by the inclusion of the subspace
\[ \xymatrix{ \partial([0,\infty)^d) \times \R^{p-d} \times G \times_H W \times [0,\infty)^k \ar[r]^-{\subseteq} & [0,\infty)^d \times \R^{p-d} \times G \times_H W \times [0,\infty)^k. } \]

	So as before, we define a \ourdefn{local neat tubular neighborhood germ} near this point to be an $H$-equivariant map
	\[ U_\alpha \ra [0,\infty)^d \times \R^{p-d} \times W \times [0,\infty)^k \]
	that
	\begin{enumerate}
		\item is the identity in the $[0,\infty)^d \times \R^{p-d}$ coordinates,
		\item is an embedding, and
		\item preserves depth (and therefore by \autoref{depth_preserving_means_open_embedding} is an open embedding).
	\end{enumerate}
	
	We extend $\phi_0$ to a smooth map
	\[ \phi_\alpha\colon U_\alpha \to [0,\infty)^d \times \R^{p-d} \times W \times [0,\infty)^k \]
	that agrees with $\phi_0$ along $\partial([0,\infty)^d) \times \R^{p-d} \times W \times [0,\infty)^k$, by the smooth extension from a corner as before. % as defined in \cite[Cor 2.32]{GIMM1}, referenced above as \autoref{smooth_extension}.
	On each point $s \in [0,\infty)^d \times \R^p$, $\phi_\alpha$ is obtained by an affine-linear combination of maps $[0,\infty)^d \times \R^{p-d} \times W \times [0,\infty)^k \to [0,\infty)^d \times \R^{p-d} \times W \times [0,\infty)^k$ obtained by projecting one of the $[0,\infty)^d$ directions to 0. This satisfies (1)-(3) after shrinking $U_\alpha$ further, by the same argument as in the previous proof.
	
%	These affine-linear combinations satisfy point (1) automatically. This implies the map is full-rank at the origin (the first derivative can be calculated by traveling only in $\partial([0,\infty)^d) \times \R^{p-d} \times W \times [0,\infty)^k \to [0,\infty)^d$), so it satisfies and (2) by the inverse function theorem.
%	
%	As before, we prove (3) by noticing that $\phi_\alpha$ preserves the subspace in which the $i$th coordinate in $[0,\infty)^k$ is zero, for each $0 \leq i \leq k$, because this subspace is preserved by the affine-linear combination. Therefore by \autoref{neat_isotopy}, if we shrink $U_\alpha$ to be a compact connected neighborhood of the origin in $[0,\infty)^d \times \R^{p-d}$ times a sufficiently small neighborhood of the origin in $W \times [0,\infty)^k$, $\phi_\alpha$ preserves depth on $U_\alpha$ and is therefore an open embedding there.
%	
%	In summary, after shrinking $U_\alpha$ to be a smaller open set about the origin, we have shown that $\phi_0$ extends to a local neat tubular neighborhood germ $\phi_\alpha$ on $U_\alpha$.
	
	Now we patch these together by the same argument, starting with all of $\Delta^p \times U$, and adding open sets $U_\alpha$ of the above form one at a time. We take a single compact neighborhood $C \subseteq U$ of the zero section of $\nu$, and put a compact set $C_\alpha$ inside each $U_\alpha$, so that $\Delta^p \times C$ and the sets $C_\alpha$ cover all of $D(\partial\Delta^p \times \nu)$. At the end we get an extension of $\phi_0$ to a neighborhood of the zero section $\Delta^p \times M$ and sides $D(\partial\Delta^p \times \nu)$ in $D(\Delta^p \times \nu)$, that is an embedding and depth-preserving as a map out of $\Delta^p \times \nu$. To modify this to a full tubular neighborhood without changing it on the sides or near the zero section, we use exactly the same technique as in the last step of \cite[Prop 2.43]{GIMM1}, as illustrated in \cite[Fig 8]{GIMM1}.
\end{proof}

\begin{proof}[Proof of \autoref{neat_tubular_nbhd_thm}]
	Given any $\Delta^p$-family of neat embeddings $M \to N$ and any $\partial \Delta^p$-family of neat tubular neighborhoods $D(\partial \Delta^p \times \nu) \to N$, we regard them as a family of germs and use \autoref{ntn_step1} to exend them to a $\Delta^p$-family of germs $\Delta^p \times U \to N$. Then we use \autoref{ntn_step2} to extend the original family of tubular neighborhoods to a $\Delta^p$-family of neat tubular neighborhoods $D(\Delta^p \times \nu) \to N$.
\end{proof}

\subsection{Spherical blowups and blowdowns}\label{blowups}

%\cmar{I think it makes sense thematically to put a short section about blowups here. Because we've just talked about neat tubular neighborhoods, and the next section defines compactified strata.}

We will be using the spherical blowup of a manifold along a submanifold. The basic example is the blowup of $\R^m$ at the origin, which is the smooth map
\begin{align*}
	S^{m-1} \times [0,\infty)  &\to \R^m  \\
	(\sigma,t) &\mapsto t\sigma.
\end{align*}
More generally, there is
\begin{align*}
	S^{m-1} \times [0,\infty) \times \R^{n} \times [0,\infty)^k &\xrightarrow{\pi} \R^m \times \R^{n} \times [0,\infty)^k \\
	(\sigma,t,v,w) &\mapsto (t\sigma,v,w).
\end{align*}
The preimage of a point whose $\R^m$-coordinate is zero is a sphere $S^{m-1}$. These are the points at which $t = 0$.

\begin{defn}\label{spherical_blowup}
A smooth map $f:M\to N$ of manifolds with corners is a \ourdefn{spherical blowup}, or \ourdefn{spherical blowdown}, if $N$ can be covered by open sets $V$ such that $V$ is diffeomorphic to an open set $U \subseteq \R^m \times \R^{n} \times [0,\infty)^k$ and the restricted map $f:f^{-1}(V)\to V$ is diffeomorphic to the restricted map $\pi: \pi^{-1}(U)\to U$. The set of points in $M$ corresponding to points where $t=0$ is called the \ourdefn{exceptional locus}. Its image in $N$, the \ourdefn{exceptional image}, is a neat submanifold $P \subseteq N$. We say that $M$ is the \ourdefn{spherical blowup of $N$ along $P$}, and write $M = B(N,P)$.
\end{defn}

\begin{rem}
	We use both ``blowup'' and ``blowdown'' to refer to the same notion, preferring the term ``blowup'' when the emphasis is on the passage from $N$ to $M$, and ``blowdown'' when the emphasis is on the passage from $M$ to $N$.
\end{rem}

Away from the exceptional locus the map $f$ is a diffeomorphism
\[ M \setminus f^{-1}(P) \congar N \setminus P, \]
while the map $f^{-1}(P) \to P$ is a fiber bundle with fiber $S^{m-1}$. In fact, it is isomorphic to $S(\nu(P \to N))$, the sphere bundle of the normal bundle of $P$ in $N$. Conceptually, passing from $N$ to its spherical blowup $B(N,P)$ means removing $P$ and then adding a point for every normal direction at every point of $P$. Note that the depth in $M$ of each point in $f^{-1}(P)$ is one higher than the depth in $N$ of the image point. 

\begin{rem}
	A spherical blowdown with $m \geq 2$ is \emph{not} a submersion of manifolds with corners. For example, in the simplest example $S^1 \times [0,\infty) \to \R^2$ the first derivative fails to have rank 2 at points in $S^1 \times \{0\}$.
\end{rem}

\begin{rem}
	We say ``spherical'' to distinguish this notion from the projective blowup, in which the exceptional locus forms a fiber bundle over $P$ with fiber $\mathbb{RP}^{m-1}$. The projective version corresponds to the usual notion of blowup in algebraic geometry.
\end{rem}

\begin{lem}
	Given any neat submanifold $P \subseteq N$, there is a spherical blowup $M \to N$ along $P$, and this is unique up to unique diffeomorphism over $N$.
\end{lem}

\begin{proof}
	The blowup exists everywhere locally along $N$, and is clearly unique away from $P$. So it suffices to show that any two blowups near the same point of $P$ are canonically diffeomorphic. This is proven in \cite[Lemma 2.5]{arone_kankaanrinta} in the case without corners. The case with corners follows immediately because the lifted map respects the subspace $[0,\infty)^k \subseteq \R^k$.
%	
%	Given a diffeomorphism $\phi$ defined on an open subset of $\R^m \times \R^{n} \times [0,\infty)^k$ and preserving the subspace $\{0\} \times \R^n \times [0,\infty)^k$, if we let $\phi_1$ be the $\R^m$-component and $\phi_2$ be the $\R^n \times [0,\infty)^k$-component,
%	
%	The proof is in Tom's notes.
\end{proof}

\begin{cor}\label{blowup_g_action}
	If $N$ is a smooth $G$-manifold with corners and $P \subseteq N$ is a neat $G$-invariant submanifold, then there is a unique smooth $G$-action on the blowup $M = B(N,P)$ such that the map $M \to N$ is equivariant.
\end{cor}

Suppose that $T$ is an equivariant neat tubular neighborhood of $P$, $U \subseteq T$ is the corresponding open tubular neighborhood whose fibers are open discs, and $S \subseteq T$ is the corresponding sphere bundle, so that $T = U \amalg S$ as sets.

\begin{lem}\label{blowup_removes_tubular_neighborhood}
	With the above notation, we have equivariant diffeomorphisms
	\[ B(N,P) \cong N \setminus U, \qquad B(N,P) \cup_{f^{-1}(P)} T \cong N, \]
	where the union is along the identification $f^{-1}(P) \cong S(\nu(P \to N)) \cong S \subseteq T$. The latter diffeomorphism can also be taken to be the identity on $P \subseteq T$.
\end{lem}

\begin{proof}
	Pick an outward-pointing vector field along $S$, and use this to extend $T$ to a slightly larger neat tubular neighborhood $T'$, such that the inclusion $T \subseteq T'$ is identified with an inclusion of two disc bundles of different radii in the same vector bundle. Without loss of generality, $T$ is identified with the disc of radius 1.
	
	The local trivializations of $T'$ identify open sets in $N$ with $\R^m \times \R^{n} \times [0,\infty)^k$ in such a way that $T$ corresponds to the open ball of radius 1 about the origin in $\R^m$, times all of $\R^{n} \times [0,\infty)^k$. In each of these coordinate charts, we map in $S^{m-1} \times [0,\infty) \times \R^{n} \times [0,\infty)^k$ by the formula
	\[ (\sigma,t,v,w) \mapsto ((t+1)\sigma,v,w). \]
	This is invariant under changes between the trivializations, since these changes apply elements of $O(m)$ in the $\R^m$ direction, so they preserve this formula. This gives the diffeomorphism $B(N,P) \to N \setminus U$. Note that this diffeomorphism sends $f^{-1}(P)$ to $S$ by the canonical identification of both with $S(\nu(P \to N))$, as required.
	
	The second diffeomorphism follows from the fact that $(N \setminus U) \cup_{S} T \cong N$. (The smooth structure on the left-hand side relies on choices of collars, so it is well-defined up to isotopy by \autoref{collar_isotopy_extension}, which is enough to get the statement that a diffeomorphism of the required form exists.)
\end{proof}

In summary, a blowup of $N$ along a neat submanifold $P$ is diffeomorphic to $N$ minus an open neat tubular neighborhood of $P$.

\begin{cor}\label{neat_complement_well_defined}
	For a given neat embedding $i\colon M \to N$ of compact $G$-manifolds with corners, the complement of a neat tubular neighborhood of $i$ is well-defined up to equivariant diffeomorphism.
\end{cor}

Note that this could also be deduced from the neat tubular neighborhood theorem (\autoref{neat_tubular_nbhd_thm}) by taking an isotopy between any two neat tubular neighborhoods and extending it to an ambient isotopy, as in the proof of \cite[Lem 2.55]{GIMM1}.
%\begin{proof}
%	Regard $D(\nu)$ as a submanifold of $N$ using one fixed neat tubular neighborhood. Any two neat tubular neighborhoods are equivariantly isotopic by \autoref{neat_tubular_nbhd_thm}. This isotopy $[0,1] \times D(\nu) \to [0,1] \times N$ can be extended to an equivariant ambient isotopy $[0,1] \times N \to [0,1] \times N$ by the proof of \cite[Lem 2.55]{GIMM1}, and this gives the diffeomorphism between their complements.
%\end{proof}

Recall that a map $g\colon L \to N$ is transverse to $P \subseteq N$ if for every $x \in L$ such that $g(x) \in P$, the tangent spaces $Dg(T_x L)$ and $T_{g(x)}P$ span $T_{g(x)}N$, or in other words if the projection of $Dg(T_x L)$ onto the normal space $\nu_{g(x)}(P \to N)$ is surjective.

\begin{lem}\label{pullback_blowup}
	If $f\colon M \to N$ is a spherical blowup along $P \subseteq N$, and $g\colon L \to N$ is smooth map transverse to $P$, then the pullback $L \times_N M$ is a smooth submanifold of $L \times M$, and it is a spherical blowup of $N$ along $g^{-1}(P) \subseteq L$.
	\[ \xymatrix{
		L \times_N M \ar[d] \ar[r] & M \ar[d]^-f \\
		L \ar[r]^-g & N
	} \]
\end{lem}

As in \autoref{blowup_g_action}, the equivariant version of this statement also follows immediately. The statement about being a submanifold of $L \times M$ allows us recognize when a given smooth structure on $L \times_N M$ agrees with the one that makes $L \times_N M$ into a blowup: we only have to check that the maps to $L$ and $M$ are smooth.

\begin{proof}
%	Much of this follows from the main result of \cite{arone_kankaanrinta}, but we also give an explicit proof here. Note that spherical blowups are not submersions, so we can't just quote away the fact that $L \times_N M$ is smooth.
	It suffices to prove this locally. Away from $g^{-1}(P)$ we are pulling back a diffeomorphism and therefore we get a unique smooth structure making the map to $L$ a diffeomorphism. Along $g^{-1}(P)$, the map $g$ can be locally modelled by
\begin{align*}
	\R^m \times \R^{n'} \times [0,\infty)^{k'} &\to \R^m \times \R^{n} \times [0,\infty)^k \\
	(x,y,z) &\mapsto (x,f_1(x,y,z),f_2(x,y,z)),
\end{align*}
by the rank theorem for submersions on the resulting map to $\R^m$. Therefore the pullback diagram has the local model
	\[ \xymatrix{
		S^{m-1} \times [0,\infty) \times \R^{n'} \times [0,\infty)^{k'} \ar[d] \ar[r] & S^{m-1} \times [0,\infty) \times \R^{n} \times [0,\infty)^k \ar[d] \\
		\R^m \times \R^{n'} \times [0,\infty)^{k'} \ar[r] & \R^m \times \R^{n} \times [0,\infty)^k
	} \]
	\[ \xymatrix{
		(\sigma,t,y,z) \ar[d] \ar[r] & (\sigma,t,f_1(t\sigma,y,z),f_2(t\sigma,y,z)) \ar[d] \\
		(t\sigma,y,z) \ar[r] & (t\sigma,f_1(t\sigma,y,z),f_2(t\sigma,y,z))
	} \]
	Therefore there is a smooth structure on the pullback (the upper-left corner of the diagram) under which the left-hand  vertical map is a spherical blowdown. On the other hand, under this smooth structure the map to the product of the bottom-left and top-right terms is a smooth embedding, split by the maps that project onto the factors $\R^{n'} \times [0,\infty)^{k'}$ and $S^{m-1} \times [0,\infty)$. Therefore this smooth structure on the pullback $L \times_N M$ is determined by its embedding into $L \times M$.
\end{proof}

\subsection{Fixed points, strata, and compactified strata}\label{strata}

In this section we establish notation and terminology for various aspects of the isotopy stratification of a $G$-manifold $M$ with trivial corners. We begin with some terms that apply to $G$-spaces in general.

\begin{defn}\label{fixed_points}
Let $X$ be any space with continuous $G$-action, and let $H \leq G$ be a subgroup of $G$.
	\begin{itemize}
		\item $X^H \subseteq X$ is the set of \emph{$H$-fixed points}, points that are fixed by every element of $H$. \vspace{0.2em}
		\item $X^{(H)} = \bigcup_{H' \sim H} X^{H'}$ is the set of \emph{$H$-conjugate fixed points}, the union of the $H'$-fixed point set as $H'$ runs over all conjugates of $H$ in $G$. \vspace{0.2em}
		\item $X_H \subseteq X^H$ is the \emph{$H$-stratum}, the set of points whose stabilizer is exactly $H$. It can be expressed as the set difference%\kmar{looks like this should be union over $K>H$}\cmar{Fixed!}
	\[ X_H = X^H \Large \setminus \bigcup_{K > H} X^K. \]
		\item $X_{(H)} = \coprod_{H'\sim H} X_{H'}$ is the \emph{$H$-conjugate stratum}, the set of points whose stabilizer is conjugate to $H$.
	\end{itemize}
\end{defn}

Note that the normalizer $NH$ preserves $X^H$, and therefore the Weyl group $WH = NH/H$ acts on $X^H$ and on $X_H$. The full group $G$ preserves $X^{(H)}$ and $X_{(H)}$.

In the special case when $X^K$ is empty for every subgroup $K$ that properly contains $H$ (and therefore also for every subgroup that properly contains a conjugate of $H$), the set $X^{(H)}$ is the \emph{disjoint} union $\coprod_{H'\sim H} X^{H'}$ over conjugates of $H$. %This is so because a point that is fixed by two different conjugates of $H$ is fixed by the subgroup generated by these.

%Note that any two conjugates $X^{H_i}$ and $X^{H_j}$ will have intersection of the form $X^K$, where $K \leq G$ is the subgroup generated by $H_i$ and $H_j$. Such a $K$ must contain $H$ up to conjugation. So if the fixed point subspaces $X^K$ are empty for all $K \geq H$, the spaces $X^{H_i}$ and $X^{H_j}$ for distinct conjugates $H_i$ and $H_j$ are disjoint, and we have that  %Note that each of the subsets $X^{H_i}$ is closed as a subset of $X^{(H)}$, but there are finitely many of them, so each of them is open as well. Therefore this is a disjoint union in the topological sense.

\begin{lem}
Suppose that $M$ is a $G$-manifold with $G$-trivial corners and $H \leq G$. Then
	\begin{itemize}
		\item The subset $M^H$ is a neat $NH$-submanifold of $M$.
		\item If $M^K = \emptyset$ for all $K > H$, then $M^{(H)}$ is a neat $G$-submanifold of $M$.
	\end{itemize}
	%\kmar{But $M^H$ is not a $G$-manifold. This is a $WH$-neat embedding.}\cmar{Fixed!} 
%	\mmar{Hmm, how does the union in the second statement become a $G$-manifold? Each piece is a $WH_i$ manifold }
\end{lem}

\begin{proof}
	For the first statement, since it is a local statement it is enough to consider the case when $M$ is a $G$-invariant open subset of $G \times_K V \times [0,\infty)^k$ for some $K\le G$ and some $K$-representation $V$. Since $(e,0,0)$ is an $H$-fixed point, $K$ must contain $H$. A neighborhood of this point is contained in $V \times [0,\infty)^k$, and on this subspace the inclusion of the $H$-fixed points is identified with
	\[ V^H \times [0,\infty)^k \subseteq V \times [0,\infty)^k \]
which is a neat embedding.
The second statement follows from the first since $M^{(H)} = \coprod_{H\sim H'} M^{H'}$ is a disjoint union of submanifolds.
\end{proof}

%We apply \autoref{neat_tubular_nbhd_thm} and conclude the following.

\begin{cor}\label{tubular_nbhd_of_fixed_pts}
Suppose that $M$ is a compact $G$-manifold with $G$-trivial corners. Then
	\begin{enumerate}
		\item The submanifold $M^H$ has a neat $NH$-tubular neighborhood in $M$.
		\item If $M^K = \emptyset$ for all $K > H$,  then $M^{(H)}$ has a neat $G$-tubular neighborhood in $M$.
	\end{enumerate}
\end{cor}

\begin{rem}
	These statements would fail if we dropped the assumption that $G$ acts trivially on the corner coordinates.
\end{rem}

There is a useful way of compactifying the strata of a compact $G$-manifold. To begin, fix a total ordering of the conjugacy classes of subgroups
$$ \{e\} =(H_0) < (H_1) < \ldots < (H_k) = G$$
extending the subconjugacy partial ordering.

\begin{figure}[h]
	\centerline{\includegraphics[scale=1.1]{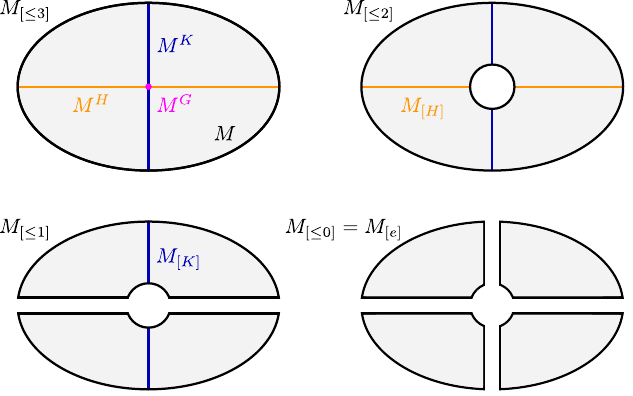}}
	\caption{Compactified strata}
	\label{defn_strata_figure}
\end{figure}

\begin{defn}\label{defn_strata_uncompactified}
For a $G$-space $M$, let $M_{\leq j}$ be the subspace of points such that the isotropy subgroup is conjugate to $H_i$ for some $i\le j$. %is at least as free as $H_j$, i.e. between $\{e\}$ and $H_j$ in the above total ordering. 
\end{defn}
These are open sets satisfying
$$
\emptyset = M_{\leq -1} \subseteq M_{\leq 0}\subseteq M_{\leq 1}\subseteq \dots \subseteq M_{\leq k} = M
$$
Thus the difference sets $M_{\le j}\backslash M_{\le j-1}$ are the invariant strata $M_{(H_j)}$.
\begin{con}\label{defn_strata}
Suppose that $M$ is a compact $G$-manifold with trivial corners. We make, for each $0\le j\le k$, a compact codimension zero $G$-invariant submanifold $M_{[\le j]}\subset M$. It will be contained in $M_{\le j}$ (that is, it will have no $H_i$-isotropy for any $i>j$), and it will be a compactification of (a smaller copy of) $M_{\le j}$. These sets will be nested:
$$
 M_{[\leq 0]}\subseteq M_{[\leq 1]}\subseteq \dots \subseteq M_{[\leq k]} = M.
$$
%that are essentially compactifications of the manifolds $M_{\leq j}$. 
%If $M$ is a compact $G$-manifold with trivial corners, we define a compactification of $M_{\leq j}$, denoted $M_{[\leq j]}$. More exactly, it is a compact submanifold of $M$ (with corners), invariant under the $G$-action, contained in $M_{\leq j}$, and having a dense open subset equivariantly diffeomorphic to $M_{\leq j}$. 
The construction is by downward induction on $j$. To begin, let $M_{[\le k]}$ be $M$. Now suppose that $M_{[\le j]}$ is given, and that it has no $H_i$-isotropy for $i>j$ and is a compactification of (a copy of) $M_{\le j}$. For $H\sim H_j$ let $M_{[H]}\subset M_{[\le j]}$ be the neat submanifold $M_{[\le j]}^{H}$. Let $M_{[(H_j)]}$ be $M_{[\le j]}^{(H_j)}$, the disjoint union of these. Using \autoref{tubular_nbhd_of_fixed_pts}, choose a neat $G$-invariant tubular neighborhood of $M_{[(H_j)]}$ in $M_{[\le j]}$. Let $M_{[\le j-1]}$ be $M_{[\le j]}\backslash T_j$, where $T_j$ is the open tubular neighborhood. The manifold $M_{[\le j-1]}$ is a compactification of $M_{[\le j]}\backslash M_{[\le j]}^{(H_j)}$ and therefore a compactification of $M_{\le j}\backslash M_{\le j}^{(H_j)}=M_{\le j-1}$. Likewise $M_{[(H)]}$ is a compactification of $M_{(H)}$, and $M_{[H]}$ is a compactification of $M_{H}$

%The difference set $M_{\leq j}\backslash M_{\leq j-1}$ is an open tubular neighborhood $T_j$ of a submanifold, denoted by $M_{[H_j]}$, that is a compactification of (something equivariantly diffeomorphic to) the stratum $M_{(H_ j)}$.

%We proceed inductively, beginning with $j=k$.
%\begin{itemize}
%\item $M_{[\leq k]}=M$.
%\item When $M_{[\leq j]} \subseteq M$ has been defined, then 
%\item Construct $M_{[\leq j-1]}$ by removing the interior of a neat $G$-tubular neighborhood of the neat submanifold $M_{[(H_j)]}$ of $M_{[\leq j]}$ using \autoref{tubular_nbhd_of_fixed_pts}.
%\item After $M_{[\leq j]} \subseteq M$ has been defined, then for any $H \in (H_j)$ we define $ M_{[H]}$ to be $(M_{[\leq j]})^H$, and $M_{[(H_j)]}$ is defined to be $(M_{[\leq j]})^{(H_j)} = \coprod_{H \in (H_j)} M_{[H]}.$
%The spaces $M_{[H]}$ for $H \sim H_j$ are disjoint from each other because all points with isotropy closer to $G$ have already been removed.
%\item $M_{[\leq j-1]}$ is obtained by removing a neat $G$-tubular neighborhood of the neat submanifold $M_{[(H_j)]}$ of $M_{[\leq j]}$ using \autoref{tubular_nbhd_of_fixed_pts}. 
%\end{itemize}
%Thus at the end $M_{[\leq 0]}$ has free $G$-action and is a compactification of $M_{\{e\}}$.

Observe that in general $M_{[\leq j]}$ will have corners, even if $M$ does. Note also that $M_{[H]}$, $M_{[(H_j)]}$, and $M_{[\leq j]}$ are embedded in $M$ but not neatly embedded. %(Parts of their boundaries map into the interior of $M$.) 
See \autoref{defn_strata_figure}.
\end{con}

In light of \autoref{blowup_removes_tubular_neighborhood}, we could instead define the stratum $M_{[\leq j-1]}$ by blowing up $M_{[\leq j]}$ along $M_{[(H_j)]}$. This would make the construction canonical. %, and makes it clear that $M_{[H]}$ is a compactification of the $H$-stratum $M_H$. %$M_{[(H)]}$ is a compactification of $M_{(H)}$.
However, we find it more convenient to proceed as above using tubular neighborhoods, partly because iterated blowups are tricky to discuss and partly because there is some advantage to having the compactified strata embedded in the manifold $M$. 

\begin{lem}\label{h_brackets_equivalent}
	There is an equivariant deformation retraction of $M_H$ onto the subspace $M_{[H]}$, and therefore the inclusion $M_{[H]} \to M_H$ is an equivariant homotopy equivalence.
\end{lem}
%	If $M \to N$ is an equivariant neat embedding with neat tubular neighborhood $D$, then $N \setminus D$ is homeomorphic to $N \setminus M$, and $\overline{(N \setminus D)} \to N \setminus M$ is an equviariant homotopy equivalence.
%\end{lem}

\begin{proof}
	We define the deformation retraction in each closed tubular neighborhood minus its zero section, by moving in the radial direction until everything goes to the boundary sphere bundle $S$. Perform this homotopy along each tubular neighborhood as we remove it from the manifold $M$ in \autoref{defn_strata}. The resulting homotopy is equivariant, is stationary on $M_{[H]}$, preserves the subspace $M_H$, and after finitely many steps sends all of $M_H$ into $M_{[H]}$, giving the desired deformation retraction. %Both of these are proven by expanding and shrinking in the radial direction of the neat tubular neighborhood.
\end{proof}

%\begin{cor}\label{h_brackets_equivalent}
%	The inclusions $M_{[H]} \to M_H$ are equivariant homotopy equivalences.
%\end{cor}

\subsection{The category of manifolds with compactified strata}

In \cite[Def 6.10]{GIMM1} we defined a category $\Man$ of manifolds and equivariant smooth embeddings:
\begin{defn}\label{m_emb}
	$\Man$ is the simplicially enriched category whose objects are compact smooth $G$-manifolds with corners (embedded in some fixed sufficiently large set), and whose $p$-simplices of morphisms are equivariant smooth embeddings $\Delta^p \times M \to \Delta^p \times N$ that commute with the projection to $\Delta^p$. Equivalently, they are smooth equivariant maps $M \times \Delta^p \to N$ that give smooth embeddings $M \to N$ for every $t \in \Delta^p$.
\end{defn}

As a tool for proving that our splitting of the equivariant $h$-cobordism space is natural, we now introduce a related category $\Manstr$, related to $\Man$ by a Dwyer-Kan equivalence, in which an object is a manifold $M$ with a choice of compactified strata $M_{[H]}$ embedded in $M$.

As before, choose a total ordering of the conjugacy classes $$ \{e\} =(H_0) < (H_1) < \ldots < (H_k) = G$$ extending the subconjugacy partial ordering.
	
\begin{defn}\label{m_strat}
	$\Manstr$ is the following simplicially enriched category. An object $M$ is a compact $G$-manifold with corners (that is, a $\Man$-object), plus a choice of submanifolds $M_{[\leq j]}$, $M_{[H]}$, $M_{[(H)]}$, and $T_j$ as in Construction \ref{defn_strata}. Let us call this a manifold with embedded compactified strata.
	
	%, for each $0 \leq j \leq k$, of a neat tubular neighborhood of the embedding
	%\[ M_{[(H_j)]} \to M_{[\leq j]}, \]
	%whose closed image we call $\overline{T}_j \subseteq M_{[\leq j]}$. These are defined recursively from $j = k$ to $j = 0$, because we define the submanifold $M_{[\leq j-1]} \subseteq M$ to be the complement of the open tubular neighborhood $T_j$.
	
	For any two objects $M$ and $N$, the morphism space $\Manstr(M,N)$ is the simplicial set in which a $p$-simplex is an equivariant smooth embedding $\Delta^p \times M \to \Delta^p \times N$ over $\Delta^p$, such that $f(\Delta^p \times M_{[\leq j]}) \subseteq \Delta^p \times N_{[\leq j]}$ for all $j$. It follows that $f(\Delta^p \times M_{[(H)]}) \subseteq \Delta^p \times N_{[(H)]}$ for all $H$ as well.
 \end{defn}

\begin{prop}\label{manstr_equivalent}
	The forgetful functor $\Manstr \to \Man$ is a Dwyer--Kan equivalence.
\end{prop}

The proof occupies the rest of this section. Notice that on objects the functor is surjective. %, because by \autoref{neat_tubular_nbhd_thm} it is possible to choose neat tubular neighborhoods of each manifold $M_{[(H_j)]}$ inside $M_{[\leq j]}$.  
We must show that it is an equivalence on mapping spaces.

Fix two manifolds $M$ and $N$ with embedded compactified strata, and consider the commuting square of simplicial sets
	\begin{equation}\label{manstr_diagram}
	\begin{tikzcd}
		\Manstr_{\textup{int}}(M,N) \dar \rar & \Manstr(M,N) \dar \\
		\Man_{\textup{int}}(M,N) \rar & \Man(M,N).
	\end{tikzcd}
	\end{equation}
The bottom horizontal map $\Man_{\textup{int}}(M,N) \subseteq \Man(M,N)$ is the inclusion of the sub-simplicial set consisting of those $p$-simplices of embeddings that land in the \emph{interior} of $N$. The top horizontal map $\Manstr_{\textup{int}}(M,N) \subseteq \Manstr(M,N)$ is the inclusion of the sub-simplicial set consisting of those $p$-simplices of embeddings such that $f(M_{[\leq j]})$ lands in the \emph{interior} of $N_{[\leq j]}$ for all $0 \leq j \leq k$. (These new simplicial sets in the left-hand column are not the morphism spaces of subcategories; the identity maps are lacking.)

\begin{lem}\label{manstr_interior}
	The horizontal inclusion maps in \eqref{manstr_diagram} are homotopy equivalences of simplicial sets.
\end{lem}

\begin{proof}
	For the lower map, this was proven in \cite[Thm 6.13]{GIMM1} using \cite[Prop 2.46]{GIMM1}. We recall the proof and then use the same idea for the upper map.
	
	The idea is to flow along a $G$-invariant vector field $\xi$ on $N$ that is inward-pointing, meaning that, in any coordinate chart with corner coordinates $[0,\infty)^k$, the vector $\xi(x)$ is positive in the $i$th corner coordinate whenever the point $x$ is zero in the $i$th corner coordinate. By flowing along $\xi$ we deform embeddings $M \to N$ to embeddings into the interior of $N$. We take this to be done all at once to all the embeddings, instead of doing it separately for each $\Delta^p$-family of embeddings rel $\partial \Delta^p$ as in  \cite[Thm 6.13]{GIMM1}. This provides both a map of simplicial sets $\Man(M,N) \to \Man_{\textup{int}}(M,N)$ and a homotopy between each composite and the identity, proving that the inclusion is a homotopy equivalence of simplicial sets.
	
	 The same method will show that the inclusion $\Manstr_{\textup{int}}(M,N) \to \Manstr(M,N)$ is also a homotopy equivalence. For this we need to choose the inward-pointing vector field on $N$ so that it is inward-pointing on $N_{[\leq j]}$ for all $j$. We will do this by first picking an inward-pointing vector field for the manifold with corners $N_{[\leq 0]}$, and then extending it to $N_{[\leq j]}$ one value of $j$ at a time, without changing the fact that it is inward-pointing for the previous values of $j$.
	
	Near any one point of $N_{[\leq j]} \setminus N_{[\leq (j-1)]}$, this means that we have a neat tubular neighborhood of the form
\[  \R^{m-k} \times D^{n-m} \times [0,\infty)^k \to  \R^{m-k} \times \R^{n-m} \times [0,\infty)^k, \]
	an inward-pointing vector field $\eta$ defined on its complement $\R^{m-k} \times (\R^{n-m} \setminus \textup{int }D^{n-m}) \times [0,\infty)^k$, and we want to extend it to an inward-pointing vector field on $\R^{m-k} \times \R^{n-m} \times [0,\infty)^k$. We do this by taking any extension of $\eta$ to a vector field on the complement of a slightly smaller closed ball $C^{n-m} \subseteq \textup{int }D^{n-m}$ that is still inward-pointing along $[0,\infty)^k$. This is possible because the condition of being inward-pointing is an open condition, and being inward-pointing along $(\R^{n-m} \setminus \textup{int }D^{n-m}) \times [0,\infty)^k$ implies being inward-pointing along $\R^{n-m} \times [0,\infty)^k$ as well.
	
	 Now let $\lambda_0$ and $\lambda_1$ be any a smooth partition of unity on $\R^{n-m}$ subordinate to the interior of $D^{n-m}$ and the complement of $C^{n-m}$. Let $\xi$ be any fixed inward-pointing vector field on $[0,\infty)^k$, and take the vector field $\lambda_0 \cdot (0,0,\xi) + \lambda_1 \cdot \eta$ on $\R^{m-k} \times \R^{n-m} \times [0,\infty)^k$. This agrees with $\eta$ on the complement of $D^{n-m}$, but now it is inward-pointing along $[0,\infty)^k$ everywhere, because inside $C^{n-m}$ it agrees with the inward-pointing field $(0,0,\xi)$, and in the remaining region $D^{n-m} \setminus C^{n-m}$ it is an average of two inward-pointing fields on $\R^{m-k} \times \R^{n-m} \times [0,\infty)^k$ and is therefore also inward-pointing.
	
	Doing this near each point and using a smooth partition of unity to solve the problem globally, we may extend every inward-pointing vector field on $N_{[\leq (j-1)]}$ to an inward-pointing field on $N_{[\leq j]}$. Repeating the argument for each $j$ gives a single vector field on all of $N$ that is inward-pointing for each of the compactified strata $N_{[\leq j]}$. The field may be made $G$-invariant, because the condition of being inward-pointing is preserved by averaging over the $G$-action.
	
	Now that we have such a field, the same argument shows that the inclusion $\Manstr_{\textup{int}}(M,N) \simeq \Manstr(M,N)$ is a homotopy equivalence, because flowing forward along this field gives an isotopy of embeddings $N \to N$, and composing with any of these embeddings preserves the property that $f(\Delta^p \times M_{[\leq j]}) \subseteq \Delta^p \times N_{[\leq j]}$ for all $j$.
\end{proof}
	
% \begin{lem}\label{manstr0}
% 	The forgetful map $\Manstr_0 \to \Man_0$ is an acyclic Kan fibration.
% \end{lem}
% 
% \begin{proof}
% 	By the neat tubular neighborhood theorem, any $\partial \Delta^p$ in $\Manstr_0$ going to the same manifold $M$ in $\Man_0$, can be filled to a $\Delta^p$ in $\Manstr_0$, working one stratum at a time. To apply the theorem, we do have to recognize that at each step the manifold $(\Delta^p \times M)_{[\leq j]}$ maps to $\Delta^p$ by a submersion, so by the Ehresmann fibration theorem, it is diffeomorphic over $\Delta^p$ to a trivial family $\Delta^p \times (M_{[\leq j]})$. Under this identification, we can regard the map $(\Delta^p \times M)_{[(H_j)]} \to (\Delta^p \times M)_{[\leq j]}$ as a $\Delta^p$-family of neat embeddings into a single manifold $M_{[\leq j]}$, so that \autoref{neat_tubular_nbhd_thm} applies. Iterating all the way down to $j = 1$ gives the complete $p$-simplex in $\Manstr_0$.
% \end{proof}
  
\begin{lem}\label{manstr1}
 	The forgetful map $\Manstr_{\textup{int}}(M,N) \to \Man_{\textup{int}}(M,N)$ is an acyclic Kan fibration between Kan complexes.
\end{lem}
 
\begin{proof}
	We first recall from \cite[Thm 6.13]{GIMM1} that $\Man_{\textup{int}}(M,N)$ is a Kan complex. To briefly recall the proof, suppose we have a horn $\Lambda^p_i \subseteq \Delta^p$ of equivariant embeddings $f\colon \Lambda^p_i \times M \to \Lambda^p_i \times (\textup{int }N)$. We can extend this to a $U$-family of embeddings for some open set $U$ about $\Lambda^p_i$ in $\Delta^p$, using smooth extension \autoref{smooth_extension} (this uses that the embeddings land in the interior of $N$) and \cite[Lem 2.18]{GIMM1}. Then we extend this to a $\Delta^p$-family of embeddings by pulling back along any smooth map $\Delta^p \to U$ that is the identity on $\Lambda^p_i$. This extends the $\Lambda^p_i$-family of embeddings to a $\Delta^p$-family of embeddings.
 	
 	Now it suffices to show that given any diagram
 	\[ \begin{tikzcd}
 		\partial \Delta[p] \dar \rar & \Manstr_{\textup{int}}(M,N) \dar \\
 		\Delta[p] \rar \ar[dashed]{ur} & \Man_{\textup{int}}(M,N)
	\end{tikzcd} \]
	the bottom map can be changed by a simplicial homotopy rel $\partial \Delta[p]$ to a map for which a dashed lift exists. Concretely, this means we are given a $\Delta^p$-family of equivariant embeddings $f\colon \Delta^p \times M \to \Delta^p \times (\textup{int }N)$, and a single collection (a $*$-family) of neat tubular neighborhoods of the strata of $M$ and $N$, such that along $\partial\Delta^p$ the embeddings all send $M_{[\leq j]}$ into the \emph{interior} of $N_{[\leq j]}$. Our goal is to deform the family of embeddings away from $\partial\Delta^p$ so that $f(M_{[\leq j]}) \subseteq \textup{int }N_{[\leq j]}$ across the entire family.
	
	This can be done by a sequence of deformations of $f$, so fix one value of $j$ and assume that $f(M_{[\leq j]}) \subseteq \textup{int }N_{[\leq j]}$ across the entire family. We want to deform the family so that $f(M_{[\leq j-1]}) \subseteq \textup{int }N_{[\leq j-1]}$ as well, without breaking the condition that $f(M_{[\leq k]}) \subseteq \textup{int }N_{[\leq k]}$ for all $k \geq j$. To do this, it suffices to find a vector field $\xi$ on $N$ that is inward-pointing for $N_{[\geq k]}$ for all $k \geq j$, and a function $\lambda \colon \Delta^p \to [0,\infty)$ supported on the interior of $\Delta^p$ such that if we flow along $\xi$ for time $\lambda$, the image $f(M_{[\leq j-1]})$ gets pushed into the interior of $N_{[\leq j-1]}$.
	
	We first define $\xi$ along $N_{[\leq j]}$. Since $f$ is an embedding, it cannot send points with stabilizer group smaller than $H_j$ to points with stabilizer group $H_j$. (In other words it is isovariant, see \autoref{isovariant} below.) Therefore $f(\Delta^p \times M_{[\leq j-1]})$ must be sent to the complement of $N_{[(H_j)]}$ in $N_{[\leq j]}$. In other words, it does not intersect the zero section of the tubular neighborhood $T_j$ about $N_{[(H_j)]}$. Since $\Delta^p \times M_{[\leq j-1]}$ is compact, there is some $\epsilon > 0$ such that $f(\Delta^p \times M_{[\leq j-1]})$ does not intersect the closed $\epsilon$-tube in $T_j$. By a tube lemma argument, there must also be some open set $U$ about $\partial\Delta^p$ in $\Delta^p$ such that $f(U \times M_{[\leq j-1]})$ is still contained in the interior of $N_{[\leq j-1]}$, in other words it does not intersect the closed tubular neighborhood $\overline{T_j}$.
	
	We take $\xi$ to be any vector field in the tubular neighborhood $\overline{T_j}$ that in the fiber direction $D^n$ points outwards in the radial direction, with magnitude 1, for every vector outside the $\epsilon$-tube, and inside the $\epsilon$-tube smoothly interpolates from this to the zero vector field along the zero section. Along the base direction of the tubular neighborhood, we allow $\xi$ to do anything, but we also ask that it be inward-pointing for $N_{[\leq j]}$. We see this is possible near every fiber of the tubular neighborhood, and these conditions are preserved by convex combinations, so we can patch them together to make such a vector field globally. Averaging over $G$ makes the field $G$-invariant as well.
	
	Then we extend $\xi$ to a $G$-invariant vector field on all of $N$ that is inward-pointing for $N_{[\geq k]}$ for all $k \geq j$ using the same technique as in the proof of \autoref{manstr_interior}. This arranges things so that flowing along $\xi$ for any amount of time will preserve the condition that $f(M_{[\leq k]}) \subseteq \textup{int }N_{[\leq k]}$ for all $k \geq j$, and flowing along $\xi$ until time 1 will be enough to push the image of any of the maps $f(\{t\} \times M_{[\leq j-1]})$ into the interior of $N_{[\leq j-1]}$. Therefore, if we take $\lambda\colon \Delta^p \to [0,1]$ to be any smooth function that is 1 on the complement of $U \subseteq \Delta^p$ and 0 on $\partial \Delta^p$, then flowing along $\xi$ for time $\lambda$ will deform our family rel $\partial \Delta^p$ so that in the interior we have $f(M_{[\leq j-1]})$ contained in the interior of $N_{[\leq j-1]}$. Repeating for all values of $j$ gives the desired deformation.
\end{proof}

\section{Isovariance and equivariant bundles}

In order to prove the splitting result for the stable space of equivariant $h$-cobordisms, we must consider \emph{isovariant} $h$-cobordisms. A splitting result for unstable spaces of isovariant $h$-cobordisms will yield the result by stabilization. In the stable limit, the distinction between equivariant and isovariant disappears.

\subsection{Isovariant and strictly isovariant maps}\label{isov}
Isovariance is a stronger condition than equivariance, in which not only fixed points but also strata are preserved. In this section we recall the definition of isovariance and related concepts, and we prove some technical results that will be needed later on.

%Suppose $f:X\to Y$ is $G$-equivariant. Then any time an element of $G$ fixes $x$, it also fixes $f(x)$. Isovariance is a stronger condition that makes this an ``if and only if.''

\begin{defn}\label{isovariant}
 Let $X$ and $Y$ be $G$-spaces. A continuous equivariant map $f\colon X \to Y$ is called \ourdefn{isovariant} if it also preserves stabilizers, that is, if $G_{f(x)} = G_x$ for all $x\in X$. Equivalently, $f$ is isovariant if $X^H=f^{-1}(Y^H)$ for every subgroup $H$. Equivalently, $f$ is isovariant if it preserves the strata of the isotropy stratification of $X$, in the sense that $f(X_H)\subseteq Y_H$ for every $H$, where $X_H$ is as in \autoref{fixed_points}.
 \end{defn}
 
In particular, an injective equivariant map is always isovariant.

\begin{defn}\label{isovariantheq}
	An isovariant map $f$ is an \ourdefn{isovariant homotopy equivalence} if it has an inverse up to isovariant homotopy, that is, if there exist an isovariant map $g:Y\to X$ and homotopies $gf \sim \id$ and $fg \sim \id$ through isovariant maps.
\end{defn}
 
 %Of course, a merely equivariant map satisfies $G_x\subseteq G_{f(x)}$ but not necessarily the reverse inclusion. Similarly, equivariance implies $X^H\subseteq f^{-1}(Y^H)$ for every subgroup $H\leq G$, and isovariance means $X^H= f^{-1}(Y^H)$. Therefore:
%\begin{lem}\label{isovariant_backwards}
%	An equivariant map $f\colon X \to Y$ is isovariant iff $f^{-1}(Y^H) \subseteq X^H$ for all $H \leq G$. In other words, $f$ is isovariant iff whenever $f(x)$ is fixed by $H$ then $x$ is also fixed by $H$.
%\end{lem}

%Yet another equivalent definition is that $f$ is isovariant iff $f$ preserves the strata of the isotropy stratification of $X$: $f(X_H)\subseteq Y_H$ for all $H \leq G$. 
 
 The next remark shows that smooth extension (\autoref{smooth_extension}) becomes false if ``equivariant'' is replaced by ``isovariant''. To address this difficulty we will introduce a stronger condition called \emph{strict isovariance} that takes tangential information into account. 
\begin{rem}
	Here is a counterexample to the na\" ive isovariant analogue of \autoref{smooth_extension}. The group is $\mathbb Z/2$. We give an isovariant smooth map $\partial M \to N$, where $M$ has corners and $N$ has no boundary at all, that cannot be extended to any isovariant smooth map defined on a neighborhood of $\partial M$ in $M$. This is closely related to the counterexample in \cite[Rem 2.10]{GIMM1}. 
	
	Let $N$ be $\R^\sigma$, the real line with $\Z/2$ acting by $x\mapsto -x$. Let 
$M$ be $[0,\infty)^3 \times \R^\sigma$, so that $\Z/2$ acts by $(x_1,x_2,x_3,x_4)\mapsto (x_1,x_2,x_3,-x_4)$. 

	Choose any number $c$ such that $1 < c < 2$, and define $f:M\to N$ by
		\[ f(x_1,x_2,x_3,x_4) = (x_1^2 + x_2^2 + x_3^2 - c(x_1x_2 + x_2x_3 + x_3x_1))x_4 + (c-1)x_4^3. \]
	The map is equivariant. To see that its restriction to $\partial M$ is isovariant, first note that
	$$
	x_1^2 + x_2^2 + x_3^2 - c(x_1x_2 + x_2x_3 + x_3x_1)\ge 0
	$$
	whenever $(x_1,x_2,x_3)\in \partial ([0,\infty)^3)$. Thus if $(x_1,x_2,x_3,x_4)\in \partial M$ is not fixed by the action then 
	$$
	x_1^2 + x_2^2 + x_3^2 - c(x_1x_2 + x_2x_3 + x_3x_1)+(c-1)x_4^2\ge (c-1)x_4^2 > 0,
	$$
so that $f(x_1,x_2,x_3,x_4)$ has the same sign as $x_4$ and is itself not fixed.
	
	Now suppose that $U$ is a $\mathbb Z/2$-invariant open neighborhood of $\partial M$ in $M$ and that $g:U\to N$ is an equivariant smooth map agreeing with $f$ on $\partial M$. We show that $g$ cannot be isovariant.
	
	Note that $g$ must agree with $f$ not only on the boundary but also at all points where $x_4$ is zero. (Equivariance forces it to vanish at such points.) Since $g$ and the cubic polynomial $f$ agree wherever any one of the four coordinates $x_i$ is zero, $f$ is the degree three Taylor polynomial of $g$ at the origin. Therefore the degree three Taylor polynomial of $g(t,t,t,t)$ is $f(t,t,t,t)=2(1-c)t^3$, which implies that $g(t,t,t,t)$ is negative for sufficiently small positive $t$. On the other hand, $g(0,0,0,t)=f(0,0,0,t)=(c-1)t^3 > 0$. It follows that for small positive $t$ the function $g$ must vanish somewhere on the line segment from $(0,0,0,t)$ to $(t,t,t,t)$, so that $g$ is not isovariant.
\end{rem}

 \begin{defn}\label{strictly_isovariant}
		 Let $M$ and $N$ be $G$-manifolds with corners. A smooth isovariant map $f\colon M \to N$ is \ourdefn{strictly isovariant} if for every $p\in M$ the $G_p$-equivariant linear map $D_pf:T_pM\to T_{f(p)}N$ is $G_p$-isovariant. Said differently, a strictly isovariant map $f:M\to N$ is any smooth isovariant map such that the induced map $Df:TM \to TN$ of total spaces of tangent bundles is also isovariant.
 \end{defn}
 
For example, if $f$ is an equivariant smooth embedding of $G$-manifolds, then $f$ is strictly isovariant. We also have the elementary observation: 

\begin{lem}\label{iso_to_trivial_product}
	If $M$, $N$, and $P$ are $G$-manifolds such that the $G$-action on $N$ is trivial, then a smooth equivariant map $M \to N \times P$ is strictly isovariant iff the component map $M \to P$ is strictly isovariant.
\end{lem}

%\begin{lem}\label{product_of_strict_iso}
%	If $f_1\colon M_1 \to N_1$ and $f_2\colon M_2 \to N_2$ are strictly isovariant then so is their product $f_1 \times f_2\colon M_1 \times M_2 \to N_1 \times N_2$.
%\end{lem}
%
%\begin{proof}
%	We observe that the product must be isovariant because for any $(x_1,x_2) \in M_1 \times M_2$, the stabilizer group is 
%	\[ G_{(x_1,x_2)}= G_{x_1} \cap G_{x_2} = G_{f(x_1)} \cap G_{f(x_2)} = G_{(f(x_1),f(x_2))}. \]
%	Then we apply this observation to the map induced on the tangent bundles, which is exactly the product map $Tf_1 \times Tf_2$. Since this is isovariant, $f_1 \times f_2$ is strictly isovariant.
%\end{proof}

The following is a useful criterion for isovariance of a \emph{linear} map.

\begin{lem}\label{isolinear}
Suppose that $V$ and $W$ are $G$-representations and $L\colon V\to W$ is a $G$-equivariant linear map. Then $L$ is isovariant if and only if every vector in the kernel of $L$ is fixed by the action of $G$. 
\end{lem}

\begin{proof}
Suppose that $L$ is isovariant. Then $G$ fixes every vector in the kernel, because if $L(v)=0$ then $G_v=G_{L(v)}=G_0=G$. 

Conversely, suppose that the action of $G$ on the kernel of $L$ is trivial. We show that $L$ is isovariant by showing that if $L(v)$ is fixed by a subgroup $H$ then $v$ is fixed by $H$. 

Assume that $L(v)$ is fixed by $H$. The vector $\bar v=\frac{1}{|H|}\sum_{h\in H}hv$ is fixed by $H$. The difference $v-\bar v$ is fixed by $G$, because it belongs to the kernel of $L$. Since $\bar v$ and $v-\bar v$ are fixed by $H$, $v$ is fixed by $H$ as well. %By the criterion of \autoref{isovariant_backwards}, therefore $T$ is isovariant.
%If $Tv=0$, then $G_{Tv}=G$, so if $T$ is isovariant we have $G_v=G$. Conversely, suppose $gTv=Tv$ so by equivariance of $T$ this is the same as $Tgv=Tv$, so $gv-v\in \ker T$. Suppose $g$ has order $n$ (we are assuming $G$ is finite throughout), and consider the telescoping sum $\sum_{i=0}^{n-1} g^i(gv-v)=0$. By assumption, since $G$ acts trivially on the kernel, $g^i(gv-v)=gv-v$, so $n(gv-v)=0$, and we get $gv-v=0$, as desired.
\end{proof}

Isovariance on the tangent space at a fixed point implies strict isovariance in a neighborhood of that point:
\begin{lem}\label{isonbhdG}
Suppose that $f\colon M\to N$ is a smooth $G$-equivariant map and that $p\in M$ is fixed by $G$. If  $D_pf\colon T_pM\to T_{f(p)}N$ is isovariant, then the restriction of $f$ to some $G$-invariant neighborhood of $p$ is strictly isovariant.

%If the kernel of $D_p f\colon T_pM\to T_{f(p)}N$ has trivial $G$-action, 
\end{lem}

\begin{proof}
The question is local, so that without loss of generality $M$ is a neighborhood of $0$ in a $G$-representation $V$ and the point $p$ is $0$. (In general we should say the intersection of such a neighborhood with a subset of the form $[0,\infty)^k \times V'$, but by the definition of smooth equivariant map $f$ may be extended to a neighborhood of $0$.) Likewise, $N$ may be assumed to be a neighborhood of $0=f(p)$ in a representation $W$.

%Choose coordinate charts around $p$ and $f(p)$ that identify each point with the origin inside $G$-representations $V$ and $W$, respectively. By the definition of a smooth map, we may without loss of generality assume that $f$ defines a smooth, equivariant map in a neighborhood of the origin, not just the intersection of this neighborhood with a subspace of the form $[0,\infty)^k \times V'$.

Write $V=V^G\times V_0$, so that $V_0$ is a representation with $V_0^G=0$, and likewise write $W=W^G\times W_0$. With respect to these decompositions, write 
$$f(v_1, v_2)= (f_1(v_1, v_2), f_2(v_1, v_2)).$$ 
Define a map $g$ from a neighborhood $0$ in $V^G\times V_0$ to $V^G\times  W_0$ by 
$$
g(v_1,v_2)=(v_1,f_2(v_1,v_2)).
$$
%Write the map $f$ with respect to the decompositions $V=V^G\times (V^G)^\perp$ and $W=W^G\times (W^G)^\perp$ as 
%$$f(v_1, v_2)= (f_1(v_1, v_2), f_2(v_1, v_2)).$$ 

By \autoref{iso_to_trivial_product}, the map $f$ will be strictly isovariant in a neighborhood of $0$ iff $f_2$ is strictly isovariant in the same neighborhood, iff $g$ is strictly isovariant in the same neighborhood. Therefore it suffices if $g$ when restricted to some neighborhood of $0$ becomes a smooth embedding. %Define a map from a neighborhood of $0$ in $V=V^G\times (V^G)^\perp$ to $V^G\times (W^G)^\perp$ by 
%$$
%g(v_1,v_2)=(v_1,f_2(v_1,v_2)).
%$$

For this it suffices if the linear map $D_0g$, derivative of $g$ at $0$, is injective, and for that it is enough if the derivative at $0$ of $v_2\mapsto f_2(0,v_2)$ is injective. This in turn follows from the fact (\autoref{isolinear}) that $\ker D_0f\subseteq (T_0M)^G=V^G\times 0$. 
% It follows that the derivative $D_0g$ is injective.  the restriction of $g$ to some neighborhood of $0$ is an embedding and therefore strictly isovariant. By \autoref{iso_to_trivial_product}, in the same neighborhood $f_2$ is also strictly isovariant, and therefore $f$ is strictly isovariant.
 %Indeed, for any fixed $v_0$ in a neighborhood of 0 in $V^G$, the derivative of $f_1(v_0, -)\colon (V^G)^\perp \to  (W^G)^\perp$ is by assumption injective, so by the inverse function theorem $f_1(v_0, -)$ is a local diffeomorphism from a neighborhood of 0 in $(V^G)^\perp$ to a neighborhood of 0 in $(W^G)^\perp$. Thus $f_1$ is a diffeomorphism on a neighborhood of 0 in $V$, and is therefore isovariant.
\end{proof}

%Generalizing the previous statement, we have:
\begin{cor}\label{isonbhdH}
Suppose that $f\colon M\to N$ is a smooth $G$-equivariant map, $p$ is a point in $M$, and $G_p= G_{f(p)}=H$. If $D_p f\colon T_pM\to T_{f(p)}N$ is $H$-isovariant, then there is a $G$-invariant neighborhood of $p$ on which the map $f$ is strictly $G$-isovariant.

\end{cor}

\begin{proof}
By \autoref{isonbhdG}, the point $p$ has a neighborhood $U$ on which $f$ is strictly $H$-isovariant. Let $g_i$ be coset representatives of $G/H$. If $U$ is chosen small enough, then the sets $g_i U$ are disjoint. In that case strict $H$-isovariance on $U$ implies strict $G$-isovariance on $GU=\bigsqcup_i g_i U$.
%Intersect this with a neighborhood of $p$ where all stabilizers with respect to the $G$-action are subgroups of $H$. 
\end{proof}

Here is another useful consequence of \autoref{isonbhdH}. It will help us to build strictly isovariant maps by an inductive process.
%\tmar{I like this more general formulation}
\begin{lem}\label{rank_isovariant_condition}
%	An equivariant map from a $G$-invariant neighborhood of $M^H$ in $M$ to $N$ 
Let $P\subset M$ and $Q\subset N$ be $G$-submanifolds. Let $f\colon M\to N$ be an equivariant map such that $f(P)\subset Q$. If the restricted map $P\to Q$ is strictly isovariant, and if for each $p\in P$ the induced map of normal spaces
\[ \nu_p(P \to M) \to \nu_{f(p)}(Q \to N) \]
is injective, then $f$ is strictly isovariant on some $G$-invariant neighborhood of $P$.
\end{lem}

\begin{proof}
By \autoref{isolinear} and \autoref{isonbhdH} it is enough if every vector in the kernel of $D_pf\colon T_pM\to T_{f(p)}N$ is fixed by $G_p$. By assumption any such vector is in $T_pP$, therefore in the kernel of $T_pP\to T_{f(p)}Q$, therefore fixed by $G_p$.
%Along the zero section, every vector in the kernel of the map on first derivatives has no component in the fiber direction, and is therefore $H$-fixed. By \autoref{isolinear} and \autoref{isonbhdG}, the map is therefore strictly isovariant in a neighborhood of each point in the zero section.
\end{proof}

\subsection{Equivariant fiber bundles and isovariant maps}\label{fiberbundles1}

Recall the notion of a smooth equivariant fiber bundle from \autoref{fiber_bundle}. In this subsection we give a couple of results on the interaction between such bundles and the concept of isovariance.

\begin{lem}\label{pullbackiso}
Suppose the following is a pullback diagram of equivariant spaces 
\[ \xymatrix{
	C \ar[d]_-{\ti g} \ar[r]^-{\ti f} & D \ar[d]^-g \\
	A \ar[r]^-f & B.
} \]
If $f$ is isovariant, then $\ti f$ is also isovariant.

If in addition all of the spaces are smooth manifolds, the maps between them are smooth, either $f$ or $g$ is a submersion, and $f$ is strictly isovariant, then $\ti f$ is also strictly isovariant.
\end{lem}

\begin{proof}
Let $(a,d)\in C=A\times_B D$. Then the stabilizer groups are
\[ G_{(a,d)}=G_a\cap G_d=G_{f(a)}\cap G_d=G_{g(d)}\cap G_d=G_d . \]
The statement for strict isovariance follows by applying the isovariance statement to the induced square of tangent bundles, which is also a pullback square of $G$-spaces since $f$ or $g$ is a submersion.
%Since $g$ is equivariant, $\Sta_d\subseteq \Sta_{g(d)}$, and since $f$ is isovariant, $\Sta_a=\Sta_{f(a)}$. By assumption $g(d)=f(a)$, thus we get that $\Sta_d\subseteq \Sta_a$. Hence $\Sta_a\cap \Sta_d=\Sta_d$ for any $a\in A$, so $\ti f$ is isovariant.
\end{proof}

\begin{lem}\label{bundlelemma}
Suppose $p\colon E\to B$ is an equivariant fiber bundle, $A$ and $B$ are $G$-CW complexes, and the following  is a pullback diagram
\[ \xymatrix{
	f^\ast E \ar[d]_-{\ti p} \ar[r]^-{\ti f} & E \ar[d]^-p \\
	A \ar[r]^-f & B.
} \]
 If $f$ is an isovariant homotopy equivalence then so is $\ti f$. %If furthermore $f$ is an isovariant deformation retraction of $A$ onto a subcomplex $B$, then $\ti f$ is also an isovariant deformation retraction of $f^*E$ onto $E$.
\end{lem}

\begin{proof}
Let $g\colon B\to A$ be an isovariant  homotopy inverse of $f$ and let $h\colon B\times I\to B$ be an isovariant homotopy from $\id_B$ to $fg$. Since the base spaces are $G$-CW complexes and are therefore equivariantly numerable, we can use the covering homotopy theorem, which gives an equivariant homotopy $\ti h$ extending $h$
\[ \xymatrix{
	E\times I \ar[d]_-{p\times I} \ar[r]^-{\ti h} & E \ar[d]^-p \\
	B\times I \ar[r]^-h & B.
} \] 
By the universal property of the pullback, the map $\ti h$ factors as
$$\ti h\colon E\times I\xrightarrow{\cong} h^\ast E \to E.$$
The first map is a bundle isomorphism and is therefore isovariant. The second map is also isovariant by \autoref{pullbackiso}, so $\ti h$ is isovariant. Restricting to $t=1$, we get a factorization
$$\ti h_1\colon E=E\times 1\xrightarrow{\cong} (fg)^\ast E\rightarrow f^\ast E \xrightarrow{\ti f} E,$$
which defines an isovariant right homotopy inverse $\ti g\colon E \to f^\ast E$ of $\ti f$.

Now, applying the same argument to the pullback diagram of equivariant bundles 
\[ \xymatrix{
	 E \ar[d]_-{} \ar[r]^-{\ti g} & f^\ast E \ar[d]^-{\ti p} \\
	B \ar[r]^-g & A.
} \]
we get that $\ti g$ has an isovariant right homotopy inverse. But we already know that $\ti f$ is its isovariant left homotopy inverse, thus $\ti f$ and $\ti g$ must be isovariant homotopy inverses.
%
%Finally, if $f$ is an isovariant deformation retraction then its homotopy inverse is an inclusion $g\colon B \to A$ such that $fg = \id$. In this case we may take the identity homotopy from $\id_B$ to $fg$, and then the resulting homotopy inverse $\ti g$ is actually a strict right inverse, in other words $\ti f \ti g = \id$. And when we go the other direction, we take the homotopy from $gf$ to $\id_A$ to be constant along $B \subseteq A$, so that in the covering homotopy theorem the bundle over $A \times I$ is identified with a trivial product along the subspace $B \times I$. This ensures that the resulting isovariant homotopy $\ti g\ti f \sim \id_{f^*E}$ is a homotopy rel $B$, in other words an isovariant deformation retraction.
\end{proof}

\section{Isovariant $h$-cobordisms and their stabilization}\label{sec:isovhcob}

In \cite{GIMM1} we defined the notion of equivariant $h$-cobordism on a compact smooth $G$-manifold $M$ with corners. By considering families of such cobordisms parametrized by a simplex, we defined a fibrant simplicial set $\mc H_{\eq,\sbt}(M)$, whose realization $\mc H_{\eq}(M)$ we called the equivariant $h$-cobordism space. We now define the isovariant $h$-cobordism space $\mc H_{\iso}(M)$, adapting the earlier definition in a straightforward way. When we wish to emphasize which group $G$ is acting on $M$, we will use the notation $\mc H_{G\dashiso}(M)$ for this space.

In \cite{GIMM1} we also considered equivariant $h$-cobordisms equipped with additional structure: ``mirror structure'' and ``encasing''. This led to spaces and homotopy equivalences 
$$
\mc H_{\eq}^c(M)\to \mc H_{\eq}^m(M)\to \mc H_{\eq}(M).
$$
The point of considering the extra structure was that $\mc H_{\eq}^c(M)$ is strictly functorial with respect to smooth embeddings $M\to M'$ (by a construction called ``polar stabilization''). This was used to make a space of equivariant $h$-cobordisms stabilized with respect to dimension. 

Here we extend all of this to the isovariant setting. Most of the work is routine, but there are a few points where we must make detailed arguments to verify that isovariance is preserved. %In order to explain this, we will briefly review the ``polar stabilization'' construction.

\subsection{Isovariant $h$-cobordisms}\label{isovhcob}
We begin by recalling the definition of equivariant $h$-cobordism from \cite{GIMM1}. Let $M$ be a compact smooth $n$-dimensional $G$-manifold with corners. The standard trivial cobordism on $M$ is taken to be $M\times [-1,0]$. (We use $[-1,0]$ instead of $[0,1]$ because later we will be doubling $W$ along the top $M\times \{0\}$ and we like $[-1,1]$ better than $[0,2]$.)

An \ourdefn{equivariant cobordism} on $M$ is a compact $(n+1)$-dimensional $G$-manifold $W$ with corners, equipped with the following structure. There is a face of $W$, called the bottom. There is a face $N$ of $W$, disjoint from $M$, called the top. The closure of the complement of $M\sqcup N$ in $\partial W$ is called the sides. There is a \ourdefn{lower collar}, a diffeomorphism between a neighborhood of the (union of the) bottom and sides of $W$ and a neighborhood of $M\times \{-1\}\cup\partial M\times [-1,0]$ in $M\times [-1,0]$. The choice of lower collar (or rather its germ along the bottom and sides) is part of the data of the cobordism. In particular the bottom and sides of $W$ are identified with those of the trivial cobordism $M\times [-1,1]$ in a definite way.

$W$ is called an \emph{equivariant $h$-cobordism} if the inclusions $M \to W$ and $N \to W$ are equivariant homotopy equivalences.

\begin{figure}[h]
	\centering
	\includegraphics{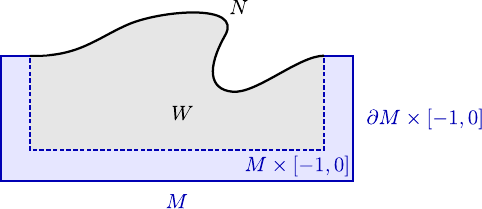}
	\vspace{-.5em}
	\caption{An $h$-cobordism on a manifold with corners $M$.}\label{fig:h-cobordism}
\end{figure}

%Recall from \autoref{isovariant} that an \ourdefn{isovariant homotopy equivalence} is an isovariant map that has an isovariant inverse up to isovariant homotopy.
\begin{defn}\label{hcob_def}
	An \ourdefn{isovariant $h$-cobordism} on $M$ is an equivariant cobordism as above such that
	the inclusions $M \hookrightarrow W$ and $N \hookrightarrow W$ are isovariant homotopy equivalences (\autoref{isovariantheq}).
\end{defn}

%See \autoref{fig:h-cobordism}.

\begin{rem}
Since every injective equivariant map is isovariant, in any equivariant cobordism $W$ the inclusions $M\hookrightarrow W$ and $N\hookrightarrow W$ are automatically isovariant. However, an isovariant map can be an equivariant homotopy equivalence without being an isovariant homotopy equivalence.

%if $W$ is an equivariant $h$-cobordism, it is not necessarily an isovariant $h$-cobordism, because the inverse homotopy equivalences and the relevant equivariant homotopies cannot necessarily be chosen to be isovariant. 
\end{rem}

Recall that a $G$-map $f\colon X \to Y$ between $G$-CW complexes is an equivariant homotopy equivalence if and only if it induces equivalences $X^H \to Y^H$ for all $H \leq G$. There is a similar criterion for isovariant homotopy equivalences in the case of $G$-manifolds in an $h$-cobordism. The following lemma is due to Hauschild, and translates easily to our particular definition of $h$-cobordism. Here $M_H$ means the $H$-stratum of $M$, from \autoref{fixed_points}.

\begin{lem}\cite[Proposition V.3]{hauschild}\label{strata_criterion}
An equivariant $h$-cobordism $(W,M,N)$ is an isovariant $h$-cobordism if and only if  $M_H \to W_H$ and $N_H\to W_H$ are  equivalences for all $H \leq G$.
\end{lem}
%\cmar{Do we need to assume $X \to Y$ is the bottom face of a cobordism for this?}\mmar{yes, we do. I don't know if the second statement is true in general.}
%The first of these claims generalizes to $G$-CW complexes, but the second is special to manifolds. 
This result is special to cobordisms. For general $G$-CW complexes one needs more data, %to measure the weak equivalences
see e.g.~\cite{sarah,inbar_sarah,dula_schultz}.

%The equivariant homotopy equivalence condition means that the inclusions of fixed-point manifolds $M^H \overset\sim\to W^H$ and $N^H \overset\sim\to W^H$ are homotopy equivalences for each $H \leq G$. The isovariant homotopy equivalence condition is equivalent to single-isotropy subspaces $M_H \overset\sim\to W_H$ and $N_H \overset\sim\to W_H$ being homotopy equivalences. (This is special to manifolds.)
Because every isovariant homotopy equivalence is an equivariant homotopy equivalence, every isovariant $h$-cobordism is an equivariant $h$-cobordism. The following example shows that not every equivariant $h$-cobordism is an isovariant $h$-cobordism. 
\begin{exmp}\label{gap_counterexample}
	Let $G=C_4 = \langle t\rangle$ be the cyclic group of order 4 and let $H$ be the subgroup $\langle t^2\rangle$. Let $G$ act freely on $S^3$ through the usual free $S^1$-action. Construct a free $G$-CW complex $\widetilde X$ by attaching cells to $S^3$: first a free orbit of 2-cells $G\times e^2$ and then a free orbit of 3-cells $G\times e^3$. The attaching map for $G\times e^2$ can be arbitrary. Choose the attaching map for $G\times e^3$ in such a way that the corresponding element of $\pi_2(S^3\cup (G\times e^2))\cong \Z[G]$ is $2 - t^2$. This element is not a unit in $ \Z[G]$, but it maps to a unit (namely $1$) in $\Z[G/H]$. 
	
	Let $X$ be the free $G/H$-space $\widetilde X/H$. The pair of free $G$-spaces $(\widetilde X,S^3)$ yields a pair of free $G/H$-spaces $(X,\R P^3)$. The inclusion $f\colon \mathbb{RP}^3 \to X$ is an integral homology equivalence, but it is not a homotopy equivalence, since the induced map of universal covers $S^3\to \widetilde{X}$ is not a homotopy equivalence; in fact, $H_2(\widetilde X,S^3)\cong\Z[G]/(2-t^2)\cong \Z/3 \oplus \Z/3$.
	
Now make a manifold model of the same kind. Beginning with the trivial cobordism on the free $G$-manifold $S^3\times I^2$, modify it by attaching handles to the top, first a free orbit of 2-handles and then a free orbit of 3-handles, in such a way that the boundary operator in the resulting relative handle chain complex is given by $2-t^2$. The result is a cobordism $\widetilde Y$ of free $G$-manifolds having $\widetilde Z=S^3\times I^2$ as its bottom face. Divide by the action of $H$ to get a cobordism $Y$ of free $G/H$-manifolds having $Z=\mathbb{RP}^3 \times I^2$ as its bottom face. The inclusion of that face is a homology equivalence but not a homotopy equivalence (for the same reason as above), so this is not an $h$-cobordism. We now enlarge $Y$ and $Z$ in a way that makes them simply connected and introduces some fixed points. 

The Hopf fibration $S^3\to S^2$ yields a circle bundle $\R P^3\to S^2$ by dividing by the $H$-action. (This is in fact the tangent circle bundle of $S^2$.) Let $D$ be the total space of the corresponding disk bundle. So $D$ is a compact 4-dimensional $G/H$-manifold whose fixed point set $D^{G/H}$ is the zero section $S^2$ of the disk bundle, and $\partial D$ may be identified (as a $G/H$-space) with $\R P^3$. Glue $D$ to $\R P^3\times I$ along $\partial D\cong \R P^3\times 0$ (the result is isomorphic to $D$), and let $M$ be the product of this with $I$, so that 
$$M = (D \times I) \cup_{(\partial D\times I)} Z\cong D\times I$$
The cobordism $Y$ was obtained by attaching handles to the top of $Z\times I$.  Let $W$ be obtained by attaching the same handles to the top of $M\times I$. In other words, $W$ is the result of gluing together the trivial cobordism $D \times I \hookrightarrow D \times I\times I$ and the nontrivial cobordism $Z \hookrightarrow Y$ along the trivial cobordism $\partial D\times I\hookrightarrow \partial D\times I\times I$. We show that $W$ is an equivariant $h$-cobordism. 

By construction, the fixed point set $W^{G/H}$ is a trivial cobordism $M^{G/H}\times I$. To see that $W$ is nonequivariantly an $h$-cobordism, first note that $M$ and $W$ are simply-connected, since $M$ is homotopy equivalent to $S^2$ and $W$ is obtained from $M\times I$ by attaching 2- and 3-handles. Therefore, since $H_j(W,M)\cong H_j(Y,Z)$ is trivial for all $j$, the inclusion $M\to W$ is a homotopy equivalence. Let $N$ be the top of $W$. It is simply-connected because the simply-connected space $W$ is obtained from $N\times I$ by attaching only 3- and 4-handles. Therefore the inclusion $N\to W$ is a homotopy equivalence, since $H_j(W,N)\cong H^{6-j}(W,M)$ is trivial.  

However, $W$ is not an isovariant $h$-cobordism, because the map on free strata $M_{\{e\}} \to W_{\{e\}}$ is not a homotopy equivalence, being
the same up to homotopy equivalence as the map $ Z\to Y$ or $\mathbb{RP}^3 \to X$ that we constructed above.

\end{exmp}

The issue raised by the above example goes away if we assume a gap hypothesis.

\begin{defn}\label{gap}
A $G$-manifold $M$ satisfies the \ourdefn{weak gap hypothesis} if for any isotropy subgroups $H\subsetneq K$ and any components $C\in \pi_0(M^H)$ and $D\in \pi_0(M^K)$ with $C_H, D_K\neq \emptyset$ and $D\subset C$, the inequality $\dim C-\dim D \geq 3$ holds.
\end{defn}

This hypothesis fails in \autoref{gap_counterexample}, where the difference in dimensions is only 2.

\begin{prop}\label{eqisiso}\cite[Theorem 4.51.]{luck}
Let $W$ be an equivariant $h$-cobordism on $M$ and assume $M$ satisfies the weak gap hypothesis. Then $W$ is an isovariant $h$-cobordism.
\end{prop}

\begin{rem}
Weak gap hypotheses are also needed for the equivariant $h$-cobordism theorem, but not the isovariant $h$-cobordism theorem. To be specific, by \cite[Theorem 4.42]{luck}, if $\dim M_H / WH \geq 5$ for all isotropy groups $H$ and $M$ satisfies the weak gap hypothesis, then an equivariant $h$-cobordism $W$ on $M$ is trivial if and only if its equivariant Whitehead torsion vanishes. On the other hand, by \cite[Theorem 5.51]{luck}, if $\dim M_H / WH \geq 5$ for all isotropy groups $H$, then an isovariant $h$-cobordism $W$ on $M$ is trivial if and only if its isovariant Whitehead torsion vanishes (without assuming any gap hypothesis on $M$).
\end{rem}

As in \cite[Definition 4.2]{GIMM1}, we use additional structures on cobordisms in defining spaces of cobordisms.

\begin{defn}\label{def_retract_height} \ %We define the following structures on isovariant $h$-cobordisms.
\begin{enumerate}
	\item The \ourdefn{double} $\bar W$ of an isovariant $h$-cobordism $W$ is two copies of $W$ glued along their common top face $N$. We think of one as ``flipped over'' and use the interval $[0,1]$ in the place of $[-1,0]$, with $M$ located at 1 and $N$ located at 0. 
	\item An isovariant \ourdefn{\mirror} $h$-cobordism is an isovariant $h$-cobordism $W$ together with a $(G \times C_2)$-equivariant smooth structure on $\bar W$ that restricts to the smooth structure of $W$ on each copy of $W$ in the double.
	
	\item Given an isovariant {\mirror} $h$-cobordism, an \ourdefn{encasing function} is a smooth, strictly $G$-isovariant map
	\[ \rho = (r,h)\colon W \to M \times [-1,0], \]
	having the following properties.
	\begin{itemize}
		\item The composition $\rho\circ c$ coincides with the identity map of $M\times [-1,0]$ in a neighborhood of the bottom and sides. In particular, $r:W\to M$ is a retraction to the bottom.
		\item The function $h$ satisfies $h^{-1}(-1) = M$, $h^{-1}(0) = N$, and the derivative $Dh$ has rank 1 along $N$. We call $h$ the \emph{height function}.
		\item The retraction $r\colon W \to M$ sends $W$ minus the sides to the interior of $M$.
		\item The double of $\rho$, in other words the map $\bar W \to M \times [-1,1]$ that extends $\rho$ and commutes with reflection, is smooth. (We may again denote this extended map by $\rho$, or by $(r,h)$.) Note that the double of $\rho$ maps a neighborhood of the boundary of $\bar W$ to a neighborhood of the boundary of $M\times [-1,1]$ by a diffeomorphism. 
	\end{itemize}
%	The function $h$ satisfies $h^{-1}(0) = N$, and the derivative $Dh$ has rank 1 along $N$; we call it the \emph{height function}. The map $ \bar W \to M \times [-1,1]$, that extends $\rho$ and commutes with reflection is smooth. (We sometimes denote this extended map by $\rho$, or by $(r,h)$.) Note that the extended $\rho$ maps a neighborhood of the boundary of $\bar W$ to a neighborhood of the boundary of $M\times [-1,1]$ by a diffeomorphism.
\end{enumerate}

\begin{rem}
By the trivial \autoref{iso_to_trivial_product}
%even though a map $(f_1, f_2)\colon X\to Y\times Z$ to a product is equivariant if and only if each factor $f_1$ and $f_2$ is equivariant, the same is not true about isovariance in general: a map to the product could be isovariant while one of the individual factors is not. However, in our construction of encasing functions,
the map $\rho = (r,h)\colon W \to M \times [-1,0]$ is strictly $G$-isovariant if and only if its first factor $r$ is strictly $G$-isovariant. And, although we will not need this fact, we observe that in $\bar W$ the only subgroups of $G\times C_2$ that can possibly occur as stabilizers of points or tangent vectors are those of the form $H \times \{e\}$ and $H \times C_2$. As a result, the doubled map $\bar W \to M \times [-1,1]$ is always strictly $(G \times C_2)$-isovariant. 
\end{rem}

An \ourdefn{encased} $h$-cobordism is a mirror $h$-cobordism equipped with \arhperiod This additional structure is needed to define stabilizations of $h$-cobordisms.

\end{defn}

\begin{figure}[h]
	\centering
	\includegraphics{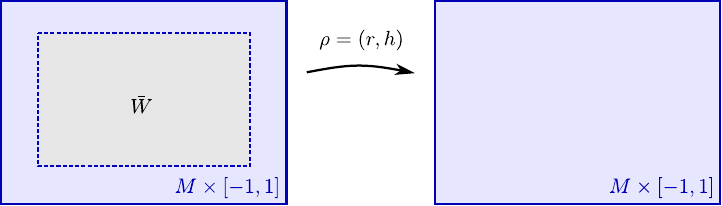}
	\vspace{-.5em}
	\caption{An encased $h$-cobordism.}\label{fig:encased}
\end{figure}
Note that when making an encased $h$-cobordism it is redundant to specify the lower collar $c$, since it must coincide with the inverse of \therh  $\rho$ near the boundary.

\subsection{Spaces of isovariant $h$-cobordisms}\label{isovhcobspaces}

\begin{defn}\label{hcob_space} A \ourdefn{family of isovariant $h$-cobordisms} parametrized by $\Delta^k$ is defined just as in \cite[Def 4.8]{GIMM1}; it is a smooth equivariant fiber bundle $E \to \Delta^k$ whose fibers are isovariant $h$-cobordisms, with a uniform choice of open set on which the lower collar is defined. 
%
%So it is a compact manifold with $G$-trivial corners $E$, an equivariant face $\Delta^k \times M \to E$, and an equivariant projection map $p\colon E \to \Delta^k$ that is a smooth fiber bundle, such that each fiber is an isovariant $h$-cobordism. Here $\Delta^k$ has the trivial $G$-action. In addition, we need a lower collar $c\colon \Delta^k \times M \times I \to E$ that satisfies the same conditions as before on each fiber of $p$, except that the open neighborhood of the bottom and sides must be open as a subset of $E$, not just on each fiber separately. For definiteness, each such family has as its underlying set a subset of $(\mc U \times \Delta^k)$, and the face and degeneracy maps handle this subset in the obvious way.
%
%A \ourdefn{family of encased $h$-cobordisms} is defined in the same way except fibers are encased $h$-cobordisms, and we assume the upper collar, retract, and height function are smooth across the family and that the germ of the upper collar is uniform across the family.
%
For set-theoretic safety, we fix a sufficiently large set $S$ containing $M$, and assume that the underlying set of each cobordism $W$ is a subset of $S$. We also assume that each $\Delta^k$-family of cobordisms is a subset of $\Delta^k \times S$, respecting the projection to $\Delta^k$, with the face and degeneracy maps acting on this subset in the obvious way.

Let $\mc H_{\iso,\sbt}(M)$ be the simplicial set in which a $k$-simplex is a family of isovariant $h$-cobordisms parametrized by $\Delta^k$. Let $\mc H_{\iso}(M)$ be its realization, the \ourdefn{space of isovariant $h$-cobordisms}. We sometimes write $\mc H_{G\dashiso}(M)$ to emphasize which group is acting on the cobordisms, as in \autoref{cobordisms_on_free_manifold} below. 

This is a subspace $\Hiso(M) \subseteq \mc H_{\eq}(M)$ of the space of equivariant $h$-cobordisms; in fact, it is the union of a set of connected components of the latter. Note that, in both the equivariant and isovariant cases, the face $M \times \Delta^k \to E$ and the lower collar are isovariant while the projection $p\colon E \to \Delta^k$ is merely equivariant.
\end{defn}

We also define the larger space of isovariant mirror $h$-cobordisms $\mc H_{\iso,\sbt}^m(M)$. As in \cite[Def 4.8]{GIMM1}, a $k$-simplex in this space is a $k$-simplex in $\Hiso(M) $ equipped with a smooth structure on the double (of the entire family). The forgetful map from this space to  $\mc H_{\iso,\sbt}(M)$ is an equivalence, by the same proof as \cite[Prop 4.13]{GIMM1}. Another way to say this is that $\mc H_{\iso,\sbt}^m(M)$ is a union of components in $\mc H_{\eq,\sbt}^m(M)$, and the equivalence from \cite[Prop 4.13]{GIMM1} can be restricted to these components.

We then define the space of isovariant encased mirror $h$-cobordisms $\mc H_{\iso,\sbt}^c(M)$. A $k$-simplex in this space has the additional structure of a strictly isovariant \rh
\[ \rho = (r,h)\colon\bar E \to M \times [-1,1] \]
having all of the properties from \autoref{def_retract_height} for each of the cobordisms $W \subseteq E$. The encasing is required to agree with the lower collar on a neighborhood of the bottom and sides in $E$ (not just on a set that intersects each $W$ in a neighborhood of the bottom and sides). This subspace  $\mc H_{\iso,\sbt}^c(M)\subset \mc H_{\eq,\sbt}^c(M)$ is not a union of components of the larger space, because an equivariant encasing function might fail to be strictly isovariant.

If $E$ is a trivial family $\Delta^k \times W$ of cobordisms (as it always is, up to isomorphism), then we may state these conditions a little more simply. We are asking for a $(G \times C_2)$-equivariant, strictly $G$-isovariant smooth map
\[ \rho = (r,h)\colon \Delta^k \times \bar W \to M \times [-1,1] \]
such that for each point of $\Delta^k$ it defines an encasing function for $\bar W$. We require there to be a single open neighborhood in which all of these \rhs agree with $c^{-1}$. Again by \autoref{iso_to_trivial_product}, $\rho$ is strictly isovariant if and only if the retraction $r$ is strictly isovariant. 

%\tmar{I would place Propositions 4.15 and 4.16 right here. Before beginning the long proof of 4.12, we should tell the reader what it is for (other than proving that something is a Kan complex, which is a mere remark)}

\begin{rem}\label{h_cob_space_kan}
	The spaces $\mc H_{\iso,\sbt}(M)$ and $\mc H_{\iso,\sbt}^m(M)$ are Kan complexes by the same proof as \cite[Lem 4.11]{GIMM1}.
\end{rem}

\begin{prop}\label{space_of_mirrorhcob}
The forgetful map from mirror to ordinary isovariant $h$-cobordisms
\[ \mc H_{\iso,\sbt}^m(M) \to \mc H_{\iso,\sbt}(M) \]
is a weak equivalence of Kan complexes.
\end{prop}

\begin{proof}
The proof is identical to the one for equivariant $h$-cobordisms in \cite[Prop 4.13]{GIMM1}.
\end{proof}

We need to know also that the forgetful map from $\mc H_{\iso,\sbt}^c(M)$ to $\mc H_{\iso,\sbt}^m(M)$ is a weak equivalence. This is important for later results because, although the encased version is the one on which the maps of our splitting are defined, in proving that they form equivalences we will use a fiber sequence in \autoref{fiberseq} that is proved using the unencased version.

To prove that these models for the $h$-cobordism space are equivalent, we will have to prove that the space of isovariant encasing functions is a contractible Kan complex (not just contractible), following the corresponding argument in \cite[Prop 4.16]{GIMM1}. The main thing to check is that the space of \emph{strictly isovariant} retractions is a contractible Kan complex. As we have indicated, this is a crucial result for the paper, but the proof is unforunately quite long. We have attempted to alleviate the experience for the reader with several illustrations.
%, which critically uses the fact that the space of encasing functions is a contractible Kan complex (not just contractible). We will now prove the corresponding fact for \emph{isovariant} encasing functions, which is unfortunately quite long.

%To show that $\mc H_{\iso,\sbt}^c(M)$ is a Kan complex, the main thing to check is that the space of \emph{strictly isovariant} retractions is contractible.

%To give the details: for a fixed {\mirror} $h$-cobordism $W$, a $k$-simplex of \rhs is a smooth map $\rho = (r,h)\colon \Delta^k \times \bar W \to M \times [-1,1]$ that for each point of $\Delta^k$ defines \arhperiod We require there to be a single open neighborhood on which the \rhs agree with $c^{-1}$, throughout the entire simplex $\Delta^k$. Note that $(r,h)$ must be strictly isovariant. By \autoref{iso_to_trivial_product}, this is equivalent to asking that the retract $r$ is strictly isovariant.

\begin{prop}\label{retracts_contractible}
	The space of strictly isovariant retractions on a given isovariant {\mirror} $h$-cobordism is a contractible Kan complex.
\end{prop}

%\cmar{problem: isovariant smooth extension is false, unless we take isovariant to mean strictly isovariant, in which case isovariant smooth approximation is false. I do think the space of strictly isovariant retracts will still be contractible.}
The proof of this uses similar ideas to \cite[Prop 4.6]{GIMM1} but is more elaborate because of the lack of an isovariant smooth approximation theorem and an isovariant smooth extension theorem. To get around this, and to ensure that the retractions we construct are isovariant, we define the retractions successively on a neighborhood of the $H$-stratum for different subgroups $H$, working one conjugacy class at a time.

\begin{proof}
	Let $W$ be an isovariant mirror $h$-cobordism on $M$ and let $\bar W$ be its double. Fix a total ordering on conjugacy classes of subgroups of $G$, $(G) = (H_k)$, $(H_{k-1})$, $\cdots$, $(H_0) = (e)$, which extends the partial ordering given by inclusion. We give a stratification of $W$, starting with the bottom and sides and successively adding sets of fixed points.  Let $W^{\geq k+1}$ be just the bottom and sides of $W$. For each $0 \leq i \leq k$, let $W^{\geq i} \subseteq W$ be the closed subspace of $W$ consisting of the bottom and sides and the fixed point set $W^K$ for all $(K) \geq (H_i)$. In the end $W^{\geq 0} $ is $W$. (This can be viewed as dual or complementary to the stratification from  \autoref{defn_strata_uncompactified}, where we were removing fixed points by larger subgroups in the ordering. See \autoref{defn_strata_tall_figure} for an illustration.) 
	
\begin{figure}[h]
	\centerline{\includegraphics[scale=.9]{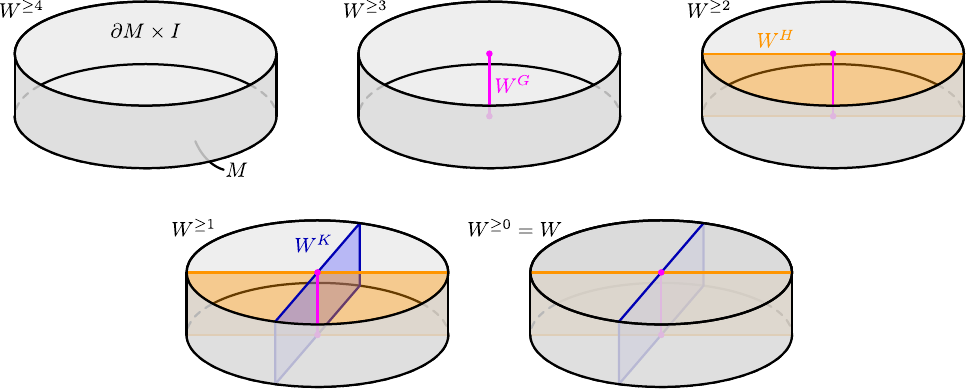}}
	\caption{A stratification of $W$ for an ordering $G>H>K>e$}
	\label{defn_strata_tall_figure}
\end{figure}

	For each value of $i$ we let $\bar W^{\geq i}$ be the corresponding subspace of the double $\bar W$. Note that $W^{\geq i}$ and $\bar W^{\geq i}$ are invariant under the $G$-action.
	
	Next, for each piece $\bar W^{\geq i} \subseteq \bar W$  in the stratification, we define the space of strictly isovariant retractions defined on a germ of an open neighborhood of $\bar W^{\geq i}$. For each $0 \leq i \leq k+1$, let $R^{\geq i}$ be the simplicial set in which a $p$-simplex is a $\Delta^p$-family of $(G \times C_2)$-equivariant, strictly $G$-isovariant smooth maps $\bar U^{\geq i} \to M$, defined on a $(G \times C_2)$-invariant open set $\bar U^{\geq i}$ containing $\bar W^{\geq i}$, which agree with the lower collar near the bottom and sides of $W$ (or equivalently near the boundary of $\bar W$), and which send the rest of $\bar W$ into the interior of $M$. The reader is encouraged to think about the lower half, where $U^{\geq i} = W \cap \bar U^{\geq i}$ is a $G$-invariant open set containing $W^{\geq i}$ and we have a strictly $G$-isovariant map $U^{\geq i} \to M$.
	
	\begin{figure}[h]
	\centerline{\includegraphics[scale=.9]{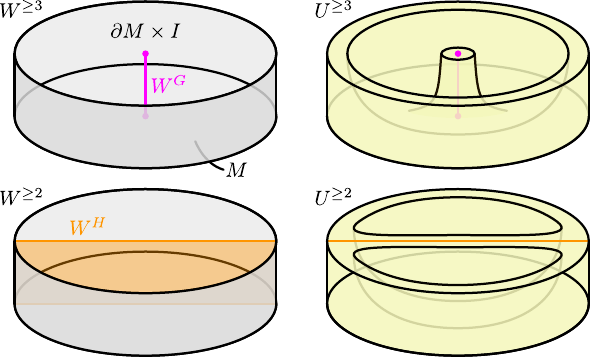}}
%	\caption{A stratification of $W$ for an ordering $G>H>K>e$}
%	\label{defn_strata_tall_figure}
\end{figure}

	Note that $R^{\geq k+1} = *$ is the one-point space because there is only one germ of a retraction defined near the boundary of $\bar W$ that agrees with the lower collar. Therefore it suffices to prove that for each $i$ the forgetful map $R^{\geq i} \to R^{\geq i+1}$ is an acyclic Kan fibration, since then we can conclude that $R^{\geq 0}$ is a contractible Kan complex, as desired. We have reduced the problem to constructing lifts of the following form.
	\begin{equation}\label{411_lift}
	\xymatrix{
		\partial \Delta^p \ar[d] \ar[r] & R^{\geq i} \ar[d] \\
		\Delta^p \ar[r] \ar@{-->}[ur] & R^{\geq i+1}
		}
	\end{equation}
	
	We are given a $\Delta^p$-family of $(G \times C_2)$-equivariant, strictly $G$-isovariant maps
	\[ r_0: \Delta^p \times \bar U^{\geq i+1} \to M, \]
	defined on some $(G \times C_2)$-invariant open set $\bar U^{\geq i+1} \subseteq \bar W$ containing $\bar W^{\geq i+1}$. We are also given a $\partial \Delta^p$-family of $(G \times C_2)$-equivariant, strictly $G$-isovariant maps
	\[ \partial \Delta^p \times \bar U^{\geq i} \to M, \]
	defined on a $(G \times C_2)$-invariant open set $\bar U^{\geq i}$ containing $W^{\geq i}$, and that along $\partial \Delta^p \times (\bar U^{\geq i+1} \cap \bar U^{\geq i})$ agrees with $r_0$. Without loss of generality $U^{\geq i+1} \subseteq U^{\geq i}$, so we have a single family of maps called
	\[ r_0\colon (\partial \Delta^p \times \bar U^{\geq i}) \cup (\Delta^p \times \bar U^{\geq i+1}) \to M. \]
	Our goal is to extend this family to $\Delta^p \times \bar U^{\geq i}$, after possibly shrinking $\bar U^{\geq i}$ to be a smaller $(G \times C_2)$-invariant open subset of $\bar W$ containing $\bar W^{\geq i}$.
	
%	When illustrating the argument, we sometimes only describe the bottom half $W$, where our open sets are $U^{\geq i} = W \cap \bar U^{\geq i}$, and our family of maps is strictly $G$-isovariant and has the form
%	\[ r\colon (\partial \Delta^p \times U^{\geq i}) \cup (\Delta^p \times U^{\geq i+1}) \to M. \]
	
	As in the proof of \cite[Prop 4.6]{GIMM1}, we first make a uniform choice of open neighborhood of the sides of $\bar W$ on which all of the maps in $r_0$ are the standard projection of the lower collar to $M$. By the tube lemma, this open set contains a product $S \times [-1,1]$ where $S \subseteq M$ is an open neighborhood of $\partial M$. Let $M_0 \subseteq M$ be a trimming with the closure of $M \setminus M_0$ contained in $S$. Let $\bar W_0 \subseteq \bar W$ be the corresponding submanifold in which we have removed the open set $(M \setminus M_0) \times [-1,1]$. Then $r_0$ is standard on $\bar W \setminus \bar W_0$, and so without loss of generality we can focus on $\bar W_0$ for the remainder of the proof.
	
	\begin{figure}[h]
	\centerline{\includegraphics[scale=.9]{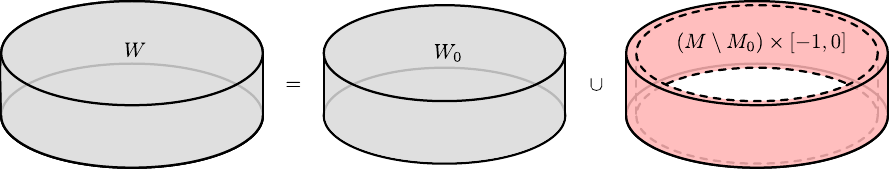}}
	%	\caption{A stratification of $W$ for an ordering $G>H>K>e$}
	%	\label{defn_strata_tall_figure}
	\end{figure}

	So let $\bar U_0^{\geq i} = \bar U^{\geq i} \cap \bar W_0$. We now have a family of smooth maps
	\[ r_0\colon (\partial \Delta^p \times \bar U_0^{\geq i}) \cup (\Delta^p \times \bar U_0^{\geq i+1}) \to \textup{int }M, \]
	and our goal is to extend them to $\Delta^p \times \bar U_0^{\geq i}$. The advantage of this is that now all of the maps land in the interior of $M$, so that we can use smooth extension and smooth approximation results.

	\begin{figure}[h]
	\centerline{\includegraphics[scale=.9]{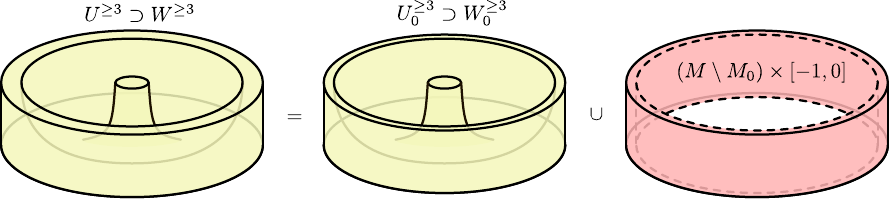}}
	%	\caption{A stratification of $W$ for an ordering $G>H>K>e$}
	%	\label{defn_strata_tall_figure}
	\end{figure}

	Next, as in \autoref{defn_strata}, we may define a descending sequence of $G$-invariant submanifolds $\bar W_{[\leq i]} \subseteq \bar W_0$ by letting $\bar W_{[\leq k]}$ be $\bar W_0$ and then successively letting $\bar W_{[\leq i]}$ be $\bar W_{[\leq i+1]}$ minus an open $(G \times C_2)$-equivariant tubular neighborhood of the manifold $\bar W_{[(H_i)]} = \bigcup_{H \in (H_{i+1})} \bar W_{[\leq i+1]}^{H}$. Intersecting with $W_0$ so that the $C_2$ action disappears, we recover \autoref{defn_strata} for the $G$-manifold $W_0$.
	
	\begin{figure}[h]
	\centerline{\includegraphics[scale=.9]{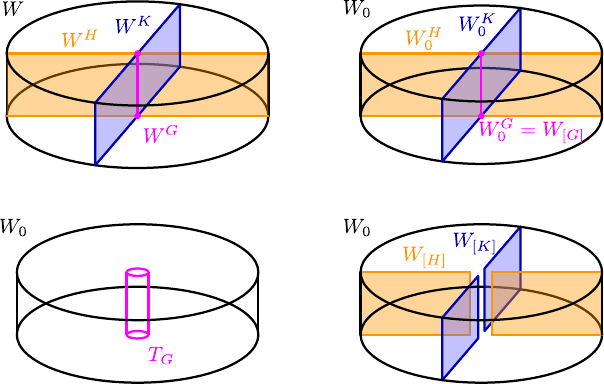}}
	%	\caption{A stratification of $W$ for an ordering $G>H>K>e$}
	%	\label{defn_strata_tall_figure}
	\end{figure}
	
	This construction depends on choices of tubular neighborhoods $\bar T_H$ about each manifold $\bar W_{[H]}$. Let us choose them so that the closure of each of the tubular neighborhoods $\bar T_H$ is contained in the open set $\bar U^{\geq i+1}$, for all $H \in (H_j)$ for $j > i$. Doing so guarantees that once we define $\bar W_{[(H_i)]}$ and its constituent pieces $\bar W_{[H]}$, the boundary $\partial \bar W_{[H]}$ is completely contained in the open set $\bar U_0^{\geq i+1}$. On the lower half, this means that the bottom and sides of $W_{[H]}$ are contained in $U_0^{\geq i+1}$, as shown.
	
	\begin{figure}[h]
	\centerline{\includegraphics[scale=.9]{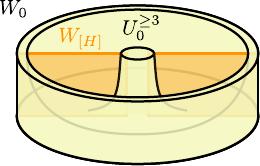}}
%	\caption{The bottom and sides of $W_{[(H_i)]}$ are contained in $U_0^{\geq i+1}$.}
	%	\label{defn_strata_tall_figure}
	\end{figure}
	
%	define a compact manifold with corners $\bar W_{[(H_i)]}$ as the $H$-fixed points of a manifold $\bar W_{[\leq i]}$, which is defined inductively by starting with $\bar W_{[\leq k]} = \bar W$ and successively removing tubular neighborhoods from the manifolds $W_{[(H_j)]}$ in $W_{[\leq j]}$ until we get to $W_{[\leq i]}$, and then we take the $H$-fixed points to get $W_{[H]}$. If we choose these tubular neighborhoods to be small enough, their closures will all be contained in the open set $U^{\geq i+1}$, since this open set already contains the sets $W_{(H_j)}$ for all $j > i$.
	
%	Let us also modify the definition of $W_{[H]}$ slightly by removing a collar of the sides of $W$ in $W^H$, chosen to be small enough that its closure is also contained in $U^{\geq i+1}$. As a result, sides of $W_{[H]}$ are completely contained in the open set $U^{\geq i+1}$, and they are also disjoint from the sides of $W$. (So the given family of retracts $r$ maps $\Delta^p \times W_{[H]}$ into the \emph{interior} of $M$.)

%	\begin{figure}[h]
%		\centerline{\includegraphics[scale=.9]{strata_nbhds.pdf}}
%		\caption{A stratification of $W$}
%		\label{defn_strata_tall_figure}
%	\end{figure}

	Next we let $\bar V_0^{\geq i} \subseteq U_0^{\geq i}$ be the smaller open set
	\[ \bar V_0^{\geq i} = \bar U_0^{\geq i+1} \cup \bigcup_{H \sim H_i} \bar T_H \]
	that is the union of $\bar U_0^{\geq i+1}$ and the open tubular neighborhoods of the strata $\bar W_{[H]}$ for all $H \sim H_i$. %This forms an open set that is ``segmented'' into pieces $U_{H_j}$, each of which is an open tubular neighborhood of a compact embedded manifold with corners $W_{[H_j]}$. We also have one piece $U_\partial$ that is the union of the collars on the sides of $W$. The open set $U^{\geq i+1^{\geq$ contains all of these pieces, save for the pieces $T_H$ for $H \sim H_i$, where it only contains the disc bundle on a neighborhood of the bottom and sides of $W_{[(H_i)]}$.
	\begin{figure}[h]
	\centerline{\includegraphics[scale=.9]{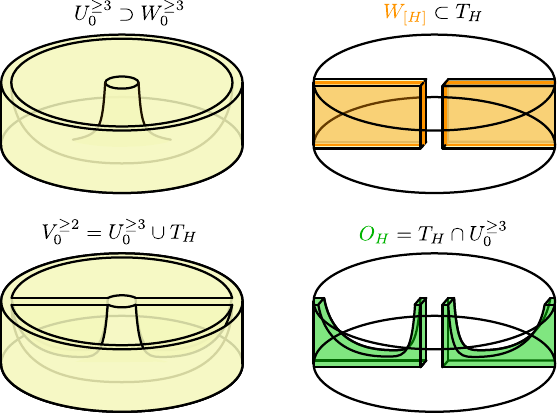}}
	%	\caption{A stratification of $W$ for an ordering $G>H>K>e$}
	%	\label{defn_strata_tall_figure}
	\end{figure}
	
	Of course $\bar W_{[(H_i)]}$ splits into disjoint pieces, one for each conjugate of $H = H_i$. Consider any one of these, $\bar W_{[H]}$, and its corresponding tubular neighborhood $\bar T_H$ inside $\bar W_0$. The manifold $\bar W_{[H]}$ has a free action by the Weyl group $WH = NH/H$, and the disc bundle $\bar T_H$ over it has a (non-free) action by the normalizer group $NH$. If we can define the family of retractions on $\bar T_H$, $NH$-isovariantly, it will extend to a family of $G$-isovariant retractions on the union of the components $\bar T_{(H_i)} = \bigcup_{H' \in (H_i)} \bar T_{H'}$.
	
	Let $\bar O_H = \bar T_H \cap \bar U_0^{\geq i+1}$; this is an open set in the disc bundle $\bar T_H$ that contains the disc bundle over $\partial\bar W_{[H]}$,
	\[ \bar T_H|_{\partial \bar W_{[H]}} \subset \bar O_H. \]
	We have reduced the problem to extending a given family of $(NH \times C_2)$-equivariant, strictly $NH$-isovariant retractions
	\[ r_0\colon (\partial \Delta^p \times \bar T_H) \cup (\Delta^p \times \bar O_H) \to \textup{int }M \]
	to the larger space $\Delta^p \times \bar T_H$, after shrinking the radius of the tube defining $\bar T_H$, and while preserving the original function in a neighborhood of $\bar T_H|_{\partial \bar W_{[H]}}$ (so that the extension on $\Delta^p \times \bar T_H$ and the given map on $\Delta^p \times \bar U^{\geq i+1}$ together give a smooth map on all of $\Delta^p \times \bar V^{\geq i}$).
	
	To do this, we first make a $(NH \times C_2)$-equivariant smooth extension to
	\[ r_1\colon (P_0 \times \bar T_H) \cup (\Delta^p \times \bar O_H) \to \textup{int }M, \]
	where $P_0 \subseteq \Delta^p$ is an open subset containing the boundary. We do this by the same method as the proof of \cite[Prop 4.6]{GIMM1}, by choosing an equivariant smooth extension from $\partial \Delta^p \times \bar T_H$ to $P_0 \times \bar T_H$ that does not necessarily agree with $r_0$ on $\Delta^p \times \bar O_H$, then interpolating between this extension and $r_0$ so that the result still agrees with $r_0$ on a smaller neighborhood of $\Delta^p \times \bar T_H|_{\partial \bar W_{[H]}}$. By \autoref{isonbhdH}, for $P_0$ sufficiently small this extension is still strictly $NH$-isovariant.
	
%	By isovariant smooth extension (\autoref{iso_smooth_extension}), we can extend to a strictly isovariant retract
%	\[ r\colon (P_0 \times \bar T_H) \cup (\Delta^p \times \bar O_H) \to \textup{int }M, \]
%	where $P_0 \subseteq \Delta^p$ is an open subset containing the boundary. Now our map is defined near the bottom and sides of $\Delta^p \times \bar T_H$, where ``sides'' refers to the boundary of $\Delta^p$ and the disc bundle over the boundary of $\bar W_{[H]}$, but not the boundary of the discs in each fiber of $\bar T_H$.
	
	The proof is therefore finished by the following lemma, where $\bar N = \Delta^p \times \bar W_{[H]}$, $D = \Delta^p \times \bar T_H$, and $\bar N_0 \subseteq \partial \bar N$ consists of $\partial \Delta^p \times \bar W_{[H]}$ and $\Delta^p \times \partial \bar W_{[H]}$. This extends $r_1$ to $\Delta^p \times \bar T_H$, and therefore to all of $\Delta^p \times \bar V^{\geq i}$, giving the desired lift in \eqref{411_lift}.
\end{proof}

Let $M$ be a $G$-manifold with corners. Let $H \leq G$ and $\bar N$ a compact smooth $(WH \times C_2)$-manifold with corners, that is the double of a free $WH$-manifold $N$ along one of its faces $N_t \subseteq \partial N$ called the ``top.'' In particular, every point and tangent vector in $\bar N$ has stabilizer group $H \times C_2$ or $H \times \{e\}$ for the given subgroup $H \leq G$.

Let $N_0 \subseteq \partial N$ be the closure of the complement of $N_t$ in $\partial N$, which we might call the ``bottom and sides'' of $N$. In particular, $\bar N_0$ is just the boundary of the double $\partial \bar N$. Let $\bar D \to \bar N$ be the disc bundle of a smooth $(NH \times C_2)$-equivariant vector bundle over $\bar N$, let $\bar D_0 = \bar D|_{\bar N_0}$ be the restricted disc bundle, and let $\bar D_1$ be any $(NH \times C_2)$-invariant neighborhood of $\bar D_0$ in $\bar D$.

	\begin{figure}[h]
	\centerline{\includegraphics[scale=.9]{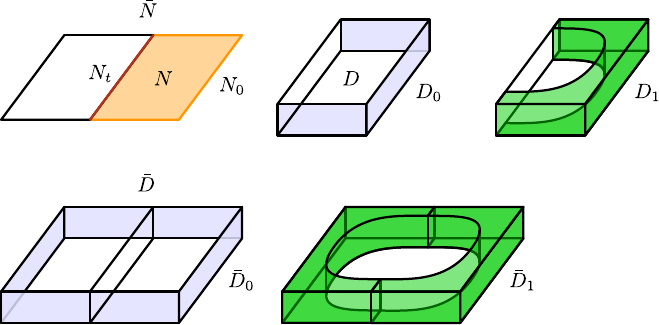}}
	%	\caption{A stratification of $W$ for an ordering $G>H>K>e$}
	%	\label{defn_strata_tall_figure}
	\end{figure}
	
\begin{lem}
	If $N_0 \to N$ is a $WH$-equivariant homotopy equivalence, then any smooth $(NH \times C_2)$-equivariant, strictly $NH$-isovariant map
	\[ r_1\colon \bar D_1 \to \textup{int }M \]
	can be extended to a smooth $(NH \times C_2)$-equivariant, strictly $NH$-isovariant map $\bar D \to \textup{int }M$, after possibly replacing $\bar D_1$ by a smaller neighborhood of $\bar D_0$, and shrinking the radius of the disc bundle $\bar D$.
\end{lem}
	
%\begin{lem}
%	Let $M$ be a $G$-manifold with corners, $H \leq G$, $N = N_H$ a compact free $WH$-manifold with corners, and $D \to N$ an $NH$-equivariant disc bundle of a vector bundle. Let $N_0 \subseteq \partial N$ be any union of faces that is equivariantly equivalent to $N$, $D_0 = D|_{N_0}$ the restricted disc bundle, and $D_1$ any neighborhood of $D_0$ in $D$. Then any strictly $NH$-isovariant map
%	\[ r\colon D_1 \to \textup{int }M \]
%	can be extended to a strictly $NH$-isovariant map $D \to \textup{int }M$, after possibly shrinking $D_1$ and the radius of $D$.
%\end{lem}

\begin{proof}
	Without loss of generality $D^H = N$, because otherwise we could re-define the zero section to be $D^H$. (In the proof just before, this was already true.)
	
	We construct the extension in two steps, extending first to the zero section $N \subseteq D$ and then to the rest of the bundle. %Let $D_1 \subseteq D$ be the open neighborhood of $D_0$ on which the map $r\colon D_1 \to \textup{int }M$ is defined.

	To extend to the zero section we use the technique from the equivariant proof in \cite[Lem 4.3]{GIMM1}. Let $\bar N_1 = \bar N \cap \bar D_1$. Then $r_1$ takes $\bar N_1$ into the stratum $\textup{int }M_H$ because $r_1$ is $NH$-isovariant. Once we think of its target as $\textup{int }M_H$, we can think of $r_1$ as just a $(WH \times C_2)$-equivariant map, since this is equivalent to being $(NH \times C_2)$-equivariant and strictly $NH$-isovariant.
	
	Take $\bar N_2 \subseteq \bar N_1$ to be a smaller closed neighborhood of $\bar N_0$ bounded by some manifold embedded in the interior of $N$, and that $(WH \times C_2)$-equivariantly deformation retracts onto $\bar N_0$. Let $\bar D_2 = \bar D|_{\bar N_2}$ be the corresponding disc bundle. If we choose this small enough then $\bar D_2 \subseteq \bar D_1$.
	
	Restricting attention to the lower half, since $N_0 \to N$ is a $WH$-equivariant homotopy equivalence and $N_2$ is an $WH$-equivariant subcomplex equivalent to $N_0$, it follows that $(N_2,N)$ is an $WH$-equivariant DR-pair. Therefore we get a continuous $WH$-equivariant map
	\[ r_2\colon N \to \textup{int }M_H \]
	that agrees with the restriction of $r_1$ to the bottom half of $\bar N_2$. Because $\textup{int }M_H$ has trivial $C_2$-action, this extends uniquely to a continuous $(WH \times C_2)$-equivariant map
	\[ r_2\colon \bar N \to \textup{int }M_H \]
	that agrees with $r_1$ on $\bar N_2$. Then we use equivariant smooth approximation (\autoref{smooth_approximation}) to deform $r_2$ to be a smooth $(WH \times C_2)$-equivariant map
	\[ r_3\colon \bar N \to \textup{int }M_H \]
	that agrees with $r_1$ on some smaller neighborhood $\bar N_3 \subseteq \bar N_2$ containing $\bar N_0$. As discussed before, $r_3$ must be strictly $NH$-isovariant because every point in both $\bar N$ and $\textup{int }M_H$ has isotropy $H$ (when we ignore the $C_2$ action).
	
	We next extend to $\bar D$. The map is already specified in the $TN$ direction, so we focus on the normal bundle $\bar \nu \to \bar N$ of the inclusion of $\bar N$ into $\bar D$, in other words the vector bundle associated to $\bar D$. Let $\xi \to \textup{int }M_H$ denote the normal bundle of the inclusion $\textup{int }M_H \to M$. The original retraction $r_1$, when restricted to $\bar D_3$, on first derivatives gives a vector bundle map
	\[ \xymatrix{
		\bar \nu|_{\bar N_3} \ar[d] \ar[r]^-{r_1} & \xi \ar[d] \\
		\bar N_3 \ar[r]^-{r_1} & \textup{int }M_H.
	} \]
	Since $r_1$ was strictly isovariant and $\bar N = \bar D^H$, by \autoref{rank_isovariant_condition} this vector bundle map is injective on each fiber. In other words, the derivative of $r_1$ defines a fiberwise linear embedding of bundles over $\bar N$
	\[ \bar \nu|_{\bar N_3} \ra ((r_3)^*\xi)|_{\bar N_3} \]
	and we want to extend it to the entire bundle, $\bar \nu \to (r_3)^*\xi$, equivalently a commuting diagram
	\[ \xymatrix{
		\bar \nu \ar[d] \ar@{-->}[r] & \xi \ar[d] \\
		\bar N \ar[r]^-{r_3} & \textup{int }M_H.
	} \]
	By \autoref{embedding_extension} applied to the lower half, since $N_0 \to N$ is an $NH$-equivariant homotopy equivalence, we get a smooth $NH$-equivariant extension of the map that we call $r_4$, which is a fiberwise linear embedding, which agrees with $r_3$ when we project down to $N$, and which agrees with the derivative of $r_1$ on some smaller neighborhood $N_4$ of $N_0$:
	\[ \xymatrix{
		\nu \ar[d] \ar[r]^-{r_4} & \xi \ar[d] \\
		N \ar[r]^-{r_3} & \textup{int }M_H.
	} \]
	The double of this is a continuous $(NH \times C_2)$-equivariant extension to all of $\bar \nu$. By \autoref{embedding_approximation}, we can deform this to a smooth $(NH \times C_2)$-equivariant fiberwise linear embedding $r_5$ that still agrees with the derivative of $r_1$ on $\bar N_4$:
	\[ \xymatrix{
		\bar \nu \ar[d] \ar[r]^-{r_5} & \xi \ar[d] \\
		\bar N \ar[r]^-{r_3} & \textup{int }M_H
	} \]
	Along a tubular neighborhood of $M^H$ in $M$, we can also think of $r_5$ as a smooth map $\bar D \to M$ that agrees with $r_3$ along $\bar N$, and that in $\bar D_4$ is given by the first derivative of $r_1$ along $\bar N_4$ in the normal direction.
	
	Using a smooth partition of unity subordinate to $\bar D_4$ and to $\bar D \setminus \bar D_0$, we can patch $r_5$ and $r_1$ together to get a smooth $(NH \times C_2)$-equivariant function $r_6\colon \bar D \to \textup{int }M$. As in the proof of \cite[Lemma 2.28 and Prop 2.42]{GIMM1} this interpolation is carried out in a representation $V$ containing a tubular neighborhood $\Omega$ of an embedded copy of $\textup{int }M$. To get a map back to $\textup{int }M$ we have to apply the projection $\Omega \to \textup{int }M$, so the interpolation is only really defined on an open subset of $\bar D$ containing the zero section $\bar N$, given by the preimage of $\Omega$. However, by simply shrinking the radius of $\bar D$ we can guarantee that this open set is all of $\bar D$.
	
	This linear interpolation process does not change the first derivative on any vector where the derivatives of $r_1$ and $r_5$ agree, so in particular the first derivative of $r_6$ is still injective in the normal direction to $\bar N$ inside $\bar D$. By \autoref{rank_isovariant_condition} therefore $r_6$ is therefore strictly isovariant in a neighborhood of $\bar N$.
\end{proof}

\begin{prop}\label{encased_to_mirror}
The forgetful map from encased to mirror isovariant $h$-cobordisms
\[ \mc H_{\iso,\sbt}^c(M) \to \mc H_{\iso,\sbt}^m(M) \]
is a weak equivalence of Kan complexes.
\end{prop}

\begin{proof}
As in \cite[Prop 4.16]{GIMM1}, this follows if we know that the space of all isovariant encasings $(r,h)$ on a fixed isomorphic $h$-cobordism is a contractible Kan complex. This is the product of the space of all isovariant retractions $r$ and the space of all height functions $h$. For the space of height functions the proof is the same as in \cite[Prop 4.6]{GIMM1}. For the space of retractions, this is \autoref{retracts_contractible}.
\end{proof}

\begin{cor}\label{heqishiso}
Suppose $M$ is a compact smooth $G$-manifold that satisfies the gap hypothesis from \autoref{gap}. Then the inclusion of the isovariant $h$-cobordism space into the equivariant space is an isomorphism
$$\Hiso^c(M)\xrightarrow{\cong} \mc H_{\eq}^c(M).$$
\end{cor}
\begin{proof}
For each family of equivariant $h$-cobordisms, by \autoref{eqisiso}, over each point of $\Delta^k$ the equivariant $h$-cobordism $W$ over that point is also an isovariant $h$-cobordism. Therefore the simplicial sets defining $\Hiso(M)$ and  $\mc H_{\eq}(M)$ consist of the same sets at every level, so that $\Hiso(M)\xrightarrow{\cong} \mc H_{\eq}(M)$ is an isomorphism. The claim now follows from \autoref{space_of_mirrorhcob} and \autoref{encased_to_mirror}. % therefore the corresponding map for encased $h$-cobordisms is a weak equivalence. (It is not an isomorphism because strictly isovariant retracts form a different space then equivariant retracts.)
\end{proof}

\begin{cor}\label{heqishisobigV}
Suppose $M$ is a compact smooth $G$-manifold. Then there exists some representation $V$ of $G$, such that$$\Hiso^c(M\times D(V))\xrightarrow{\cong} \mc H_{\eq}^c(M\times D(V)).$$
\end{cor}
\begin{proof}
For a large enough representation $V$, the compact smooth manifold $M\times D(V)$ will satisfy the gap hypothesis, thus the result follows by the previous corollary.\end{proof}

The proof of the following is a simple inspection of the definition. We state it for encased $h$-cobordisms, but the analogous statements for plain and mirror $h$-cobordisms are also true.

\begin{lem}\label{disj_unions_to_products}
	The isovariant $h$-cobordism space sends disjoint unions to products
\[ \Hiso^c\left( \coprod_{i=1}^n M_i \right) \cong \prod_{i=1} \Hiso^c(M_i) \]
	and similarly for the equivariant $h$-cobordism space.
\end{lem}%\kmar{It seems to me that the product is homeomorphic to the H-iso of the coproduct. It is the same space. Also, we wrote it that way in (2) before Thm 6.3.}

%As a result of \autoref{ehresmann_with_corners}, every $\Delta^k$-family of $h$-cobordisms can be trivialized, along with the lower collars.

%For a family of encased $h$-cobordisms, we don't assume that the upper collar and retract and height function can be trivialized as well, i.e. we don't assume the family is of the form $\Delta^k \times W$ for a single retract and height function on $W$.

When the $G$-action on $M$ is free, or more generally has a single isotropy type, the space simplifies to a nonequivariant one:
\begin{lem}\label{cobordisms_on_free_manifold}
	If $M$ is a smooth manifold with $G$-trivial corners and the $G$-action on $M$ is free then there is a canonical homotopy equivalence %\mmar{I'm not sure about this notation with the upperscript for the group... let's think about a good notation}
	\[ \mc H_{G\dashiso}^c(M) \xrightarrow{\cong} \mc H_{G\dasheq}^c(M) \xrightarrow{\simeq} \mc H^c(M/G). \]
	More generally, if $M$ has one orbit type, that is, if $M = M_{(H)}$ for some $H \leq G$, so that the isotropy group of every point is conjugate to $H$, then there is a canonical equivalence to the space of non-equivariant encased $h$-cobordisms
	\[ \xymatrix{
		\mc H_{G\dashiso}^c(M) \ar[d]^-\simeq \ar[r]^-\cong & \mc H_{G\dasheq}^c(M) \ar[d]^-\simeq \\
		\mc H_{WH\dashiso}^c(M^H) \ar[r]^-\cong & \mc H_{WH\dasheq}^c(M^H) \ar[r]^-\simeq & \mc H^c((M^H)/WH)
	} \]
\end{lem}

\begin{proof}
	Since $M$ has a single isotropy group, the gap hypothesis is trivially satisfied and so the isovariant and equivariant spaces coincide.
	
	In the free case we map $\mc H_{G\dashiso}^c(M)  \to \mc H^c(M/G)$ by sending each cobordism $W$ to the quotient $W/G$, and similarly for families of cobordisms over $\Delta^k$. This preserves the collars, retractions, and height functions. The inverse map sends a cobordism $W'$ over $M/G$ to the pullback
	\[ \xymatrix{
		W \ar[d] \ar[r] & W' \ar[d]^-r \\
		M \ar[r] & M/G.
	} \]
	Again this receives canonical collars, retractions, and height functions. Going back and forth gives the same cobordism, and all of its attached data, up to canonical isomorphism. These canonical isomorphisms give a homotopy between each composite and the identity, so that $\mc H_{G\dashiso}^c(M) \to \mc H^c(M/G)$ is a homotopy equivalence.
	
	In the single-orbit-type case, note that $M$ is a disjoint union of spaces $M^{H'}$ for $H' \sim H$, and so each cobordism $W$ is a disjoint union of cobordisms $W^{H'}$ over each of these components. Therefore if we map $\mc H_{G\dashiso}^c(M) \to \mc H_{WH\dashiso}^c(M^H)$ by restricting the cobordism $W$ over $M$ to $W^H$ over $M^H$, this operation has an inverse sending $W'$ over $M^H$ to $G \times_{NH} W'$ over $M$. Going back and forth gives the same cobordism up to canonical isomorphism, so this restriction map $\mc H_{G\dashiso}^c(M) \to \mc H_{WH\dashiso}^c(M^H)$ is a homotopy equivalence. We follow up with the equivalence from first part of the proof, $\mc H_{WH\dashiso}^c(M^H) \to \mc H((M^H)/WH)$, and we are done.
\end{proof}

%\begin{lem}
%	If $M$ is a manifold with $G$-trivial corners and $H \leq G$ then restriction to the $H$-fixed points
%	\[ \Hiso^G(M) \to \Hiso^{WH}(M^H) \]
%	is a Kan fibration whose fiber is equivalent to isovariant $G$-$h$-cobordisms on the compactification of $M - M^H$.
%\end{lem}
%\cmar{Doesn't apply to equivariant $h$-cobordisms}
%
%\begin{proof}
%\cary{Showing it's a Kan fibration is easy, but the fiber is not the same. To show it's the same we use that the space of tubular neighborhoods is weakly contractible.}
%\end{proof}

\subsection{The stabilization of an  isovariant $h$-cobordism}\label{isovhcobstab}
In \cite{GIMM1} we gave a detailed treatment of functoriality for the space of encased mirror equivariant $h$-cobordisms. We now restrict that construction to the subspace of isovariant $h$-cobordisms. We must show that the main construction in the proof, the stabilization of an $h$-cobordism, preserves isovariant $h$-cobordisms.

%Before we recall the definition, let us fix some notation that will appear in several of the proofs later in this paper.
%\begin{notn}\label{complement}
%	If $\nu$ is any vector bundle and we add a trivial line bundle to make $\nu \times \R$, we define the following subspaces of the product $D(\nu) \times [-1,1]$:
%	\begin{align*}
%		D(\nu \times \R) &= \{ (x,u) \in \nu \times \R : \|x\|^2 + u^2 \leq 1 \} \\
%		S(\nu \times \R) &= \{ (x,u) \in \nu \times \R : \|x\|^2 + u^2 = 1 \} \\
%		C(\nu \times \R) &= \{ (x,u) \in \nu \times \R : \|x\|^2 + u^2 \geq 1, \ \|x\| \leq 1, \ |u| \leq 1 \}
%	\end{align*}
%	In other words, $D(\nu \times \R)$ is just the unit disc and $S(\nu \times \R)$ is the unit sphere, but $C(\nu \times \R)$ is the complement of the interior of $D(\nu\times \R)$ inside $D(\nu) \times [-1,1]$. Note that the entire product $D(\nu) \times [-1,1]$ is the union of $C(\nu\times \R)$ and $D(\nu\times \R)$ along $S(\nu\times \R)$.
%\end{notn}

The stabilization takes each $h$-cobordism on $M$ and smooth embedding $e\colon M\to M'$ to an $h$-cobordism on $M'$ . When the embedding has codimension 0, the extension $\St^e(W)$ merely tacks on a trivial cobordism over $M' \setminus M$, see \cite[Definition 4.13]{GIMM1}. It is easy to see that this preserves isovariant $h$-cobordisms.

The only other case we need to consider is when we embed $M$ into the unit disc bundle $D(\nu)$ of a vector bundle $p\colon \nu \to M$. Let $S(\nu)\to M$ be the corresponding unit sphere bundle.
%Then the stabilization $\St^\nu(W)$ is We recall the definitions of an encased isovariant $h$-cobordism $\St^\nu(W)$ on the manifold $D(\nu)$  from \cite[\S 5]{GIMM1}, which is defined by defining its double $\St^\nu(\bar{W})$. Then the  the double $\St^e(\bar W)$ of the stabilization $\St^e(W)$ along the embedding $e$ is defined by first forming $\St^\nu(\bar W)$ and then further extending along the codimension 0 embedding $D(\nu)\to M'$ to get the double of a mirror $h$-cobordism over $M'$.
Let $\nu \times \R$ be the product of $p$ with the trivial bundle.

\begin{notn}\label{complement}
	We let $C(\nu \times \R)$ denote the complement of the interior of $D(\nu\times \R)$ inside $D(\nu) \times [-1,1]$, so that $D(\nu) \times [-1,1]$ is the union of $C(\nu\times \R)$ and $D(\nu\times \R)$ along $S(\nu\times \R)$. Explicitly:
	\begin{align*}
		D(\nu \times \R) &= \{ (x,u) \in \nu \times \R : \|x\|^2 + u^2 \leq 1 \} \\
		S(\nu \times \R) &= \{ (x,u) \in \nu \times \R : \|x\|^2 + u^2 = 1 \} \\
		C(\nu \times \R) &= \{ (x,u) \in \nu \times \R : \|x\|^2 + u^2 \geq 1, \ \|x\| \leq 1, \ |u| \leq 1 \}
	\end{align*}
\end{notn}

The disk bundle $D(\nu \times \R) \to M$ has a map to $[0,1]$, namely the norm in the fiber over $M$. On the other hand, $W$ also has a map to $M\times [0,1]$, namely the encasing function $\rho$ to $M \times [-1,0]$ flipped upside down. We can therefore form the fiber product
\[ W \times_{M \times [0,1]} D(\nu \times \R). \]
%This space will be the essential part of $St^\nu ( \bar{W})$.
%which as in  \cite[\S 5.2]{GIMM1} can be rewritten as the colimit
%\[ \xymatrix{
%	N \times_M S(\nu \times \R) \ar[d] \ar[r] & N \\
%	W \times_M S(\nu \times \R). &
%	} \]

%Let $C(\nu \times \R)$ denote the complement of the interior of $D(\nu\times \R)$ inside $D(\nu) \times [-1,1]$, so that $D(\nu) \times [-1,1]$ is the union of $C(\nu\times \R)$ and $D(\nu\times \R)$ along $S(\nu\times \R)$.

\begin{defn}\label{hcob_stab}
	As a topological space with $G$-action, the double of the stabilization $\St^\nu(\bar W)$ is the colimit of the diagram
%\[ \xymatrix{
%	& N \times_M S(\nu \times \R) \ar[d] \ar[r] & N \\
%	S(\nu \times \R) \ar[d] \ar[r]^-{} & W \times_M S(\nu \times \R) & \\
%	C(\nu \times \R), & &
%	} \]
%	or more simply the colimit of
	\[ \xymatrix{
	S(\nu \times \R) \ar[d] \ar[r]^-{} & W \times_{M \times [0,1]} D(\nu \times \R) \\
	C(\nu \times \R). &
	} \]
The encasing
$ \St^\nu \rho\colon\  \St^\nu(\bar W) \to D(\nu) \times [-1,1] $
is defined by projecting away from $W$ to the pushout
\[ \xymatrix{
	S(\nu \times \R) \ar[d] \ar[r]^-{} & D(\nu \times \R) \\
	C(\nu \times \R), &
} \]
which is $D(\nu) \times [-1,1]$.

	\end{defn}
%\noindent	The middle horizontal map is induced by the inclusion of $M\cong M \times \{-1\}$ into $W$. 

By \cite[Proposition 5.10]{GIMM1}, the space $\St^\nu(\bar W)$ has a well-defined smooth structure and by \cite[Proposition 5.13]{GIMM1}, the encasing function is smooth. Indeed, the entire reason for working with doubles of $h$-cobordisms is to give this stabilization a canonical smooth structure. We refer the reader to \cite[\S 5]{GIMM1} for details of these proofs and further discussion of the subtleties that arise and the ideas needed to fix them.

\begin{figure}[h]
	\centering
	\includegraphics{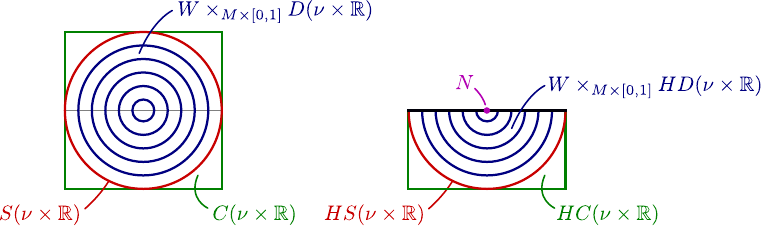}
	\vspace{-.5em}
	\caption{The doubled stabilization $\St^\nu(\bar W)$ and its lower half $\St^\nu(W)$.}\label{stabfig}
\end{figure}

The group $C_2$ acts on $\St^\nu(\bar W)$ by flipping the $\R$ direction. The set of $C_2$-fixed points $\St^\nu(\bar W)^{C_2}$ is the fiber product
\[ W \times_{M \times [0,1]} D(\nu), \]
%which can also be expressed as the pushout of the diagram 
%\[ \xymatrix{
%	N \times_M S(\nu) \ar[d] \ar[r] & N \\
%	(\bar W \times_M S(\nu))/C_2 &
%	} \]
%and
which in \autoref{stabfig} corresponds to the center line. The whole is the double of the lower half along its top, and this lower half is called $\St^\nu(W)$. The following notation helps us describe this lower half $\St^\nu(W)$ more effectively. %We will show that  $\St^\nu(\bar W)$ is the double of a mirror $h$-cobordism from $D(\nu)$ to $\St^\nu(\bar W)^{C_2}$.

\begin{notn}\label{half_disc}
	As in \autoref{complement}, we let $HD(\nu \times \R)$ refer to the lower half-disc of the bundle $\nu \times \R$, and similarly $HS$ is the lower half-sphere and $HC$ is the lower half-complement:
	\begin{align*}
		HD(\nu \times \R) &= \{ (x,u) \in \nu \times \R : \|x\|^2 + u^2 \leq 1, \ u \leq 0 \} \\
		HS(\nu \times \R) &= \{ (x,u) \in \nu \times \R : \|x\|^2 + u^2 = 1, \ u \leq 0 \} \\
		HC(\nu \times \R) &= \{ (x,u) \in \nu \times \R : \|x\|^2 + u^2 \geq 1, \ \|x\| \leq 1, \ u \in [-1,0] \}
	\end{align*}
\end{notn}

Then with this notation, the lower half $\St^\nu(W)$ can be written as
\[ \St^\nu(W) = \big( W \times_{M \times [0,1]} HD(\nu \times \R) \big) \cup_{HS(\nu \times \R)} HC(\nu \times \R) \]
or as the pushout of the diagram
	\[ \xymatrix{
	HS(\nu \times \R) \ar[d] \ar[r]^-{} & W \times_{M \times [0,1]} HD(\nu \times \R) \\
	HC(\nu \times \R). &
	} \]

%Next, we prove that $\St^\nu(\bar W)$ is indeed the double of an $h$-cobordism. Note that the space $\St^\nu(\bar W)$ has a $C_2$-action which flips the sign of the $\R$ coordinate. 

%We start with some preliminary lemmas.
%\cmar{Moved to section 3}

We already proved in \cite[Proposition 5.8]{GIMM1} that for an equivariant $h$-cobordism $W$, the stabilization $\St^\nu(\bar W)$ is the double of an equivariant mirror $h$-cobordism on $D(\nu)$. So to prove the isovariant version of this claim, all that remains is to show that the inclusions of the top and bottom of $\St^\nu(W)$ are isovariant homotopy equivalences, and that the encasing function for $\St^\nu(W)$ is \emph{strictly} isovariant. We begin with this result:
 
\begin{lem}\label{cylinder_isovariant}
	If $\nu \to M$ is an equivariant vector bundle, then the mapping cylinder $C(D(\nu) \to M)$ admits a deformation retraction onto its end $D(\nu)$ that is both isovariant and fiberwise over $M$.
\end{lem}

\begin{proof}
	Let $tv$ denote the scaling of the vector $v \in D(\nu)$ by the number $t \in [0,1]$. Note that when $t = 0$ this produces the corresponding point in the zero section of $D(\nu)$ in the same fiber as $v$.
	
	Define the retraction
	\[ r\colon M \times_{(D(\nu) \times \{0\})} D(\nu) \times [0,1] \to D(\nu) \]
	by the formula $r(v,t) = (tv,1)$, for $v \in D(\nu)$ and $t \in [0,1]$. Since the $v$ is being scaled, its projection to $M$ is preserved, and since points with both $v \neq 0$ and $t > 0$ are sent to points with $v \neq 0$, the isotropy of each point is preserved, so this map is isovariant.
	
	We deform this continuously back to the identity via two straight-line homotopies
	\[ (tv,1) \sim (tv,t) \sim (v,t). \]
	Again, since the vectors are merely being scaled the projection to $M$ is preserved, and since points with both $v \neq 0$ and $t > 0$ are preserved, the map is isovariant.
%	\mmar{double check the formulas, we want to go in to the pink line, then collapse that line so that stabilizers don't change}
\end{proof}
	
\begin{prop}\label{stabiso}
	Let $W$ be an encased isovariant $h$-cobordism with isovariant retract $r$. Then $\St^\nu(\bar W)$ is the double of an encased isovariant $h$-cobordism on $D(\nu)$ with top given by $\St^\nu(\bar W)^{C_2}$.	\end{prop}

\begin{proof}
We first check that the stabilization is an isovariant $h$-cobordism. The stabilization $\St^\nu(W)$ is built out of three pieces, as depicted in \autoref{isovdef}. The lowest piece is $HC(\nu \times \R)$, which has an isovariant deformation retract to its bottom $D(\nu)$, and also to its top $HS(\nu \times \R)$, by moving along the $\R$ direction. Note that the top boundary of this region $HS(\nu\times \R)$ is also $G$-homeomorphic to $D(\nu)$.

\begin{figure}[h]
	\centerline{\includegraphics[scale=1.3]{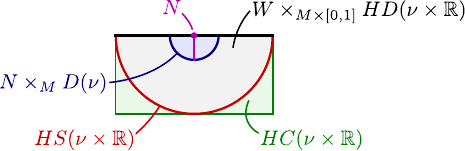}}
	\caption{Isovariant deformation}
	\label{isovdef}
\end{figure}

The second piece is not all of the fiber product $W \times_{M \times [0,1]} HD(\nu \times \R)$, only the part of it lying over a sub-interval of the form $[\epsilon,1]$. This region of the half-disc is homeomorphic to the product $[\epsilon,1] \times D(\nu)$, and so this piece of the fiber product is equivariantly homeomorphic to $W \times_M D(\nu)$.

Applying \autoref{bundlelemma} to the square
\[ \xymatrix{
	W \times_M D(\nu) \ar[d]^-{\ti p} \ar[r]^-{\ti r} & D(\nu) \ar[d]^-p \\
	W \ar[r]^-r & M,
} \]
since $W$ is an isovariant $h$-cobordism over $M$ and $r$ is an isovariant retract, we conclude that the natural projection of $W \times_M D(\nu)$ to $D(\nu)$ is an isovariant homotopy equivalence.

The final piece is the shaded piece in \autoref{isovdef}, the fiber product $W \times_{M \times [0,1]} HD(\nu \times \R)$ restricted to $[0,\epsilon]$. In this region we may assume we are in the top collar of $W$, so that the cobordism is the trivial cobordism $N \times [-1,0]$. Then the fiber product simplifies up to $G$-homeomorphism to the mapping cylinder of the projection
\[ N \times_M D(\nu) \to N. \]
This can be rewritten as the pullback $N \times_M (-)$ applied to the mapping cylinder from \autoref{cylinder_isovariant}. By that lemma, the inclusion of the bottom $N \times_M D(\nu)$ into this mapping cylinder is an isovariant homotopy equivalence.

Composing these all together, the inclusion of the middle region into the entire cobordism admits an isovariant deformation retract, and is therefore an isovariant equivalence. Therefore the inclusion of $HS(\nu \times \R)$ into the entire cobordism is an isovariant homotopy equivalence. This is isovariantly homotopic to the inclusion of the bottom $D(\nu)$, so the inclusion of the bottom is an isovariant homotopy equivalence as well.
%\cmar{I feel a little iffy about the way these equivalences are composed here. We don't quite know that a pushout of an isovariant h.e. is an isovariant h.e.. We do know this for maps that admit deformation retracts though.}

%We also have to check that the retract $\St^\nu r\colon \St^\nu(W) \to D(\nu)$ is an isovariant homotopy equivalence, but this follows since its precomposition with the isovariant homotopy equivalence $D(\nu) \to \St^\nu(W)$ is just the identity map of $D(\nu)$.

It remains to show that the inclusion of the top $\St^\nu(\bar W)^{C_2}$ into $\St^\nu(W)$ is also an isovariant homotopy equivalence. We proceed in three steps as before. We start with the top shaded region; relative to its top, this is also homeomorphic to the mapping cylider of $N \times_M D(\nu) \to N$, so by the same argument the inclusion of this top has an isovariant deformation retract. Next we consider the cobordism $W \times_M D(\nu)$ relative to its top and sides, not just its top. However we already know the inclusion of the top is an equivalence, and the inclusion of the top of the sides into the sides is an equivalence by the same argument, so by 2 out of 3 the inclusion of all of these into the cobordism is an equivalence. Finally, as we have already observed, the region $HC(\nu \times \R)$ has an isovariant deformation retract to its top $HS(\nu \times \R)$. Composing these together gives the desired isovariant homotopy equivalence.

It remains to show that the encasing function $\St^\nu(\bar W) \to D(\nu) \times [-1,1]$ is strictly isovariant if the original encasing function $(r,h)\colon \bar W \to M \times [-1,1]$ was. Using \cite[Lemma 5.15]{GIMM1}, we can locally model this encasing function by either
\begin{itemize}
	\item an identity map (in the trivial region and the frontier),
	\item the projection
	\begin{align*}
		(W \setminus (M \amalg N)) \times_{M \times [0,1]} D(\nu \times \R) &\to (M \times [0,1]) \times_{M \times [0,1]} D(\nu \times \R) \\
		&\cong D(\nu \times \R)
	\end{align*}
	 (in the nontrivial region), or
	\item the product of $r|_U\colon U \to M$ and the identity of $D(\nu \times \R)$, where $U$ is an open subset of the top of the cobordism $N \subseteq W$ (in the cone locus).
\end{itemize}
The first of these is trivially a strictly isovariant map, The second of these is a composition of a smooth equivariant inclusion and a pullback of the strictly $G$-isovariant map $(r,h)\colon W \to M \times [0,1]$, so it is strictly isovariant by \autoref{pullbackiso}. The last of these is similarly an inclusion $U \to W$ followed by a pullback of the strictly isovariant map $r\colon W \to M$, so it is strictly isovariant as well. The remaining properties of $\St^\nu(\bar W)$ being encased are already checked in \cite[Lemma 5.18]{GIMM1}.
\end{proof}

Now that we have shown in \autoref{stabiso} that the stabilization construction for $h$-cobordisms preserves isovariance, the rest of the results in \cite{GIMM1} go through and allow us to conclude the main result \cite[Theorem 6.18]{GIMM1} for spaces of isovariant $h$-cobordisms. Recall from \autoref{m_emb} the definition of the simplicial category $\Man$ of smooth compact $G$-manifolds with corners and smooth embeddings.

\begin{thm}\label{functor}
There is a functor $$\mc H_{G\dashiso} (-)\colon \Man \to \sSet\, ,$$ sending each compact $G$-manifold with corners $M$ to a space equivalent to $\mc H_{G\dashiso}^c(M)$, and each homotopy class of smooth equivariant embeddings $M_0 \to M_1$ to the homotopy class of maps $\mc H_{G\dashiso}^c(M_0) \to \mc H_{G\dashiso}^c(M_1)$ given by the stabilization defined in this section.
\end{thm}

\begin{prop}\label{stabilization_commutes_with_fixed_points}\label{stabilization_commutes_with_orbits}
	The $H$-fixed points and $G$-orbits define natural transformations
	\[ \mc H_{G\dashiso}(M) \ra \mc H_{WH\dashiso}(M^H), \qquad \mc H_{G\dashiso}(M) \simar \mc H_{1\dashiso}(M/G), \]
	where the second one is only defined on the subcategory where the $G$-action is free. Furthermore, for each individual manifold $M$ this map is homotopic to the map of $h$-cobordism spaces induced by taking the fixed points or orbits, respectively.
\end{prop}

Therefore by \autoref{cobordisms_on_free_manifold}, the first natural transformation is an equivalence whenever $M$ has isotropy concentrated in the conjugacy class $(H)$, while the second natural transformation is always an equivalence.

\begin{proof}
	The first natural transformation follows from \cite[Corollary 6.22]{GIMM1} by restricting to those components that are isovariant $h$-cobordisms. The second natural transformation also follows from \cite[Corollary 6.22]{GIMM1} since on a free $G$-manifold every equivariant $h$-cobordism is an isovariant $h$-cobordism.
\end{proof}

\begin{prop}\label{stabilization_commutes_with_coproducts}
	There is a natural weak equivalence of functors $\Man \times \Man \to \sSet$
	\[ \mc H_{G\dashiso}(M \amalg N) \simar \mc H_{G\dashiso}(M) \times \mc H_{G\dashiso}(N). \]
\end{prop}

\begin{proof}
	The proof is the same as the proof of \cite[Corollary 6.22]{GIMM1}, only much easier. In the language of that proof, we form a pullback square of left fibrations
	\[ \xymatrix{
		\Hanst \times \Hanst \ar[d]^-{p \times p} \ar[r]^-{\amalg} & \Hanst \ar[d]^-{p} \\
		\Manst \times \Manst \ar[r]^-{\amalg} & \Manst.
	} \]
	Here $\Manst$ is the category of manifolds and stabilization data from \cite[Def 6.11]{GIMM1}, where an object is a compact smooth $G$-manifold $M_0$ and a morphism is a zig-zag $M_0 \leftarrow D_0 \subseteq M_1$, where $D_0 \to M_0$ has a certain structure called a ``round bundle'' that remembers the distance to the origin \cite[Def 5.2]{GIMM1}. The functor called $\amalg$ on the bottom of our diagram sends $(M_0,N_0) \mapsto M_0 \amalg N_0$, and sends two round tubular neighborhoods $M_0 \leftarrow D_0 \subseteq M_1$ and $N_0 \leftarrow E_0 \subseteq N_1$ to their disjoint union
	\[ M_0 \amalg N_0 \leftarrow D_0 \amalg E_0 \subseteq M_1 \amalg N_1. \]
	The map on the left still has the structure of a ``round bundle'' in a canonical way. This similarly extends to the $p$-simplices of morphisms as well, which are defined as maps $M_0 \times \Delta^p \leftarrow D_0 \subseteq M_1 \times \Delta^p$ with the same properties.
	
	The functor called $\amalg$ on the top of our square is given on the simplicial set of objects by the canonical isomorphism $\mc H_{G\dashiso}(M \amalg N) \cong \mc H_{G\dashiso}(M) \times \mc H_{G\dashiso}(N)$ from \autoref{disj_unions_to_products}. On the simplicial set of morphisms, we apply the disjoint union functor $\amalg$ as above to turn every $p$-simplex of morphisms in $\Hanst \times \Hanst$ to a $p$-simplex of morphisms in $\Hanst$, as defined in \cite[Def 6.15]{GIMM1}. A routine check shows that this respects composition and the simplicial structure. Since our square gives an isomorphism between the simplicial sets of objects in $\Hanst$ over each object in $\Manst$, after passing from left fibrations to simplicial functors, it gives a weak equivalence of simplicial functors $\mc H_{G\dashiso}(M \amalg N) \simeq \mc H_{G\dashiso}(M) \times \mc H_{G\dashiso}(N)$.
	%From \cite[Def 6.15]{GIMM1}, these $p$-simplices of morphisms consist of $p$-simplices in $\Manst(M_0,M_1)$ and $\Manst(N_0,N_1)$ along with $h$-cobordisms over $M_0 \times \Delta^p$from $M_0$ to $M_1$ and from $N_0$ to $N_1$ to the corresponding 
%	consists of data   given by a pair of isovariant $h$-cobordisms $X$ on $M_0$ and $Y$ on $N_0$, and $p$-tuple of composable morphisms $M_0 \to \cdots \to M_p$ and $N_0 \to \cdots \to N_p$ in $\Manst$ We extend this to every stabilization in the obvious way. \cmar{Think about this more.}
\end{proof}

\section{Equivariant and isovariant Whitehead torsion}\label{sec:torsion}

In this section we recall and extend the theory of equivariant and isovariant Whitehead torsion. Equivariant Whitehead torsion is defined for all maps of finite $G$-CW complexes, but isovariant torsion is more restrictive, and in the literature is only defined for compact smooth $G$-manifolds. We recall these notions and extend the notion of isovariant torsion to compact smooth $G$-manifolds \emph{with facets}, a special kind of manifolds with corners where the boundary has slightly better behavior. We will use this in the next section to prove some key facts about $h$-cobordisms that are trivial on one side, which will lead to a key fiber sequence in the proof of the main theorem.

\subsection{Manifolds with facets and equivariant triangulations}

From this point onwards we make the following restriction for the manifolds $M$ that we consider. Recall from \autoref{df:face} that a face of $M$ is a union of components of $\ti\partial M$ that is $G$-invariant and such that the map back to $M$ is injective. We made that definition because it's the setting in which a $G$-invariant collar can be defined, but now we want to focus on the individual components of these faces.
\begin{defn}
	A \ourdefn{facet} $A \subseteq \ti\partial M$ is a single component with the property that $i\colon \ti\partial M \to M$ is injective on $A$. By abuse of notation we also use $A$ to refer to the corresponding subspace of $\partial M \subseteq M$.
\end{defn}

In general, a facet will not be $G$-invariant, because its orbit will consist of several facets. These may not form a face because the different facets may overlap in $M$, for instance if $M$ is a square and $G = C_4$ acts by rotation by $90^\circ$.

\begin{defn}
	A \ourdefn{compact $G$-manifold with facets} is a compact $G$-manifold with $G$-trivial corners $M$ such that $\partial M$ is a union of facets.
\end{defn}

The following lemma is a helpful consistency check that as we pass up to cobordisms or back down to their faces, we never leave the category of manifolds with facets.

\begin{lem}\label{facets_have_facets}
	If $M$ is a manifold with facets, then every face $F \subseteq \partial M$ is also a manifold with facets, as is every $h$-cobordism $W$ on $M$.
\end{lem}

\begin{proof}
	Since $F$ is a face, it has a collar by \autoref{collars_unique}, so without loss of generality $M = F \times I$. But then if $F$ failed to be a manifold with facets, it would have smooth boundary component $\ti A$ such that $\ti A \to F$ is not injective. But then $\ti A \times I$ would be a smooth boundary component of $M = F \times I$, whose map back to $M$ would not be injective. This would contradict the assumption that $M$ is a manifold with facets.
	
	For cobordisms, we first check the trivial cobordism $M \times I$, whose boundary consists of two faces $M \times \{0\}$ and $M \times \{1\}$ and the products $F \times I$ for each facet $F \subseteq \partial M$. Since $F \to M$ is injective, so is $F \times I \to M \times I$, so these are all facets as well, making $M \times I$ a manifold with facets.
	
	Now if $W$ is an arbitrary $h$-cobordism on $M$, then by assumption $\partial W$ consists of $M$ and $N$, which are assumed to be faces, and $\partial M \times I$, which is composed of facets by the previous analysis, so $W$ is a manifold with facets as well.
\end{proof}

Our goal in this first subsection is to deduce an equivariant triangulation theorem for compact $G$-manifolds with facets, from the corresponding theorem of Illman for closed $G$-manifolds \cite{illman_finite}.

To set this up, we first analyze how collars on the different facets should interact. Recall from \autoref{df:collar} that for any facet $F \subseteq M$, a \ourdefn{collar} is a smooth map $F \times [0,1] \to M$ that is depth-preserving on $F \times [0,1)$ (and is therefore an open embedding on that set), and that is the identity of $F$ on $F \times \{0\}$. The depth-preserving condition implies the following.
\begin{lem}
	For any point in the intersection of two facets $x\in F \cap F'$, any collar for $F$ will send $x \times [0,1]$ into $F'$.
\end{lem}

\begin{proof}
	If any point $x$ fails this condition, take the point $x \in \partial F$ for which the condition fails on $x \times [t,t+\epsilon)$ for all $\epsilon > 0$ and for the smallest possible value of $t \geq 0$. So when $s < t$, $c(s,t)$ is in all of the same facets as $x$ (except for $F$ when $s > 0$), but for infinitely many $s > t$ that are arbitrarily close to $t$, $c(s,t)$ is in some facet that $x$ was not in.
	
	If $c(x,t)$ fails to lie in one of the facets containing $x$, then one some interval $(t-\epsilon,t)$ it also fails to lie in that facet, so $t$ wasn't minimal. Therefore $c(x,t)$ is in all of the facets containing $x$, but there are infinitely many points $s \in (t,t+\epsilon)$ such that $c(x,s)$ is outside one of these facets, no matter how small $\epsilon$ is. In order for $c$ to preserve depth, each of these points $c(x,s)$ must also lie in some other facet $F''$ that does not contain $x$. Since there are only finitely many possible values for $F''$, we can take a sequence of points $s_i \to t$ in which each point $c(x,s_i)$ lies in the same facet $F''$ not containing $x$. Since $F''$ is closed, we conclude that $c(x,t) \in F''$ as well. But $c(x,t)$ is also in the same facets as $x$ (other than $F$), so the depth of $c(x,t)$ is greater than what is required by the collar condition (equal to the depth of $x$ if $t = 0$, or one less than the depth of $x$ when $t > 0$). This is a contradiction.
\end{proof}

\begin{defn}\label{df:collars_compatible}
	Given a compact $G$-manifold with facets $M$, and collars $c_i\colon F_i \times [0,1] \to M$ for each facet $F_i$, we say the collars are \ourdefn{compatible} if for each pair of facets $F_i$ and $F_j$, we have an equality of functions $c_i(c_j(x,t_j),t_i) = c_j(c_i(x,t_i),t_j)$ defined on $(x,t_i,t_j) \in (F_i \cap F_j) \times [0,1]^2$, 
%	\[ (F_i \cap F_j) \times [0,1]^2 \to M \]
%	\[ (x,t_i,t_j) \mapsto c_i(c_j(x,t_j),t_i), c_j(c_i(x,t_i),t_j) \]
and whenever $gF_i = F_j$, we have an equality of functions $gc_i(x,t) = c_j(gx,t)$ defined on $(x,t) \in F_i \times [0,1]$.
\end{defn}

Iteratively, the first condition implies that for any larger set of intersecting facets, $F_1$, $\ldots$, $F_k$, if we choose any $k$ numbers $t_1, \ldots, t_k \in [0,1]$ and apply the operations $c_i(-,t_i)$ successively to $x$ in any order, we always get the same point in $M$. So for instance
\[ c_1(c_2(c_3(x,t_3),t_2,t_1) = c_2(c_1(c_3(x,t_3),t_1),t_2). \]
Equivalently, there is a well-defined ``depth $k$ collar'' of the form
\[ c\colon (F_1 \cap \cdots \cap F_k) \times [0,1]^k \to M \]
that is obtained by applying the collars $c_i$ in any order, cf.~\cite[Prop 7.1]{qin_gluing}. Without loss of generality, we may as well take $k$ to be the depth of $x$, i.e. these are all of the facets that contain $x$. But the condition also holds for any set of facets containing $x$, even if it isn't maximal.

\begin{lem}\label{compatible_collars_exist}
	For any compact $G$-manifold with facets $M$, there exists a collection of collars $c_i\colon F_i \times [0,1] \to M$ for each facet $F_i$ that are compatible in the sense of \autoref{df:collars_compatible}.
\end{lem}

\begin{proof}
	It is enough to construct collar germs, i.e. depth-preserving embeddings $F_i \times [0,\epsilon) \to M$ that agree with $F_i$ on $F_i \times \{0\}$, with the same compatibility as above. Because then we can compose all of them with the same fixed smooth embedding $[0,1] \to [0,\epsilon)$ sending $0 \mapsto 0$, to get actual collars. For any open set $U \subseteq M$ containing the compact subset of the boundary $C \subseteq U \cap \partial M$, we will consider such germs that are defined on some open subset of $\partial M$ containing $C$ and landing in $U$.
	
	Near each point $x \in \partial M$ with depth $k$ and stabilizer $G_x \leq G$, we can choose any coordinate chart $U_\alpha \subseteq M$ containing $x$ identified with an open set $V_0 \times V_1 \subseteq V \times [0,\infty)^k$, where $V$ is a $G_x$-representation, $V_0 \subseteq V$ is a $G_x$-invariant convex open set that contains the origin, and $V_1 \subseteq [0,\infty)^k$ is a convex open set that contains $[0,1]^k$. Let $\partial^- [0,1]^k$ be the half of the boundary corresponding to the faces of the cube in which a coordinate is zero. Let $C_\alpha \subseteq U_\alpha \cap \partial M$ be a compact neighborhood of $x$ in $\partial M$ corresponding to a product $K_0 \times \partial^- [0,1]^k$, where $K_0 \subseteq V_0$ is a compact neighborhood of the origin.
	
	We consider the collection of all smooth maps $V_0 \times V_1 \to V_0 \times V_1$ that are defined in a neighborhood of $K_0 \times \partial^- [0,1]^k$, and that define depth $k$ collar germs on that neighborhood. A partially-defined smooth map $V_0 \times V_1 \to V_0 \times V_1$ is in this space iff it sends if it sends $V_0 \times (V_1 \cap \partial([0,\infty)^k)) \to V_0 \times V_1$ by the identity, and if for every point in this boundary of depth 1 whose $i$th coordinate is 0, the first derivative applied to the $i$th coordinate direction in $[0,\infty)^k$ gives a vector whose $i$th coordinate in $V_1 \subseteq [0,\infty)^k$ is positive. (At every point of higher depth, the first derivative is necessarily the identity.) These are preserved by convex combination in $V \times \R^k$, and therefore we can make new maps in this collection by taking convex combinations.
	
	As in the proof of \autoref{ntn_step1}, suppose we have a collection of $G$-equivariant collar germs $c_i\colon (F_i \cap U_\beta) \times [0,\epsilon) \to M$, defined on the $G$-invariant open set $U_\beta \cap \partial M$, and let $C_\beta \subseteq U_\beta \cap \partial M$ be any compact subset. Choose a smooth $G_x$-equivariant bump function $\lambda\colon M \to [0,1]$ that is 1 on a neighborhood of $C_\alpha \cong K_0 \times \partial^- [0,1]^k$ and 0 on a neighborhood of $U_\alpha^c$. Inside the chart $U_\alpha$, the given collars define a single function
	\[ c_\beta\colon V_0 \times V_1 \to V_0 \times V_1 \]
	defined on a neighborhood of $(\textup{supp }\lambda) \cap C_\beta$. (The neighborhood might be smaller than $U_\alpha \cap U_\beta \cap \partial M$ because the given collars might leave the open set $V_0 \times V_1$, on points in the domain that are sufficiently far away from the compact set $(\textup{supp }\lambda) \cap C_\beta$.)
	
	We define the new $G_x$-equivariant collar germs to be
	\[ c_{\alpha\beta} = \lambda(\id) + (1-\lambda)c_\beta \]
	where the linear combination is taken inside the single coordinate chart $U_\alpha$. Equivalently, we perform this interpolation to each collar $c_i\colon (F_i \cap U_\alpha) \times [0,1] \to V \times [0,\infty)^k$. As in the proof of \autoref{ntn_step1}, the interpolation is well-defined, smooth, and defines a collar germ
	\begin{itemize}
		\item in a neighborhood of $C_\alpha$ since $\lambda = 1$ there,
		\item in $(\textup{supp }\lambda)^c \cap U_\beta \supseteq (\textup{supp }\lambda)^c \cap C_\beta$ since $\lambda = 0$ there, and
		\item in a neighborhood of $(\textup{supp }\lambda) \cap C_\beta$ since both $\id$ and $c_\beta$ are defined there.
	\end{itemize}
	It is therefore well-defined and smooth on an open subset of $(U_\alpha \cup U_\beta) \cap \partial M$ that contains $C_\alpha \cup C_\beta$. We then extend this to a $G$-equivariant collar germ on the $G$-orbit of this open set, which is the union of the $G$-invariant open set $U_\beta$ and an open subset of $G \times_{G_x} U_\alpha$.
	
	Now as in the proof of \autoref{ntn_step1}, we cover $\partial M$ by the interiors of compact sets $G \times_{G_x} C_\alpha$ with this property, choose a finite list of them $G \times_{G_1} C_1$, $G \times_{G_2} C_2$, $\ldots$, whose interiors cover $\partial M$, and use the above procedure to construct a $G$-invariant collar germ defined on a neighborhood of $\bigcup_{i=1}^k G \times_{G_i} C_i$, by induction on $k$. After finitely many stages, this union is all of $\partial M$, giving our desired compatible $G$-equivariant collars.
\end{proof}

The following construction will allow us to efficiently prove claims about manifolds with facets by reducing them to claims about closed manifolds. Note that for a finite set $S$, the group $C_2^{\times S}$ can be identified with the power set $\mathcal{P}(S)$ with the operation on subsets given by symmetric difference, denoted $A\Delta B$.

\begin{defn}\label{hyperdouble}
	Let $M$ be a compact $G$-manifold with facets, and let $S$ be the set of facets, and for a point $x\in M$, let $S_x=\{F\in S\ |\ x\in F\}$. Define the \ourdefn{hyperdouble} $\overline{M}$, as a topological space, to be  the quotient of $C_2^{\times S} \times M$ by the relation $(A, x)\sim (B, x)$ when $A\Delta B\subseteq S_x$. Note that this equivalence relation is generated by $(A, x)\sim (A\Delta \{F\}, x)$ whenever $x\in F$. Informally, this means that when $x \in F$, the equivalence relation allows us to add or remove the element $F$ from the set $A$.
	
	%Thinking of the point in $C_2^{\times S} \times M$ as a subset $T \subseteq S$ and a point $x \in M$, any time $x$ lies in a face $F \in S$, the face $F$ can be added or removed from the set $T$.
	
	Define the semidirect product $C_2^{\times S} \rtimes G$ using the left action of $G$ on the set of facets $S$. So the multiplication is given by $(T_1,g_1)(T_2,g_2) = (T_1 \Delta  (g_1(T_2)), g_1g_2)$. This semidirect product has a well-defined continous action on $\overline{M}$ by the formula $(A,g) \cdot (B,x) = (A \Delta B, gx)$.
\end{defn}

Notice the following facts about $\overline{M}$:
\begin{itemize}
	\item There is a quotient map $\overline{M} \to M$ by forgetting the element of $C_2^{\times S}$, and this is $G$-equivariant.
	\item $\overline{M}$ contains $2^{\times S}$ embedded copies of $M$; they are embedded because the above equivalence relation never changes the point $x$.
	\item For every point $(A,x) \in \overline{M}$, the stabilizer group is always a semidirect product subgroup of the form $C_2^{\times S_x} \rtimes G_x$, where $S_x$ is the set of facets containing $x$ and $G_x$ is the stabilizer group of $x$ itself.
\end{itemize}

\begin{lem}\label{hyperdouble_top_mfld}\label{hyperdouble_smooth_mfld}
	$\overline{M}$ is a closed topological manifold. Furthermore, there exists a $C_2^{\times S} \rtimes G$-equivariant smooth structure on $\overline{M}$ such that $M \times \{e\}$ is a smoothly embedded $G$-invariant submanifold with corners.
\end{lem}

\begin{proof}
	Observe that if $x \in M$ has depth $k$, then it is in $k$ different facets at the same time. We form a neighborhood of $(A,x)$ by taking a neighborhood of $x \in M$ of the form $U \subseteq [0,\infty)^k \times \R^{n-k}$, then taking all subsets of $S$ that agree with $A$ except at these $k$ places, forming a subset $C_2^{\times k} \times U \subseteq C_2^{\times S} \times M$. Indexing the element of $C_2^{\times k}$ as a function $s\colon \{1,\ldots,k\} \to \{-1,1\}$, we map
	\[ C_2^{\times k} \times U \to \R^n \]
	\[ (s,y) \mapsto (s(1) \cdot y_1, \ldots, s(k) \cdot y_k, y_{k+1},\ldots, y_n). \]
	In other words, the function $s$ tells us which of the first $k$ coordinates to negate. This respects the equivalence relation because when $y_i = 0$, $(+1)y_i = (-1)y_i$. Furthermore it gives a bijection from a compact neighborhood of $(s,x) \in (C_2^{\times k} \times U)/{\sim}$ to a neighborhood of $x \in \R^n$, because the signs of the first $k$ coordinates determine which element of $C_2^{\times k}$ to pick, and once those are chosen the preimage $y \in U$ is uniquely determined. Therefore it gives a homeomorphism from a smaller open neighborhood of $(s,x) \in (C_2^{\times k} \times U)/{\sim}$ to a neighborhood of $x \in \R^n$.
	
	To make this into a smooth structure, we choose compatible equivariant collars as in \autoref{compatible_collars_exist}. Then near every depth $k$ point we only choose coordinate charts in which the first $k$ coordinates in $[0,\infty)^k$ come from the given collars. So at a point of depth $k$ with neighborhood $U \subseteq F_1 \cap \cdots \cap F_k$, we get a quotient map
	\[ C_2^{\times k} \times [0,\epsilon)^k \times U \to \R^k \times \R^{n-k} \]
	by $(s,y,u) \mapsto (s(1) \cdot y_1, \ldots, s(k) \cdot y_k, u)$. Since we are forcing ourselves to use these particular collars, the only choice we have is to change $U$ by a diffeomorphism, but that changes the coordinate chart by a diffeomorphism, so it gives a well-defined smooth structure near this point. It is also compatible with nearby points of depth $< k$ because when one of the first $k$ coordinates is positive, we can absorb it into the set $U$ and the rule does not change.
\end{proof}

This smooth structure is not canonical, as it depends on the choice of compatible collars. However we will see that it does not matter, any smooth structure such that $M$ is smoothly embedded will work for the following arguments.

Now we can move on to our main goal, which is to show that manifolds with facets have smooth equivariant triangulations, by reducing the claim to closed manifolds.

Recall that an equivariant simplicial complex $K$ is a $G$-space that is a simplicial complex, in which each element of $G$ acts by a simplicial homeomorphism. In other words, each $g \in G$ acts on each cell $\Delta^n$ either by fixing it or sending it to another cell by a linear homeomorphism. Equivalently, the complex is made of cells of the form $G/H \times \Delta^n$, attached linearly along faces.\footnote{In \cite{illman_finite} Illman imposes two more technical conditions on these complexes, but those extra conditions automatically hold after two barycentric subdivisions, so it is not necessary to state them.}

\begin{defn}\label{df:equivariant_triangulation}
	If $M$ is a compact $G$-manifold with facets, a \ourdefn{smooth equivariant triangulation} of $M$ is an equivariant homeomorphism $K \to M$, where $K$ is an equivariant simplicial complex, such that each of the resulting maps $\Delta^n \to M$ is a smooth embedding of manifolds with corners, and such that each facet of $M$ is identified with a subcomplex of $K$.
\end{defn}

\begin{thm}\label{equivariant_triangulation_unique}
	Every compact $G$-manifold with facets has an equivariant triangulation, which is unique up to subdivision and isotopy.
\end{thm}

\begin{proof}
	In the case of closed manifolds and manifolds with boundary, this is a theorem of Illman \cite{illman_finite}. We adapt his theorem to this case using the hyperdouble $\overline{M}$. The only thing to prove is that smooth $G$-equivariant triangulations of $M$ correspond to smooth $C_2^{\times S} \rtimes G$-equivariant triangulations of $\overline{M}$, preserving subdivision and isotopy.
	
	To go from a smooth $G$-equivariant triangulation of $M$ to a smooth $C_2^{\times S} \rtimes G$-equivariant triangulations of $\overline{M}$, we take each simplex $\sigma\colon \Delta^n \to M$, consider it as a simplex in $M \times \{e\} \subseteq \overline{M}$, and take its orbit under the action of $C_2^{\times S} \rtimes G$. If $S_\sigma$ is the set of facets that contain $\sigma(\Delta^n)$ and $G_\sigma$ is the stabilizer group of $\sigma$, then this orbit is of the form $(C_2^{\times S} \rtimes G) / (C_2^{\times S_\sigma} \rtimes G_\sigma) \times \Delta^n \to \overline{M}$. Thinking of these as individual simplices, they cover $\overline{M}$, and any overlap between two of them in one copy of $M$ is related by an element of $C_2^{\times S}$ to an overlap of two simplices in the original simplicial complex that covers $M$. Therefore the intersection must be a common face, and therefore we have defined a simplicial complex homeomorphic to $\overline{M}$. By construction, it defines a $C_2^{\times S} \rtimes G$-equivariant smooth triangulation of $\overline{M}$.
	
	To go the other way, we take any $C_2^{\times S} \rtimes G$-equivariant smooth triangulation of $\overline{M}$ and restrict it to $M \times \{e\} \subseteq \overline{M}$. This must give a subcomplex, since the interior of any simplex must have stabilizer group of the form $C_2^{\times S_\sigma} \rtimes G_\sigma$, and therefore must be entirely contained in one of the embedded copies of $M$, not straddled between two or more of them. This defines a $G$-equivariant smooth triangulation of $M$. The original triangulation of $\overline{M}$ must be constructed from this triangulation of $M$ since everything is in the image of the $C_2^{\times S} \rtimes G$-action, and the action of this group has to preserve the simplices.
\end{proof}

\subsection{Equivariant Whitehead torsion}
In proving that $\mc H_{\iso}^{V}(M)$ is equivalent to the usual $\mc H_\iso(M)$, we will use Whitehead torsion to give an indirect proof that two manifolds are equivariantly diffeomorphic. Therefore we begin by reviewing equivariant Whitehead torsion and isovariant Whitehead torsion. Our primary reference is \cite{luck}, though we will extend these notions to $G$-manifolds with facets, using the unique equivariant smooth triangulation from \autoref{equivariant_triangulation_unique}.

Recall that a pair $(X,Y)$ of $G$-spaces is a relative $G$-CW complex if $X$ has been obtained from $Y$ by repeatedly attaching collections of $G$-cells $G/H\times D^n$ along their boundaries $G/H\times S^{n-1}$, in order of increasing dimension. It is a finite relative $G$-CW complex if only finite many $G$-cells have been used. It is important here to actually fix the choice of cell complex structure, because the equivariant Whitehead torsion will depend on the choice that we use.

When $X$ is obtained from $Y$ by attaching $G/H \times D^n\times I $ along $G/H \times D^n \times \{0\}$, then the inclusion $Y\to X$ is called an  \emph{elementary expansion}.  A composition of finitely many elementary expansions is called an expansion. When $Y\to X$ is an expansion, then one says that $X$ collapses to $Y$. In this case $Y$ is an equivariant deformation retract of $X$ and in particular the map $Y\to X$ has an equivariant left inverse. A left inverse of an expansion is called a collapse. 
\begin{defn}
%An \emph{elementary expansion} is a map of finite $G$-CW complexes $X \to Y$ that attaches a single cell of the form $(G/H \times D^n \times \{0\}) \to (G/H \times D^n \times I)$.
 %An \emph{elementary collapse} is a map of finite $G$-CW complexes $Y \to X$ that is an equivariant homotopy inverse rel $X$ to an elementary expansion.
An \emph{equivariant simple homotopy equivalence} is an equivariant map of $G$-CW complexes that is $G$-homotopic to the composition of a finite sequence of expansions and collapses.
\end{defn}

\begin{defn}
	For any $G$-CW complex $A$, the \emph{equivariant Whitehead group} $\Wh^G(A)$ is the set of equivalence classes of finite relative $G$-CW complexes $(X, A)$ such that the inclusion $A \to X$ is an equivariant homotopy equivalence, with respect the following equivalence relation. Say that $(X,A) \sim (Y,A)$ if there is an equivariant simple homotopy equivalence $X \to Y$ relative to $A$, that is, a map that is a finite composition of expansions and collapses between finite $G$-CW complexes relative to $A$. (This is equivalent to saying that there exists a finite relative $G$-CW-complex $(Z,A)$ that collapses to both $(X,A)$ and $(Y,A)$.)
\end{defn}

It is shown in \cite[\S 4]{luck} that the rule $(X,A)+(Y,A)=(X\cup_A Y, A)$ gives a well defined abelian group structure on $\Wh^G(A)$. The inverse of $(X,A)$ is given by $(S_AX,A)$, where $S_AX$ is made by choosing a left inverse $X \to A$ for the inclusion and forming the fiberwise suspension 
$$
S_A X = (X \times I) \cup_{X \times 0\cup A \times I\cup X \times 1} A.
$$

%Fix a $G$-CW complex $A$ and consider the set of $G$-CW complexes $X$ that contain $A$ as a subcomplex, such that $A \to X$ is an equivariant homotopy equivalence.

%\begin{defn}
% The \emph{equivariant Whitehead group} $\Wh^G(A)$ is the set of such complexes $(X,A)$, up to the relation that $(X,A) \sim (Y,A)$ if there is an equivariant simple homotopy equivalence $X \to Y$ that is rel $A$, meaning the expansions, collapses, and homotopies all keep $A$ fixed.  Addition in $\Wh^G(A)$ is defined by pushout along $A$, $X \cup_A Y$. Negation is by picking a homotopy inverse $r\colon X \to A$ rel $A$, and using it to define fiberwise suspension $\Sigma_A X = (X \times I) \cup_{X \times 0, A \times I, X \times 1} A$.
% \end{defn}

\begin{exmp}
If $G$ is trivial then this is the usual Whitehead group $\Wh(A)$. When $A$ is path-connected this is a quotient of $K_1(\Z[\pi_1(A)])$, and in general it is direct sum of such groups over the path-components of $A$.
\end{exmp}

It is shown in \cite[\S 14]{luck} that there is a natural splitting 
\begin{equation}\label{wh_splitting}
\Wh^G(A) \cong \bigoplus_{(H)} \Wh(A^H_{hWH}),
\end{equation}
over conjugacy classes of subgroups of $G$. Here $\Wh(A^H_{hWH})$ is the non-equivariant Whitehead group for the space $A^H_{hWH}$. When $A$ is $G$-connected, i.e.~when $A^H$ is connected for all $H \leq G$, then each of the groups $\Wh(A^H_{hWH})$ is isomorphic to $K_1(\Z[\pi_1 A^H_{hWH} ])$. In \cite[\S 14]{luck} the splitting is proved using an algebraic description of the equivariant Whitehead group, but it could also be proved directly from the additivity property of torsion in \cite[Thm 4.8]{luck}.

\begin{exmp}\label{whitehead_free_case}
	As a special case, If $A$ is a free finite $G$-CW complex, any finite $G$-CW complex $X$ that is equivalent to $A$ must also be free. As in \autoref{cobordisms_on_free_manifold}, such complexes are determined by the corresponding quotient complexes $(X/G,A/G)$. So the above splitting has only one nonzero summand, the nonequivariant Whitehead group of the quotient:
\[ \Wh^G(A) \cong \Wh(A/G). \]
\end{exmp}

\begin{rem}In \cite{ccm}, the authors define a genuine equivariant algebraic $K$-theory $K_G(\Z[\underline{\pi_1(X)}])$, whose $G$-fixed points satisfy a tom Dieck style splitting $\prod_{(H)} K(\Z[\pi_1 X^H_{hWH} ])$, and there is an equivariant lift of the linearization map from the equivariant $A$-theory spectrum $A_G(X)$ from \cite{CaryMona}. As a corollary of the main result of this paper, the equivariant Whitehead group is the $\pi_0$ level of an equivariant Whitehead spectrum defined as the cofiber of the map $\Sigma_G^\infty X\to A_G(X)$. In \cite{ccm},  the authors define an algebraic equivariant Whitehead spectrum as the cofiber of the map $\Sigma_G^\infty X\to K_G(\Z[\underline{\pi_1(X)}])$, and they show that the map from the geometric to the algebraic equivariant Whitehead spectra is 2-connected, in particular recovering L\" uck's isomorphism of geometric and algebraic Whitehead groups.
\end{rem}

	The element of $\Wh^G(A)$ represented by $(X,A)$ is the \emph{equivariant Whitehead torsion}, denoted by $\tau^G(X,A)$, or by $\tau^G(i)$ where $i:A\to X$ is the inclusion. More generally, the torsion of any equivariant homotopy equivalence of finite $G$-CW complexes $f: X\to Y$ may be defined by replacing $f$ by a $G$-homotopic $G$-cellular map, replacing $Y$ by the mapping cylinder $Mf$, and taking the torsion of the inclusion $X \to Mf$.
	
	Any equivariant map of $G$-CW complexes $f\colon A \to B$ induces a map $f_*\colon \Wh^G(A) \to \Wh^G(B)$ by the rule $(X,A)\mapsto (X \cup_A B,B)$. 
			 
	\begin{prop}\cite[Thm 4.8]{luck}\label{equiv_torsion_properties} The equivariant Whitehead torsion has the following properties.
	\end{prop}
	\begin{itemize}
		\item Equivariantly homotopic maps $f \sim g$ have the same torsion $\tau^G(f) = \tau^G(g)$.
		\item $\tau^G(f) = 0$ iff $f$ is an equivariant simple homotopy equivalence.
		\item $f_*$ is a group homomorphism and respects composition, i.e. $\Wh^G$ is a functor to abelian groups.
		\item $f_*$ depends only on the equivariant homotopy class of $f$. Therefore if $f$ is an equivariant homotopy equivalence then  $f_*$ is an isomorphism.
		\item Given finite relative complexes $(X_j,A_j)$ with $A_j\to X_j$ a $G$-homotopy equivalence, and given maps
		\[ (X_1,A_1) \longleftarrow (X_0,A_0) \longrightarrow (X_2,A_2) \]
		of which at least one is an inclusion, we have
		 \[\tau^G(X_1 \cup_{X_0} X_2,A_1 \cup_{A_0} A_2) = i_{1\ast}\tau^G(X_1,A_1) - i_{0\ast}\tau^G(X_0,A_0)  + i_{2\ast}\tau^G(X_2,A_2), \]
		where $ i_j:A_j\to A_1\cup_{A_0}A_2$ is the inclusion.
		\item Composition of two equivalences $f\colon X \to Y$ and $g\colon Y \to Z$, satisfies the logarithmic rule
		\[ \tau^G(gf) = g_*\tau^G(f) + \tau^G(g). \]
	\end{itemize}

It may come as a surprise that equivariant Whitehead torsion is \emph{not} invariant under equivariant homeomorphism. So the choice of $G$-CW structure on the space matters. However, we do still have invariance under subdivision, because this can be reduced to the nonequivariant statement:

	\begin{lem}\label{torsion_subdivision}
		For any finite equivariant simplicial complex $K$ and any subdivision $K'$, the identity map $K \to K'$ has zero equivariant torsion.
	\end{lem}
	
	\begin{proof}
		It suffices to do this one cell of $K$ at a time. Suppose $K$ is some finite equivariant simplicial complex with subdivision $K'$, and we consider a new complex $L$ obtained from $K$ by adjoining a single equivariant simplex $G/H \times \Delta^n$ along its boundary $G/H \times \partial \Delta^n$. Furthermore, suppose $L'$ is a subdivision that restricts to the given subdivision $K'$ of $K$, and that on the new simplex is a subdivision that we call $(\Delta^n)'$. Then $(L',L)$ is the pushout of the three maps
		\[ (K',K) \longleftarrow (G/H \times (\partial\Delta^n)',G/H \times \partial\Delta^n) \longrightarrow (G/H \times (\Delta^n)',G/H \times \Delta^n). \]
		We compute that the equivariant torsion of the middle and right-hand pairs is zero. Along the splitting \eqref{wh_splitting} and canonical equivalence $(G/H \times X)^H_{hWH} \simeq X$, this becomes the claim that the non-equivariant torsion of the pairs $((\partial\Delta^n)',\partial\Delta^n)$ and $((\Delta^n)',\Delta^n)$ vanish. This in turn follows from the classical result that nonequivariant subdivisions have zero torsion \cite{milnor_torsion}, or the more powerful result that nonequivariant homeomorphisms have zero torsion \cite{chapman}.
		
		By the additivity formula in \autoref{equiv_torsion_properties}, if the equivariant torsion of $(K',K)$ vanishes then so does the equivariant torsion of $(L',L)$. Now by a simple induction the torsion always vanishes for any finite equivariant simplicial complex.
	\end{proof}
	
%\cmar{Commented out a paragraph that I don't think we need.}
%We point out some consequences. For any $G$-representation $V$ and any finite $G$-CW complex $A$, the inclusion $A\times 0\to A \times D(V)$ is the composition of elementary expansions, so that the projection $A\times D(V)\to A$ has zero torsion. Therefore any pushout of this map also has zero torsion. In particular, in any $G$-CW complex $X$, identifying a subcomplex of the form $A \times D(V)$ to $A$ induces a simple $G$-homotopy equivalence.

Since equivariant Whitehead torsion depends on the choice of $G$-CW structure, we need to take some care when adapting this definition to smooth $G$-manifolds, which do not have a canonical cell complex structure on them. We will use the following convention.
\begin{defn}
	If $M$ is any $G$-space, a \ourdefn{simple structure} on $M$ is an equivalence class of pairs $(X,f\colon X \to M)$, where $X$ is a finite $G$-CW complex and $f$ is a $G$-homotopy class that is an equivariant homotopy equivalence. Two of these pairs are equivalent if the equivarant homotopy class of the composition of one with the inverse of the other, $X \to M \to X'$, is equivariantly simple.
\end{defn}

Evidently $M$ has no simple structure unless it is equivariantly homotopy equivalent to some finite $G$-CW complex $X$, in which case the simple structures on $M$ are in bijection with $\Wh^G(X)$. If we fix a simple structure on each of the $G$-spaces $M$ and $N$, then we can define the torsion $\tau^G(f)$ for any $G$-homotopy equivalence $f\colon M \to N$, by composing with the simple structures to get a map of finite $G$-CW complexes. This torsion then satisfies the same properties as in \autoref{equiv_torsion_properties}, provided we consistently use the same simple structure for each space.

\begin{defn}
	For any smooth compact $G$-manifold with facets $M$, the \emph{canonical simple structure} on $M$ is the one represented by any smooth equivariant triangulation, in the sense of \autoref{df:equivariant_triangulation}.
\end{defn}

This is well-defined because any two smooth equivariant triangulations agree up to subdivision and isotopy, by \autoref{equivariant_triangulation_unique}. The isotopy does not change the torsion by \autoref{equiv_torsion_properties}, and the subdivision does not change the torsion by \autoref{torsion_subdivision}. Therefore every smooth compact $G$-manifold with facets $M$ does indeed have a canonical simple structure.

In conclusion, we can now speak unambiguously of the torsion $\tau^G(f)$ when $f$ is any $G$-homotopy equivalence of compact smooth $G$-manifolds with facets. Furthermore, by the uniqueness of this simple structure, every equivariant \emph{diffeomorphism} has $\tau^G(f) = 0$. (An equivariant \emph{homeomorphism} could still have nonzero equivariant torsion.)

\subsection{Isovariant Whitehead torsion}

Next we recall the isovariant Whitehead group and isovariant Whitehead torsion. In the literature these are not defined for $G$-spaces or $G$-CW complexes, only for $G$-manifolds. We extend the definition to $G$-manifolds with facets, so that the isovariant Whitehead torsion is defined for any isovariant $h$-cobordism $M \to W$ on a manifold with facets $M$, and more generally for any $G$-manifold with facets $W$ and a face (i.e.~a disjoint union of facets) $M$ such that the inclusion $M \to W$ is an isovariant homotopy equivalence.

%If $M$ is such a manifold then the group $\Wh^G_{\iso}(M)$ is defined, and if $M$ is a face of $W$ such that the inclusion $M \to W$ is an isovariant homotopy equivalence then the torsion $ \tau^G_{\iso}(W,M)$ is defined. In particular, it is defined if $W$ is an isovariant $h$-cobordism with bottom $M$.

%Now we define the isovariant Whitehead torsion. We only do this for a $G$-manifold with corners $W$, and a subset of the boundary $A \subseteq \partial W$ that is a union of faces, such that $A \to W$ that is an isovariant homotopy equivalence. Essentially, we are defining it for isovariant $h$-cobordisms, but we only need the equivalence on one side for the definition to make sense, and the sides don't need to be trivial.

\begin{defn}

		 The \emph{isovariant Whitehead group} of the $G$-manifold with facets $M$ is the direct sum\footnote{In \cite{luck} the space $M_H$ is used instead of $M_{[H]}$, but these are equivariantly homotopy equivalent by \autoref{h_brackets_equivalent}, so they have the same Whitehead group.}
		\[ \Wh^G_{\iso}(M) := \bigoplus_{(H)} \Wh^{WH}(M_{[H]}) \cong \bigoplus_{(H)} \Wh(M_{[H]}/WH). \] %\cong \bigoplus_{(H)} \Wh^1(M_{[H]}/WH). \]
		 If $W$ is any $G$-manifold with facets, and $M \to W$ is any face that is an isovariant homotopy equivalence, we define its \emph{isovariant torsion} as the sum
		\[ \tau^G_{\iso}(W,M) := \sum_{(H)} \tau^{WH}(W_{[H]},M_{[H]}). \]

\end{defn}

Some explanation is needed. We defined $M_{[H]}$ in \autoref{defn_strata} as the $H$-fixed points minus neat tubular neighborhoods of the lower strata of $M$. We define $W_{[H]}$ in exactly the same way, except that we arrange for the tubular neighborhoods in $W$ to restrict to give the chosen tubular neighborhoods in $M$. As a result, $W_{[H]}$ is a compact $G$-manifold with corners, and $M_{[H]}$ is a face. Most components of $\ti\partial W_{[H]}$ are contained in a component of $\ti\partial W$, and so the map to $W_{[H]}$ is still injective. The only other components come from the sphere bundle of a smoothly embedded submanifold, namely one of the fixed-point manifolds minus a neighborhood of the lower strata. The sphere bundle of this manifold is smoothly embedded into $W$, so the inclusion of this piece of $\ti\partial W_{[H]}$ is also injective. Therefore $W_{[H]}$ is a manifold with facets as well.

We then take the equivariant torsion $\tau^{WH}(W_{[H]},M_{[H]})$. Since the tubular neighborhoods we removed were neat, the $W_{[H]}$ we obtain is unique up to diffeomorphism by \autoref{neat_complement_well_defined}, so this torsion is well-defined. To make the definition more definite, we could take spherical blowup in the sense of \autoref{spherical_blowup} instead of removing the tubular neighborhoods.

\begin{exmp}
As a result of \autoref{whitehead_free_case}, when $M$ is any compact free $G$-manifold with facets we have
\[ \Wh^G_{\iso}(M) \cong \Wh(M/G). \]
\end{exmp}

\begin{rem}
	It is possible to prove that the inclusion of a zero section of an equivariant vector bundle, $M \to D(\nu)$, has zero equivariant torsion. However, it does not have zero isovariant torsion. In fact, its isovariant torsion is not even defined, because $M \to D(\nu)$ is not even an isovariant homotopy equivalence. Similarly, if we stabilize an isovariant $h$-cobordism $M \to W$, the inclusion $W \to \St^\nu(W)$ is not an isovariant homotopy equivalence, so it does not have a torsion. (We can deduce this from the previous claim by considering trivial cobordisms.)
	
	We do get an isovariant $h$-cobordism $M' \to \St^\nu(W)$ by \autoref{stabiso}, so its isovariant torsion is defined, but in general it does not agree with the isovariant torsion we started with,
	\[ \tau^G_\iso(\St^\nu(W),M') \neq i_*\tau^G_\iso(W,M), \]
	where $i_*\colon \Wh_\iso(M) \to \Wh_{\iso}(M')$ is the induced map of Whitehead groups. The issue is that the strata of the stabilization $St^\nu(W)$ are not the stabilizations of the strata of $W$. Only in the special case that the $G$-action on the fibers of $\nu$ is trivial, do the strata commute with stabilization, and in this case the cobordisms have the same isovariant torsion.
%	(relative to its bottom) as $W$, as in \autoref{stabilization_zero_torsion}. In fact, the groups they lie in are not even isomorphic in general: $\Wh_{\iso}(M') \not\cong \Wh_\iso(M)$.
\end{rem}

%\begin{rem}
%	The torsion $\tau^G_{\iso}(W,M)$ does \emph{not} agree in general with the torsion of the stabilization $\tau^G_{\iso}(\St^\nu(W),D(\nu))$, even though by \autoref{stabilization_zero_torsion} the map $W \to \St^\nu(W)$ is simple. The issue is that the strata of the stabilization $St^\nu(W)$ are not the stabilizations of the strata of $W$. However, in the special case that the $G$-action on the fibers of $\nu$ is trivial, the strata do commute with stabilization, and therefore the cobordisms have the same isovariant torsion.
%\end{rem}
\begin{thm}[Isovariant $s$-cobordism theorem]\label{luck_442}\cite{equivsurgery,rothenberg},\cite[Thm 4.42]{luck}
	Suppose for all $H \leq G$ either $M_H = \emptyset$ or $\dim M_H \geq 5$.
	\begin{enumerate}
		\item[(a)] If $W$ is an isovariant $h$-cobordism on $M$, the $\tau^G_{\iso}(W,M) = 0$ iff $W$ is diffeomorphic to a trivial $h$-cobordism.
		\item[(b)] Every element of $\Wh^G_{\iso}(M)$ is realized as $\tau^G_{\iso}(W,M)$ for some isovariant $h$-cobordism $W$ on $M$. This $W$ is unique up to diffeomorphism.
	\end{enumerate}
\end{thm}

The proof of (a) a simple induction up the strata, using the fact that a cobordism on a free $WH$-manifold is trivial iff its equivariant Whitehead torsion is trivial, using \autoref{cobordisms_on_free_manifold} and \autoref{whitehead_free_case}. Doing the lower strata first allows us to assume that each $W_{[H]}$ has been trivialized on the sides, which we need to know before we can conclude that $W_{[H]}$ is trivial.

The proof of (b) is by taking cobordisms representing each element $\tau^{WH}$ as a cobordism over $M_{[H]}$, stabilizing to make a cobordism on the tubular neighborhood of $M_{[H]}$, then gluing the results together. See \autoref{sec:splitting} for more details. Uniqueness follows by embedding one of the cobordisms into the lower collar of the other, and seeing that the difference has zero torsion and is therefore trivial.

Using \autoref{luck_442}, we can interpret the addition in $\Wh^G_{\iso}(M)$ by stacking: realize one element as an isovariant $h$-cobordism $W$ on $M$ with top $M'$, identify the groups $\Wh^G_{\iso}(M) \cong \Wh^G_{\iso}(M')$ along the equivalences $M \to W \leftarrow M'$, realize the second element as an isovariant $h$-cobordism $W'$ on $M'$, then define the sum to be the torsion of $(W' \cup_{M'} W,M)$.

In fact, more generally, given any element of $\Wh^G_{\iso}(M)$ represented by a manifold $W$ with face $M$, if $M'$ is any other face of $W$ disjoint from $M$, and we represent another element by a manifold $W'$ with $M'$ as a face, then we can define the sum of the torsions by the manifold $(W' \cup_{M'} W,M)$.

\begin{lem}\label{add_in_iso_wh_top}
	This coincides with the coordinate-wise addition in $\Wh^G_{\iso}(M) = \bigoplus_{(H)} \Wh^{WH}(M_{[H]})$.
\end{lem}

\begin{proof}
	It suffices to observe that we can choose tubular neighborhoods of the cobordisms and their faces consistently. This allows us to interpret
	\[ (W' \cup_M' W)_{[H]} \cong W'_{[H]} \cup_{M'_{[H]}} W_{[H]}. \]
	The torsion of this stacked cobordism agrees with the sum of the torsions in the equivariant Whitehead group $\Wh^{WH}(M_{[H]})$, by the logarithmic rule in \autoref{equiv_torsion_properties}.
\end{proof}

In the special case that the base manifold is a product $M \times I$ and the cobordism has one trivial side, we can also interpret addition by stacking side-to-side. Given two manifolds $W$ and $W'$, each of which has two adjacent faces of the form $M \times I$ joined along $M$, we can glue them along one of the faces to get a manifold $W \cup_{(M \times I)} W'$, which then gets a face of the form $(M \times I) \cup_M (M \times I) \cong (M \times I)$.

\begin{lem}\label{add_in_iso_wh_side}
	In this case we have
	\[ \tau^G_{\iso}(W \cup_{(M \times I)} W', (M \times I) \cup_M (M \times I)) = \tau^G_{\iso}(W, M \times I) + \tau^G_{\iso}(W, M \times I). \]
\end{lem}

\begin{proof}
	The construction $W \cup_{(M \times I)} W'$ commutes with passing to compactified strata $(-)_{[H]}$, so without loss of generality there is only one isotropy type and we can work with equivariant Whitehead torsion of free complexes. Then we collapse the $M \times I$ to $M$ to give a simple map
	\[ W \cup_{(M \times I)} W' \to W \cup_{M} W', \]
	showing that the torsion of $W \cup_{(M \times I)} W'$ relative to $M$ is the sum $\tau^G(W,M) + \tau^G(W',M)$. Then we note that the inclusions $M \to M \times I$ are simple, so this is identified with the sum of the torsions $\tau^G(W, M \times I) + \tau^G(W, M \times I)$.
\end{proof}

We can more generally consider the case that the base manifold $M$ has a face $A \subseteq \partial M$ and we are gluing on a cobordism $W'$ over $A \times I$ to the side of the cobordism, again stacking them side-to-side. Let $W$ be a manifold with two faces $M$ and $V$ joined along the common face $A$, and $W'$ a manifold with facets $A \times I$ and $V$ joined along the common face $A$. Gluing these together along $V$ gives the manifold $W \cup_{V} W'$, which then gets a face of the form $M \cup_A (A \times I) \cong M$.

\begin{lem}\label{add_in_iso_wh_side_2}
	In this case we have
	\[ \tau^G_{\iso}(W \cup_{V} W', M \cup_A (A \times I)) = \tau^G_{\iso}(W, M) + \tau^G_{\iso}(W',V). \]
\end{lem}

\begin{proof}
	As before we pass to compactified strata $(-)_{[H]}$, which commutes with the above gluing construction, and then work with the equivariant Whitehead torsion. Then we collapse $A \times I$ to $A$, and see that we have stacked two complexes along $V$, and so their equivariant torsions add.
\end{proof}

\section{Isovariant $h$-cobordisms rel $V$}

In this section we work with $h$-cobordisms on $M$ of a slightly more general kind: instead of requiring the sides of $W$ to constitute a trivial cobordism $\partial M \times I$, we allow one side of $W$ to be a nontrivial $h$-cobordism $V$, whose bottom is one face $A$ of $M$. (The remaining sides are still trivial).

For a fixed choice of $V$, we introduce a space $\mc H_{\iso}^{V}(M)$ of such $h$-cobordisms rel $V$. We show that under a suitable dimension hypothesis this is equivalent to the usual space $\mc H_\iso(M)$; intuitively this is because the side that is $V$ can be pulled up to the top, creating a trivial side $A \times I$ beneath it. We will spell out a more precise version of this argument where we glue on a ``basic $V$-block'' to go from $\mc H_\iso(M)$ to $\mc H_{\iso}^{V}(M)$, and then glue on an inverse block to go back. See \autoref{fig:v_block_and_inverse} for an illustration.

Finally, we use the equivalence between $\mc H_{\iso}^{V}(M)$ and $\mc H_\iso(M)$ to develop a fiber sequence relating $h$-cobordisms on $M$ to $h$-cobordisms on the fixed-point manifold $M^H$. This fiber sequence is the main ingredient in the proof of the splitting theorem in the next section.

\subsection{Relative $h$-cobordisms}\label{isovhcobrel}

Suppose that $M$ is a compact smooth $n$-dimensional $G$-manifold with facets and $A$ is an $(n-1)$-dimensional face of $M$. (So $A$ is a disjoint union of facets.) Suppose that $V$ is an $n$-dimensional equivariant $h$-cobordism on $A$; this is a manifold with facets by \autoref{facets_have_facets}. Denote the top of $V$ by $B$. By assumption, $V$ is equipped with a lower collar $c$, a partially-defined map $A \times I \to V$ that identifies a neighborhood of the bottom and sides of $A \times I$ with a neighborhood of the bottom and sides of $V$.

\begin{defn}
Given $M$, $A$, and $V$ as above, an \ourdefn{equivariant cobordism rel $V$} of $M$ is a compact $(n+1)$-dimensional $G$-manifold $W$ with facets, equipped with the following structure. There are two faces $M$ and $N$, called the bottom and top. There is a face that is identified with $V$, such that $V\cap M=A$ and $V\cap N=B$.
\end{defn}

By the \emph{sides} we mean the union of the remaining faces of $W$, that is, the closure of the complement of $M\cup V \cup N$ in $\partial W$. In the present context, the ``sides'' of the trivial cobordism $M \times [-1,0]$ means the faces other than $M \times \{-1\}$, $M \times \{0\}$, and $A \times [-1,1]$.

\begin{figure}[h]
	\centerline{\includegraphics[scale=1]{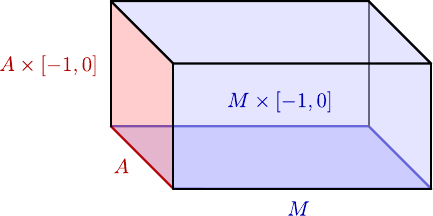}\hfill \includegraphics[scale=1]{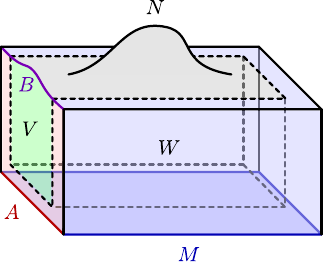}}
	\caption{A trivial $h$-cobordism rel $A\times I$ and a nontrivial $h$-cobordism rel $V$. They are identified near the bottom and sides.}
	\label{hcobside}
\end{figure}

We also require a \ourdefn{lower collar}. This now means an equivariant diffeomorphism $c$ between a neighborhood of the bottom and sides of $M\times [-1,0]$ and a neighborhood of the bottom and sides of $W$, taking the bottom to the bottom and the sides to the sides. %$M\times \{-1\}$ to the bottom of $W$, and $M\times \{0\}$ to the top of $W$. The remainder of the sides of $M \times [-1,0]$ therefore goes to the sides of $W$, and a neighborhood of $(\partial M \setminus A)\times \{0\}$ in $ M\times \{0\}$ goes to a neighborhood of $\partial N \setminus B$ in $N$. We refer to this embedding of a neighborhood of $M\times \{-1\}\cup ((\partial M \setminus \textup{int A})\times [-1,0])$ in $W$ as the and denote it by $c$. 
See \autoref{hcobside}.

\begin{defn}\label{hcob_def_rel}
If $V$ is an equivariant $h$-cobordism, an \ourdefn{equivariant $h$-cobordism on $M$ rel $V$} is a cobordism rel $V$ as above such that the inclusions $M \to W$ and $N \to W$ are equivariant homotopy equivalences. If $V$ is an isovariant $h$-cobordism on $A$, then an \ourdefn{isovariant $h$-cobordism on $M$ rel $V$} is a cobordism rel $V$ such that the inclusions $M \to W$ and $N \to W$ are isovariant homotopy equivalences.\end{defn}

%See \autoref{fig:h-cobordism}.

%\begin{figure}[h]
%	\centering
%	\includegraphics{h-cobordism.pdf}
%	\vspace{-.5em}
%	\caption{An $h$-cobordism on a manifold with corners $M$.}\label{fig:h-cobordism}
%\end{figure}

As before, we fix a sufficiently large set $\mc U$ containing $M$ and assume that the underlying set of each cobordism $W$ is a subset of $\mc U$. 

The \ourdefn{double} $\bar W$ of an $h$-cobordism $W$ rel $V$ is two copies of $W$ glued along their common top face $N$. This has $\bar V$, the double of $V$, as a face. %We think of one as ``flipped over'' and use the interval $[0,1]$ in the place of $[-1,0]$, with $M$ located at 1 and $N$ located at 0. 

\begin{figure}[h]
	\centerline{\includegraphics[scale=1]{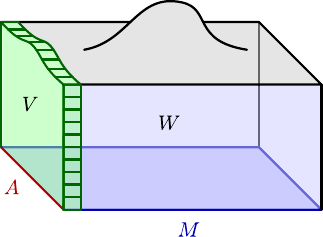}}
	\caption{A side collar on an $h$ cobordism rel $V$ }
	\label{sidecollar}
\end{figure}

\begin{defn}
We define some kinds of structure that a relative $h$-cobordism can have.
\begin{enumerate}
	\item A \ourdefn{side collar} is a collar $c_s\colon V \times I \to W$ that is compatible with the lower collar $c$ in the sense of \autoref{df:collars_compatible}. Spelling this out, $c_s(c(a,s),t) = c(c_s(a,t),s)$ for all $a \in A$ and $s,t \in I$ for which both sides are defined. See \autoref{sidecollar}.
	\item If $V$ has a mirror structure as in \autoref{def_retract_height}, a \ourdefn{\mirror} $h$-cobordism rel $V$ is a mirror structure on $W$ extending the existing one on the side collar $V \times I$. 
	\item If $W$ is mirror and $V$ has an encasing function as in \autoref{def_retract_height}, an \ourdefn{\rh} on $W$ is one that extends the given encasing function on $V \times I$. (The definition of the encasing function on $W$ is the same as in \autoref{def_retract_height}, using our new notion of ``sides'' that does not include $V$ as a side.)
\end{enumerate}
If $V$ is an encased isovariant $h$-cobordism, an \ourdefn{encased} isovariant $h$-cobordism rel $V$ is one equipped with a side collar and compatible mirror structure and \rhperiod As before, this additional structure is useful for stabilizing $h$-cobordisms.

\end{defn}

%\begin{figure}[h]
%	\centering
%	\includegraphics{encased.pdf}
%	\vspace{-.5em}
%	\caption{An encased $h$-cobordism.}\label{fig:encased}
%\end{figure}
As before, when making an encased $h$-cobordism it is redundant to specify the lower collar $c$, since it must coincide with the inverse of \therh  $\rho$ near the bottom and sides.

\begin{prop}\label{collars_contractible}
	The space of side collars is a contractible Kan complex.
\end{prop}

Note that the side collar is chosen before the mirror structure and \rh and not after.
\begin{proof}
	A side collar is defined by a vector field defined near $V$ that is inward-pointing along $V$, in the sense defined in \cite[Lem 2.31]{GIMM1}. Essentially, it has to preserve the stratification of $W$ by depth, except that the face $V$ is not counted when calculating the depth of a point.
	
	In each coordinate chart the space of such vector fields is preserved by affine linear combinations. We can therefore extend a  $\partial \Delta^p$ of such fields smoothly to a $\Delta^p$-family of such fields, by a variant of the proof of \autoref{neat_tubular_nbhd_thm} -- we make the extension on $\Delta^p$ minus one face, repeat for each face, then use a smooth partition of unity to add the extensions together. (We could alternatively smoothly extend the fields to an open neighborhood of $\partial\Delta^p$, take a constant field in the interior of $\Delta^p$, then smoothly interpolate between the resulting two fields.)
	
	Flowing along this $\Delta^p$-family of vector fields now gives the desired $\Delta^p$-family of collars.
\end{proof}

To define a space of isovariant $h$-cobordisms rel $V$, we consider families of $h$-cobordisms rel $V$ over $\Delta^k$. The definition is the expected one: a family is a compact $G$-manifold $E$ with facets, an equivariant bundle $p\colon E \to \Delta^k$, equivariant faces $\Delta^k \times M \to E$ and $\Delta^k \times V \to E$ intersecting at $\Delta^k \times A$, and a lower collar %$c\colon \Delta^k \times M \times I \to E$
 whose domain must, as in \cite[Def 4.8]{GIMM1} and \autoref{hcob_space}, be an open neighborhood of the bottom and sides, and not merely open on each fiber separately.

\begin{defn}
 The \ourdefn{space of isovariant $h$-cobordisms over $M$ rel $V$}  is the simplicial set $\mc H_{\iso,\sbt}^{V}(M)$ whose $k$ simplices are families of isovariant $h$-cobordisms rel $V$ over $\Delta^k$, \emph{with a side collar}. Similarly, we define the space of encased isovariant $h$-cobordisms $\mc H_{\iso,\sbt}^{V,c}(M)$. The face and degeneracy maps are the same as before.
\end{defn}

\begin{lem}\label{h_cob_space_rel_v_kan}
	The space $\mc H_{\iso,\sbt}^{V}(M)$ is a Kan complex.
\end{lem}

\begin{proof}
	As in \autoref{h_cob_space_kan}, this follows the same proof as in \cite[Lem 4.11]{GIMM1}.
\end{proof}

\begin{prop}\label{encased_forgetful_side}
	The forgetful map $\mc H_{\iso,\sbt}^{V,c}(M) \to \mc H_{\iso,\sbt}^V(M)$ that forgets the encasing function and mirror structure is an equivalence of Kan complexes.
\end{prop}

\begin{proof}
	For the encasing function the argument is as in \autoref{retracts_contractible}. For the mirror structure, as in \autoref{space_of_mirrorhcob} the argument is as in \cite{GIMM1}. %The side collar is handled by \autoref{collars_contractible}.
\end{proof}

To build our next example, we first introduce notation that extends \autoref{half_disc} by defining a ``quarter-disc'' for any vector bundle that contains a trivial plane bundle, in other words for any vector bundle of the form $\nu \times \R^2$.

\begin{notn}\label{quarter_disc}
	Building on \autoref{half_disc}, we let $QD(\nu \times \R^2)$ be the lower-right quarter-disc of the bundle $\nu \times \R^2$, and similarly let $QS$ be the quarter-sphere and $QC$ the quarter-complement:
	\begin{align*}
		QD(\nu \times \R^2) &= \{ (x,(t,u)) \in \nu \times \R^2 : \|x\|^2 + t^2 + u^2 \leq 1, \ t \geq 0, \ u \leq 0 \} \\
		QS(\nu \times \R^2) &= \{ (x,(t,u)) \in \nu \times \R^2 : \|x\|^2 + t^2 + u^2 = 1, \ t \geq 0, \ u \leq 0 \} \\
		QC(\nu \times \R^2) &= \{ (x,(t,u)) \in \nu \times \R^2 : \|x\|^2 + t^2 + u^2 \geq 1, \ \|x\| \leq 1, \ t \in [0,1], \ u \in [-1,0] \}
	\end{align*}
	In the case when $\nu$ is the trivial 0-dimensional vector bundle over a point, $QD(\R^2)$ is the lower-right quadrant of the unit disc in $\R^2$, $QS(\R^2)$ is the corresponding arc of the unit circle, and $QC(\R^2)$ is the closed complement of the quarter-disc in the quarter-square $[0,1] \times [-1,0]$.
\end{notn}

\begin{exmp}\label{basic_v_block}
	Let $V$ be an encased isovariant $h$-cobordism on $A$. We define an encased isovariant $h$-cobordism $B(V)$ on $A \times [0,1]$, that is $V$ on the face $A \times \{0\}$, as follows.
%	As in \autoref{half_disc}, let $HD(\R^2)$, $HS(\R^2)$, and $HC(\R^2)$ refer to the lower half of, respectively, the unit disc, the unit sphere, and the complement of the open unit disc in the square:
%	\begin{align*}
%		HD(\R^2) &= \{ (x,y) \in \R^2 : x^2 + y^2 \leq 1, \ y \leq 0 \} \\
%		HS(\R^2) &= \{ (x,y) \in \R^2 : x^2 + y^2 = 1, \ y \leq 0 \} \\
%		HC(\R^2) &= \{ (x,y) \in \R^2 : x^2 + y^2 \geq 1, \ x \in [-1,1], \ y \in [-1,0] \}
%	\end{align*}
%	Similarly we let $QD(\R^2)$, $QS(\R^2)$, and $QC(\R^2)$ refer to the right half of each of these regions, formed from a quarter of the disc:
%	\begin{align*}
%		QD(\R^2) &= \{ (x,y) \in \R^2 : x^2 + y^2 \leq 1, \ x \geq 0, \ y \leq 0 \} \\
%		QS(\R^2) &= \{ (x,y) \in \R^2 : x^2 + y^2 = 1, \ x \geq 0, \ y \leq 0 \} \\
%		QC(\R^2) &= \{ (x,y) \in \R^2 : x^2 + y^2 \geq 1, \ x \in [0,1], \ y \in [-1,0] \}
%	\end{align*}
%	We also let $HD(A \times \R^2) = A \times HD(\R^2)$ refer to the half-disc in the trivial bundle $A \times \R^2$, and similarly for all of the other notions.
	We first stabilize $\bar V$ along the embedding
	\[ A \to D(A \times \R) = A \times [-1,1], \]
	producing the pushout
	\[ \St(V) = \big( V \times_{A \times [0,1]} HD(A \times \R^2) \big) \cup_{HS(A \times \R^2)} HC(A \times \R^2). \]
	Then take $B(V)$ to be the right half, which is the pushout
	\[ B(V) = \big( V \times_{A \times [0,1]} QD(A \times \R^2) \big) \cup_{QS(A \times \R^2)} QC(A \times \R^2). \]
	This is an isovariant $h$-cobordism over $A \times [0,1]$ that is $V$ over the face $A \times \{0\}$. Forget its encasement, and take a side collar along $V$, then choose a new encasement that is compatible with this side collar and the lower collar.
	
\begin{figure}[h]
	\centerline{\includegraphics[scale=.75]{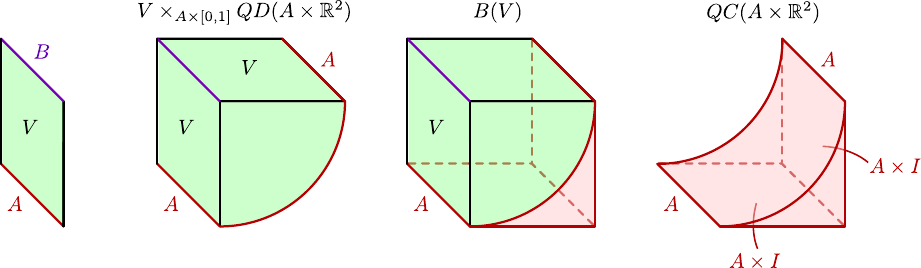}}\caption{A basic $V$-block $B(V)$ and its pieces.}
\end{figure}

	We call $B(V)$ the \ourdefn{basic $V$-block}. It is an isovariant $h$-cobordism rel $V$ in which $A \times I$ is the bottom and $V$ is the top (and another copy of $V$ is the only non-trivial side). %It is the simplest example of an $h$-cobordism rel $V$.
\end{exmp}

\begin{figure}[h]
	\centerline{\includegraphics[scale=0.6]{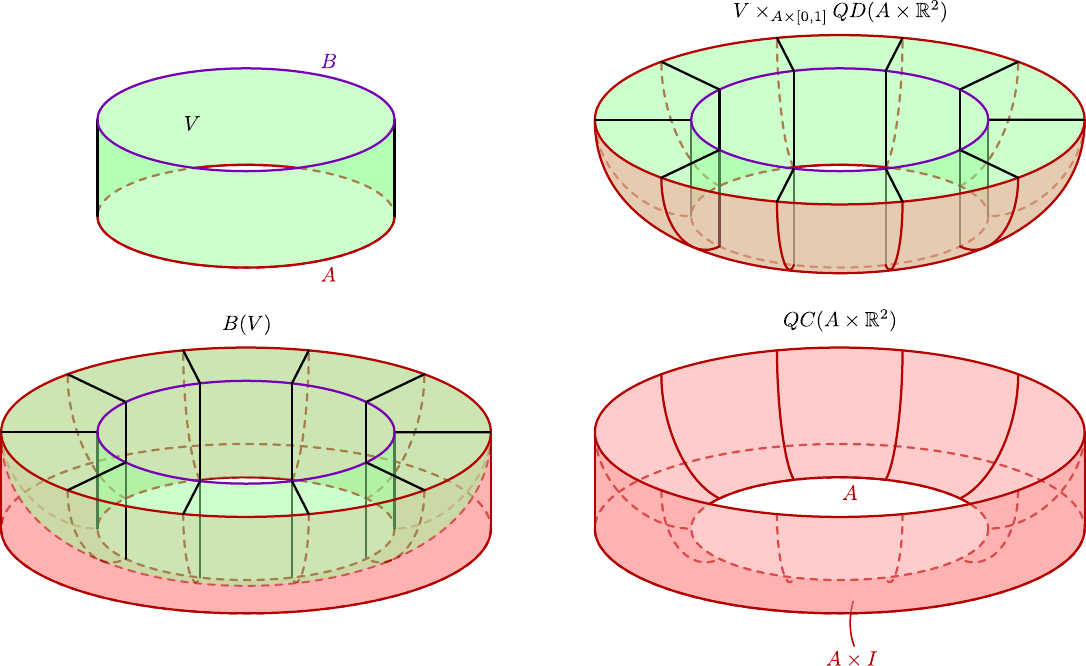}}\caption{A basic $V$-block when $A$ is a sphere.}
\end{figure}

\begin{rem}
	We believe that the isovariant torsion of $B(V)$ relative to each of its faces is:
	\[ \tau^G_{\iso}(B(V),V) = 0, \qquad \tau^G_{\iso}(B(V),A \times I) = \tau^G_{\iso}(V,A). \]
	However, we will not need these calculations in the proofs that follow. We only need to know that isovariant torsion is a well-defined concept, and that it adds when we stack manifolds together.
\end{rem}

%\begin{prop}\label{block_torsion}
%	The isovariant torsion of $B(V)$ relative to each of its faces is:
%	\[ \tau^G_{\iso}(B(V),V) = 0, \qquad \tau^G_{\iso}(B(V),A \times I) = \tau^G_{\iso}(V,A). \]
%\end{prop}
%
%\begin{proof}
%	Inspecting \autoref{hcob_stab} in the case when the normal bundle $\nu$ is a trivial bundle $M \times \R \to M$, with trivial $G$-action on $\R$, we find that
%	\[ B(V)_H \cong B(V_H), \quad B(V)_{[H]} \cong B(V_{[H]}). \]
%	Observe that the proof of \autoref{stabilization_zero_torsion} works just as well for this ``half-stabilization,'' showing that the inclusion $V_{[H]} \to B(V_{[H]})$ is $WH$-simple and therefore has zero equivariant torsion. %(This uses the fact that $G$ acts trivially in the stabilization direction.)
%	This implies both claims.
%\end{proof}

We will use this basic $V$-block to prove a version of the isovariant $s$-cobordism theorem (\autoref{luck_442}) that works rel $V$. We begin with the existence statement:
\begin{thm}\label{isovariant_s_cobordism_rel_V_existence}
	Suppose that $V$ is an isovariant $h$-cobordism on a face $A$ of $M$, and that for all $H \leq G$ either $M_H = \emptyset$ or $\dim M_H \geq 5$. Then every element of $\Wh^G_{\iso}(M)$ is realized as $\tau^G_{\iso}(W,M)$ for some isovariant $h$-cobordism $W$ on $M$ rel $V$.
\end{thm}

Note that by the logarithmic rule in \autoref{equiv_torsion_properties},
	\[ \tau^G_{\iso}(W,M) = \tau^G_{\iso}(W,M \cup_A V) + \tau^G_{\iso}(V,A) \]
so that we could just as well specify the torsion relative to $M \cup_A V$.

\begin{proof}
	Fix an element $\alpha \in \Wh^G_{\iso}(M)$. By \autoref{luck_442} there exists a cobordism $W'$ on $M$ with trivial sides such that
	\[ \tau^G_{\iso}(W',M) = \alpha - \tau^G_{\iso}(B(V),A \times I). \]
	Now take $W = B(V) \cup_{A \times I} W'$, whose base is $(A \times I) \cup_A M$, which is diffeomorphic to $M$. Along this diffeomorphism, $W$ is an $h$-cobordism on $M$, and by \autoref{add_in_iso_wh_side_2} its isovariant torsion is the desired value
	\begin{eqnarray*}
		\tau^G_{\iso}(W,M) &=& \tau^G_{\iso}(B(V),A \times I) + \tau^G_{\iso}(W',M) \\
		&=& \tau^G_{\iso}(B(V),A \times I) + (\alpha - \tau^G_{\iso}(B(V),A \times I)) \\
		&=& \alpha.
	\end{eqnarray*}
\end{proof}

\begin{prop}\label{block_inverse}
	If every nonempty stratum $A_H$ has dimension $\geq 5$, then the basic $V$-block $B(V)$ has an inverse $I(V)$, a cobordism over $A \times I$ that is $V$ on the face $A \times \{1\}$, such that
	\[ I(V) \cup_V B(V) \cong A \times I^2, \qquad B(V) \cup_{A \times I} I(V) \cong V \times I. \]
\end{prop}

\begin{figure}[h]
	\centerline{\includegraphics[scale=1.0]{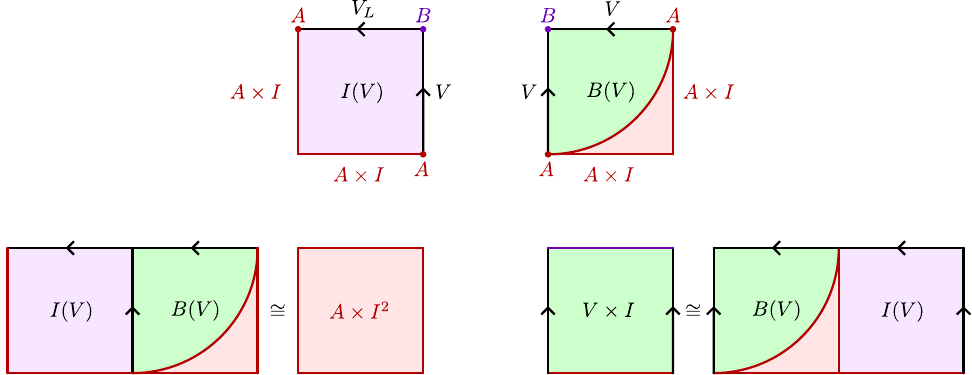}}\caption{The inverse of a basic $V$-block.}\label{fig:v_block_and_inverse}
\end{figure}

\begin{proof}
	It suffices to separately construct a left inverse $L(V)$ and a right inverse $R(V)$ such that
	\[ L(V) \cup_V B(V) \cong A \times I^2, \qquad B(V) \cup_{A \times I} R(V) \cong V \times I, \]
	because then it follows that
	\[ L(V) \cong L(V) \cup_V (V \times I) \cong L(V) \cup B(V) \cup R(V) \cong (A \times I^2) \cup_{A \times I} R(V) \cong R(V). \]
	
	Embed $V$ into $A \times I$ and let $V_L$ be its complement, so that $V_L \cup_B V \cong A \times I$. To construct $L(V)$ we choose any cobordism whose bottom is $A \times I$, whose right is $V$, whose top is $V_L$, and whose torsion rel $V$ is $-\tau^G_{\iso}(B(V),A \times I)$. (Such a cobordism exists by \autoref{isovariant_s_cobordism_rel_V_existence}, regarding $V$ as the bottom and $V_L$ as the nontrivial side.) By construction, the manifold $L(V) \cup_V B(V)$ has trivial bottom, top, and right side, and relative to the right side, its torsion is
	\[ \tau^G_{\iso}(B(V),A \times I) - \tau^G_{\iso}(B(V),A \times I) = 0. \]
	By \autoref{luck_442}, therefore $L(V) \cup_V B(V) \cong A \times I^2$.
	
%	first take a copy of $B(V)$ and flip it horizontally, so that its left edge is $(A \times I)$ and its right edge is $V$. Then we glue to its left an isovariant $h$-cobordism $W$ over $A \times I$ whose isovariant torsion is the negative of
%	\[ \tau^G_{\iso}(V,A) = \tau^G_{\iso}\St^{A \times \R}V. \]
%	For instance, we could take a cobordism $V'$ on $A$ whose torsion is $-\tau^G_{\iso}(V,A)$ and then stabilize it along $A \times \{1/2\} \to A \times I$. We need to know that $\dim A_H \geq 4$ in order to use the isovariant $s$-cobordism theorem (\autoref{luck_442}) to guarantee that such a cobordism $W$ exists. After gluing, the resulting cobordism $L(V) = W \cup_{A \times I} B(V)$ has the property that $L(V) \cup_V B(V)$ has zero isovariant torsion, using \autoref{add_in_iso_wh_side}. It is therefore diffeomorphic to the trivial cobordism $A \times I^2$, as desired.
	
	Now embed $V$ into $B \times I$ and let $V_R$ be its complement, so that $V \cup_A V_R \cong B \times I$. To construct $R(V)$ we choose a cobordism whose bottom is $A \times I$, whose left is $A \times I$, whose top is $V_R$, and whose torsion rel $V$ is $-\tau^G_{\iso}(B(V),V)$, again using \autoref{isovariant_s_cobordism_rel_V_existence}. By construction, the manifold $B(V) \cup_{A \times I} R(V)$ has trivial bottom and top, and relative to the left side $V$, its torsion is
	\[ \tau^G_{\iso}(B(V),V) - \tau^G_{\iso}(B(V),V) = 0. \]
	By \autoref{luck_442}, therefore $B(V) \cup_{A \times I} R(V) \cong V \times I$.
\end{proof}

\begin{rem}
	This result could be improved slightly to $\dim A_H \geq 4$, if the involution on the isovariant Whitehead group were a homomorphism. We know this by \cite{connolly_luck} under weak gap hypotheses on $M$, but not in general. If we did know it, we could define $I(V)$ to be $\St^{A \times \R}(V') \cup_{A \times I} B(V)$, and prove this is a two-sided inverse. When proving it is a right inverse, we get a cobordism that along the top is two copies of $V$ and two copies of $V'$, two of which are flipped. If the involution is a homomorphism then the torsion of this is trivial, so the cobordism is trivial along the top, and then as before we conclude that the entire thing is diffeomorphic to $V \times I$.
\end{rem}

\begin{cor}\label{isovariant_s_cobordism_rel_V_uniqueness}
	Suppose for all $H \leq G$ either $M_H = \emptyset$ or $\dim M_H \geq 6$. Then the isovariant $h$-cobordism $W$ rel $V$ representing a given element $\Wh^G_{\iso}(M)$ is unique up to equivariant diffeomorphism.
\end{cor}

\begin{proof}
	The dimension assumption on $M_H$ implies that $\dim A_H \geq 5$, so that \autoref{block_inverse} applies. Given two $h$-cobordisms $W_1$ and $W_2$ rel $V$ with the same isovariant torsion, gluing on the inverse $V$-block from \autoref{block_inverse} gives $h$-cobordisms with trivial sides that have the same torsion, by \autoref{add_in_iso_wh_side_2}. Therefore by \autoref{luck_442} these are diffeomorphic:
	\[ I(V) \cup_V W_1 \cong I(V) \cup_V W_2. \]
	Gluing $B(V)$ to each one and cancelling, we conclude that $W_1 \cong W_2$.
\end{proof}

\subsection{Equivalence of $\Hiso^{V}(M)$ to $\Hiso(M)$}\label{HequalsHV}

In this section we show that, if $M$ has a face $A$,  the dimensions $\dim A^H = \dim A_H$ for each $H \leq G$ are all at least 5, and $V$ is any isovariant $h$-cobordism on $A$, then the space of isovariant $h$-cobordisms on $M$ rel $V$ is equivalent to the ordinary space of isovariant $h$-cobordisms,
\[ \Hiso^{V}(M) \simeq \Hiso(M). \]
We have just shown that these two spaces have the same connected components in \autoref{isovariant_s_cobordism_rel_V_existence} and \autoref{isovariant_s_cobordism_rel_V_uniqueness}, and the idea will be to lift the argument to the rest of the space. We first consider the version in which there is no mirror structure or encasing function. However, $W$ is always assumed to be equipped with a side collar.

%Fix a neat collar neighborhood $A \times I \to M$ of $A$ in $M$, so that the complement is identified with $M$ again, giving a diffeomorphism $\phi\colon M \to (A \times I) \cup_A M$. Note that this diffeomorphism is isotopic as an embedding to the inclusion $i\colon M \to (A \times I) \cup_A M$. As the encased $h$-cobordism space $\Hiso^c(M)$ is an enriched functor by \autoref{functor}, the maps $h$ and $i$ induce homotopic maps. The same is therefore true on the spaces of ordinary $h$-cobordisms. (In fact, since $i$ and $h$ are codimension-0 embeddings, we don't need encasement to define the stabilization along these maps.)

The desired homotopy equivalence will be constructed from the following four maps.
\begin{defn}\label{plus_maps}
	Let $M$ be any compact $G$-manifold with facets, let $A$ be a face of $M$ (i.e.~a union of facets), and let $V$ be an isovariant $h$-cobordism on $A$. We define the maps
	\[ \xymatrix @C=4em @R=0.5em{
		\Hiso(M) \ar[r]^-{+B(V)} & \Hiso^{V}((A \times I) \cup_A M), \\
		\Hiso^{V}(M) \ar[r]^-{+I(V)} & \Hiso((A \times I) \cup_A M), \\
		\Hiso(M) \ar[r]^-{+(A \times I^2)} & \Hiso((A \times I) \cup_A M), \\
		\Hiso^{V}(M) \ar[r]^-{+(V \times I)} & \Hiso^{V}((A \times I) \cup_A M).
	} \]
	They are all defined the same way: glue the specified cobordism to the side of each cobordism $W$ over $M$, using the lower collar or side collar of $W$ to define the smooth structure on the result. Note that the map $+I(V)$ is only defined when $I(V)$ exists, which happens if every nonempty stratum $M_H$ has dimension $\geq 6$.
	
\begin{figure}[h]
	\centerline{\includegraphics[scale=1.0]{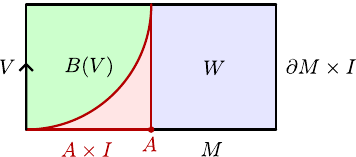} \hspace{2em} \includegraphics[scale=1.0]{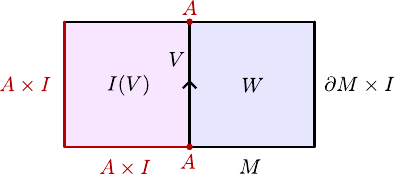}}
	\caption{The maps $+B(V)$ and $+I(V)$ send $W$ to these cobordisms.}
\end{figure}

\end{defn}

\begin{lem}
	The maps $+(A \times I^2)$ and $+(V \times I)$ from \autoref{plus_maps} are homotopy equivalences.
\end{lem}

\begin{proof}
	The map $+(A \times I^2)$ is another name for the map of $h$-cobordism spaces induced by the codimension zero embedding $\phi_0\colon M \to (A \times I) \cup_A M$. This embedding is homotopic to the map induced by the diffeomorphism $\phi_1\colon M \cong (A \times I) \cup_A M$ that stretches out a collar of $A$ in $M$ to cover $A \times I$, along an isotopy of codimension-zero embeddings $\phi_t$, each of which sends $A$ to the slice $A \times \{t\}$. As the encased $h$-cobordism space $\Hiso^c(M)$ is an enriched functor by \autoref{functor}, the maps $\phi_0$ and $\phi_1$ induce homotopic maps on the $h$-cobordism space. (In fact, since they are codimension-0 embeddings, we don't need encasement to define the stabilization along these maps.) But a diffeomorphism $\phi_1$ induces a homotopy equivalence, so the embedding $\phi_0$ must induce a homotopy equivalence as well. 
	
	For the second map $+(V \times I)$, we have to unpack the homotopy constructed implicitly in the first half of this proof and write it out explicitly, so that it will work even when the side is nontrivial. Each $\Delta^k$-family of cobordisms $E \to \Delta^k \times M$ is sent to a $\Delta^k \times \Delta^1$-family of cobordisms on $(A \times I) \cup_A M$, constructed for each $t \in \Delta^1$ by applying the codimension-zero embedding $\phi_t$ to get a cobordism on $(A \times [0,t]) \cup_{A \times \{0\}} M$ and then extending by $V \times [t,1]$. We make this a family of smooth manifolds over all of $\Delta^k \times \Delta^1$ by defining it on the image of the track of the isotopy $\phi\colon \Delta^1 \times M \to \Delta^1 \times ((A \times I) \cup_A M)$ as the pullback of $E \to \Delta^k \times M$ along the map $\id_{\Delta^k} \times \phi^{-1}$, and then taking the trivial cobordism $\Delta^k \times \Delta^1 \times V \times I$ on the rest, with smooth structure defined using the side collar of $E$. This gives a homotopy from the map $+(V \times I)$ to the map induced by the diffeomorphism $\phi_1$, and shows that $+(V \times I)$ is a homotopy equivalence.
\end{proof}

Now we define our homotopy inverses to be the following two zig-zags:
\[ \xymatrix @C=5em @R=0.8em{
	\Hiso(M) \ar[r]^-{+B(V)} & \Hiso^{V}((A \times I) \cup_A M) & \ar[l]^-\sim_-{+(V \times I)} \Hiso^{V}(M), \\
	\Hiso^{V}(M) \ar[r]^-{+I(V)} & \Hiso((A \times I) \cup_A M) & \ar[l]^-\sim_-{+(A \times I^2)} \Hiso(M).
} \]
These give well-defined maps in the homotopy category of spaces.

\begin{prop}\label{fan_away_v}
	If every nonempty stratum $M_H$ has dimension $\geq 6$, so that $I(V)$ exists, then these composites are homotopy inverses, and therefore $\Hiso^{V}(M) \simeq \Hiso(M)$.
\end{prop}

\begin{proof}
	The composite
	\[ \xymatrix @C=4em{
	\Hiso(M) \ar[r]^-{+B(V)} & \Hiso^{V}((A \times I) \cup_A M) \ar[r]^-{+I(V)} & \Hiso((A \times I) \cup_A (A \times I) \cup_A M),
	} \]
	sends each $h$-cobordism $W$ to
	\[ I(V) \cup_V B(V) \cup_{A \times I} W. \]
	By \autoref{block_inverse}, this is naturally diffeomorphic to $(A \times I^2) \cup_{A \times I} W$, and therefore this composite is homotopic to the one that attaches $A \times I^2$ twice to the side of $W$. Therefore the composite along the top and right of the following diagram is homotopic to the identity.
	\[ \xymatrix @C=4em{
	\Hiso(M) \ar[r]^-{+B(V)} & \Hiso^{V}((A \times I) \cup_A M) \ar[r]^-{+I(V)} & \Hiso((A \times I) \cup_A (A \times I) \cup_A M) \\
	& \ar[u]_-{+(V \times I)} \Hiso^{V}(M) \ar[r]^-{+I(V)} & \ar[u]_-{+(A \times I^2)} \Hiso((A \times I) \cup_A M) \\
	&& \ar[u]_-{+(A \times I^2)} \Hiso(M).
	} \]
	The middle square commutes up to homotopy because there is a diffeomorphism
	\[ I(V) \cup_V (V \times I) \cong I(V) \cong (A \times I^2) \cup_{A \times I} I(V) \]
	using the side collar and lower collar of $I(V)$. %(We could also use a combination of \autoref{add_in_iso_wh_top} and \autoref{luck_442}.) 
	Therefore the zig-zag along the lower-left is also homotopic to the identity, which gives the first half of the proof.
	
	The proof for the composite in the other direction is essentially the same. The composite
%	\[ \xymatrix @C=4em{
%	\Hiso^{V}(M) \ar[r]^-{+I(V)} & \Hiso((A \times I) \cup_A M) \ar[r]^-{+B(V)} & \Hiso^{V}((A \times I) \cup_A (A \times I) \cup_A M),
%	} \]
%	sends each $h$-cobordism $W$ to
%	\[ B(V) \cup_{A \times I} I(V) \cup_V W. \]
%	By \autoref{block_inverse}, this is diffeomorphic to $(V \times I) \cup_V W$, naturally in $W$. Therefore the composite
	 along the top and right of the following diagram is homotopic to the identity, and the square region commutes up to homotopy, so the zig-zag along the lower-left is also homotopic to the identity, finishing the proof.
	\[ \xymatrix @C=4em{
	\Hiso^{V}(M) \ar[r]^-{+I(V)} & \Hiso((A \times I) \cup_A M) \ar[r]^-{+B(V)} & \Hiso^{V}((A \times I) \cup_A (A \times I) \cup_A M) \\
	& \ar[u]_-{+(A \times I^2)} \Hiso(M) \ar[r]^-{+B(V)} & \ar[u]_-{+(V \times I)} \Hiso^{V}((A \times I) \cup_A M) \\
	&& \ar[u]_-{+(V \times I)} \Hiso^{V}(M).
	} \]
%	The square commutes up to homotopy because we can form a diffeomorphism
%	\[ B(V) \cup_{A \times I} (A \times I^2) \cong B(V) \cong (V \times I) \cup_V B(V). \]
%	Therefore the zig-zag along the lower-left is also homotopic to the identity, finishing the proof.
\end{proof}

The previous proof passes to encased $h$-cobordisms as well.
\begin{cor}\label{fan_away_v_encased}
	Under the same assumptions as in \autoref{fan_away_v}, the map $%(+(V \times I))^{-1} \circ 
	+B(V)$ also defines a homotopy equivalence on the spaces of encased $h$-cobordisms
	\[ \Hiso^c(M) \simeq \Hiso^{V,c}((A \times I) \cup_A M). \]
\end{cor}

\begin{proof}
It is straightforward to go from an encasement on $W$ to an encasement on $W \cup_{A \times I} B(V)$, so that we actually get a map of encased $h$-cobordism spaces. By \autoref{encased_forgetful_side} the encased space is equivalent to the unencased one. This gives a commuting square
	\[ \xymatrix @C=5em{
	\Hiso^c(M) \ar[d]_-\sim \ar[r]^-{+B(V)} & \Hiso^{V,c}((A \times I) \cup_A M) \ar[d]^-\sim \\ % & \ar[l]^-\sim_-{+(V \times I)} \Hiso^{V,c}(M) \ar[d]^-\sim \\
	\Hiso(M) \ar[r]^-{+B(V)}_-\sim & \Hiso^V((A \times I) \cup_A M) %& \ar[l]^-\sim_-{+(V \times I)} \Hiso^{V}(M) \\
	} \]
with the vertical maps equivalences. By \autoref{fan_away_v} the bottom map is an equivalence, so the top map is an equivalence as well.
\end{proof}

\subsection{The fiber sequence}\label{fiberseq}

Let $H \leq G$. Let $M$ be any compact $G$-manifold with facets, such that $M^K = \emptyset$ for all $K > H$. Equivalently,
\[ M^{H} \cap M^{H'} = \emptyset \]
for any nontrivial conjugate $H' = gHg^{-1} \neq H$. Let $M^{(H)} \subseteq M$ be the (disjoint) union of the fixed point submanifolds $M^{H'}$ as $H'$ ranges over all conjugates of $H$.

\begin{defn}\label{restriction_map}
	Under this assumption, we define the restriction map
	\[ \xymatrix{ \Hiso^c(M) \ar[r] & \Hiso^c(M^{(H)}) } \]
	by sending each cobordism $W$ to its $H$-conjugate fixed points $W^{(H)}$. We do the same to families of cobordisms indexed by $\Delta^k$, giving the required map of simplicial sets.
\end{defn}

Let $\nu$ denote the normal bundle of the neat embedding $M^{(H)} \subseteq M$. Fix a neat tubular neighborhood $D(\nu) \subseteq M$. As in \autoref{defn_strata}, let $M_{[<(H)]}$ denote the closure of the complement of $D(\nu)$ in $M$. We think of it as a compactification of $M \setminus M^{(H)}$ in $M$. Note that the unit sphere bundle $S(\nu)$ is a face of $D(\nu)$ (not the whole boundary) and is also a face of $M_{[<(H)]}$.

We will show that $\Hiso^c(M_{[<(H)]})$ is the homotopy fiber of the  the restriction map from \autoref{restriction_map}. First we compute the composite map
\begin{equation}\label{eq_fiber_composite}
	\xymatrix{ \Hiso^c(M_{[<(H)]}) \ar[r] & \Hiso^c(M) \ar[r] & \Hiso^c(M^{(H)})}.
\end{equation}
It takes each cobordism $W$ over $M_{[<(H)]}$ to its stabilization along the codimension 0 embedding $M_{[<(H)]} \to M$ (which just extends by the trivial cobordism) and then restricts to $M^{(H)}$. This always gives the same trivial cobordism over $M^{(H)}$. In other words the composite \eqref{eq_fiber_composite} is a constant map sending everything to the trivial cobordism. % \autoref{restriction_map} is the constant map sending any cobordism on $M_{[<(H)]}$  to the trivial cobordism on $M^{(H)}$.

Note that the space $\Hiso^c(M^{(H)})$ is not connected in general. Thus we also need to compute the fibers of the map $\Hiso^c(M) \to \Hiso^c(M^{(H)})$ over the other components. For the 0-simplex in $\Hiso^c(M^{(H)})$ corresponding to a cobordism $V$ on $M^{(H)}$, pick a strictly isovariant retract $V \to M^{(H)}$ (possible by \autoref{retracts_contractible}) and define
\begin{itemize}
	\item $D(U) = V \times_{M^{(H)}} D(\nu)$ is the fiber product of $V$ and $D(\nu)$ over $M^{(H)}$,
	\item $S(U) = V \times_{M^{(H)}} S(\nu)$ is the sides of $U$, a cobordism over $S(\nu)$.
\end{itemize}

\begin{defn}\label{glue_in_du}
As in \autoref{plus_maps}, we define the map
	\[ \xymatrix @C=4em @R=0.5em{
		\Hiso^{S(U),c}(M_{[<(H)]}) \ar[r]^-{+D(U)} & \Hiso^c(M)
	} \]
	by gluing the space $D(U)$ to the face $S(U)$ of each cobordism $W_{[<(H)]}$, 
%	as follows. A point $\Hiso^{S(U),c}(M_{[<(H)]})$ is a cobordism $W_{[<(H)]}$ on $M_{[<(H)]}$ with a nontrivial face over $S(\nu) \subseteq \partial M_{[<(H)]}$, that is identified with $S(U)$. We take such a cobordism to the union $W_{[<(H)]} \cup_{S(U)} D(U)$,
	using the side collar of $W_{[<(H)]}$ and a radial collar on $D(\nu)$ to define the smooth structure.
\end{defn}

\begin{thm}\label{fiber_seq_thm}
For each cobordism $V$ on $M^{(H)}$, the maps in \autoref{glue_in_du} and \autoref{restriction_map} define a fiber sequence
\[ \xymatrix @C=4em{ \Hiso^{S(U),c}(M_{[<(H)]}) \ar[r]^-{+D(U)} & \Hiso^c(M) \ar[r] & \Hiso^c(M^{(H)}). } \]
\end{thm}
Again, the composite of these maps is the constant map to the point represented by $V$, so this theorem is claiming that the map from $\Hiso^{S(U),c}(M_{[<(H)]})$ to the homotopy fiber of the restriction map is an equivalence. The proof of \autoref{fiber_seq_thm} occupies the rest of this section.

For simplicity, we use the version of $\Hiso$ with lower collars only, and the version of $\Hiso^{S(U)}$ with lower and side collars only, since this immediately implies the same theorem with encasements. (None of the maps in the fiber sequence require the encasements to define.)

\begin{rem}
	This theorem does not assume anything about the dimension of the strata $M_H$. The main theorem of the paper will also use \autoref{fan_away_v} however, so it will have a dimension restriction.
\end{rem}
	
\begin{lem}\label{fiber_seq_part1}
	The restriction map of \autoref{restriction_map} is a Kan fibration on the underlying simplicial sets.
\end{lem}

\begin{proof}
	\begin{equation}\label{fibration}
	\xymatrix{
		\Lambda^k_i \ar[d] \ar[r] & \mc H_{\iso,\sbt}(M) \ar[d] \\
		\Delta^k \ar[r] \ar@{-->}[ur] & \mc H_{\iso,\sbt}(M^{(H)})
		}
	\end{equation}
	Given a horn in $\mc H_{\iso,\sbt}(M)$ and an extension to a simplex in $\mc H_{\iso,\sbt}(M^{(H)})$, by the Ehresmann fibration theorem for $G$-manifolds with corners (\cite[Lem 2.34]{GIMM1}), the simplex in $\mc H_{\iso,\sbt}(M^{(H)})$ describes a family that is diffeomorphic to a trivial family $\Delta^k \times W^{(H)}$, and the horn describes a family diffeomorphic to $\Lambda^k_i \times W$, both for the same cobordism $W$ on $M$. Up to diffeomorphism, the desired lift is described then by the family $\Delta^k \times W$.
	
	To get the triangles to commute strictly, it is necessary to be careful about how we pick the underlying set of $\Delta^k \times W$. On the subspace $\Delta^k \times W^{(H)}$ we pick it to agree with the underlying set specified by $\Delta^k \to \mc H_{\iso,\sbt}(M^{(H)})$, and similarly on $\Lambda^k_i \times W$ we have it agree with the underlying set specified by $\Lambda^k_i \to \mc H_{\iso,\sbt}(M)$. These subspaces overlap, but the underlying sets agree on that overlap $\Lambda^k_i \times W^{(H)}$, by the commutativity of \eqref{fibration}.
\end{proof}

\begin{defn}\label{restriction_fiber}
	For a fixed cobordism $V$ on $M^{(H)}$, let $\mc H_{\iso,\sbt}(M)_V \subseteq \mc H_{\iso,\sbt}(M)$ be the subspace of cobordisms on $M$ whose $H$-fixed points are identically the cobordism $V$. This is the fiber of the restriction map $\mc H_{\iso,\sbt}(M) \to \mc H_{\iso,\sbt}(M^{(H)})$ at the point $V$. 
\end{defn}
To finish the proof of the fiber sequence we need to show:

\begin{lem}\label{fiber_seq_part2}
	The map $+D(U)$ induces an equivalence
	\[ \xymatrix{
		\mc H_{\iso,\sbt}^{S(U)}(M_{[<(H)]}) \ar[r]^-\sim & \mc H_{\iso,\sbt}(M)_V.
	} \]
\end{lem}

\begin{proof}
	Both spaces are Kan complexes by \autoref{h_cob_space_rel_v_kan} and \autoref{fiber_seq_part1}. So it suffices to take a square of the form
	\begin{equation}\label{fiber_seq_equivalence}
	\xymatrix{
		\partial \Delta^k \ar[d] \ar[r] & \mc H_{\iso,\sbt}^{S(U)}(M_{[<(H)]}) \ar[d] \\
		\Delta^k \ar[r] \ar@{-->}[ur] & \mc H_{\iso,\sbt}(M)_V
		}
	\end{equation}
	and modify the horizontal maps by a simplicial homotopy to a square in which a lift exists.
	
	The bottom map gives a $\Delta^k$-family of cobordisms on $M$ that on $M^{(H)}$ are identically $V$, and the top map gives a further identification of these cobordisms with $D(U)$ on $D(\nu)$, for every cobordism in the boundary family $\partial \Delta^k$. To give a lift is to give an identification of the cobordisms with $D(U)$ on $D(\nu)$, across the whole family. (The complement of this is then an isovariant $h$-cobordism in $\mc H_{\iso,\sbt}^{S(U)}(M_{[<(H)]})$ by \autoref{strata_criterion} because the inclusion of the top and bottom is an equivalence on each stratum. This is the part of the proof that would fail for equivariant $h$-cobordisms.)
	
	This can't be done strictly, but we are allowed to modify the maps up to homotopy. This corresponds to changing the family up to diffeomorphism.
	
	We do the lift in two stages. For each $0 < r \leq 1$ let $D_r(\nu) \subseteq D(\nu)$ be the smaller tubular neighborhood corresponding to the radius $r$ disc bundle inside $D(\nu)$. Let $D_r(U) \subseteq D(U)$ be the fiber product $D_r(\nu) \times_{M^{(H)}} V$. Define an intermediate space $\mc H_{\iso,\sbt}^{S_r(U)}(M_{[<(H)],r})$ of $h$-cobordisms that are identified with $D_r(U)$ over $D_r(\nu)$, not necessarily over all $D(\nu)$. The ``$r$'' subscript in $M_{[<(H)],r}$ indicates that we are using the version of the compactification where we remove the open disc bundle of radius $r$ around $M^{(H)}$. Now we construct lifts of the form
	\begin{equation}\label{fiber_seq_equivalence_2}
	\xymatrix @C=4em{
		\partial \Delta^k \ar[dd] \ar[r] & \mc H_{\iso,\sbt}^{S(U)}(M_{[<(H)]}) \ar[d] \\
		& \mc H_{\iso,\sbt}^{S_r(U)}(M_{[<(H)],r}) \ar[d] \\
		\Delta^k \ar[r] \ar@{-->}[ur] \ar@/^1.5em/@{-->}[uur] & \mc H_{\iso,\sbt}(M)_V.
		}
	\end{equation}
	
	To do the bottom lift for some value of $r \in (0,1]$, we observe that the data gives a $\partial \Delta^k$s worth of neat tubular neighborhoods of $V$ in the cobordisms $W$. By \autoref{neat_tubular_nbhd_thm} this can be extended to a $\Delta^k$s worth of germs of neat tubular neighborhoods. This means that we have an identification of the cobordisms with $D_r(U)$ near $V$. Changing the family up to diffeomorphism by relabeling the underlying sets, we can assume that in the whole family the cobordisms have a subset that is identical to $D_r(U)$, whose intersection with the bottom is $D_r(\nu)$ (and the original family has not been changed along $\partial \Delta^k$). This gives the first lift.
	
	Given this, to get the top lift, we pull back along an isotopy of diffeomorphisms $M \to M$ that expands $D_r(\nu)$ to $D(\nu)$. This induces a homotopy of the horizontal maps while preserving commutativity. At the end of the homotopy, the map now lands in the subspace $\mc H_{\iso,\sbt}^{S(U)}(M_{[<(H)]})$ so it has a unique lift making the diagram commute.
\end{proof}

\begin{rem}
	\autoref{fiber_seq_part2} fails for equivariant $h$-cobordisms. The problem is when we don't have an isovariant $h$-cobordism, deleting fixed point subspaces may produce a cobordism that is no longer homotopy equivalent to the base.
\end{rem}

\begin{rem}
	There is a similar fiber sequence
	\[ \xymatrix{ \mc H(Q) \ar[r] & \mc H(M,P) \ar[r] & \mc H(P) } \]
	for any neatly embedded $G$-manifold $P \subseteq M$, where $Q$ is the complement of an open tubular neighborhood of $P$, and $\mc H(M,P)$ is defined so that each point is a pair of an $h$-cobordism on $P$, and an extension to a cobordism on $M$, with the additional property that it is also an $h$-cobordism on the compactification of $M \setminus P$.
	%\cmar{Is this true? How do we know we still have an $h$-cobordism after deleting the cobordism over $P$?}\mmar{In Tom's setting, $P$ was a codimension 0 part of the boundary}
\end{rem}

\section{The splitting theorem}\label{sec:splitting}

In this concluding section we prove the main theorems, using the fiber sequence from the previous section for the inductive step. We first prove the unstable splitting theorem for the space of isovariant $h$-cobordisms, and then use the functoriality of $\mc H_\iso^c(-)$ to deduce a splitting theorem for the stabilized space of isovariant (or equivalently equivariant) $h$-cobordisms. Finally, we extend the splitting to be natural with respect to the manifold $M$ as well.

\subsection{Unstable splitting}\label{unstablesplitting}

Suppose $M$ is a compact $G$-manifold with $G$-trivial corners. Take a total ordering on the conjugacy classes $(H_0) = \{e\} < (H_1) < \ldots < (H_k) = G$ extending the subconjugacy partial ordering.

Recall the definition of compactified strata $M_{[(H_j)]}$ and $M_{[\leq j]}$ introduced in \autoref{defn_strata}. By \autoref{cobordisms_on_free_manifold}, we have a homotopy equivalence
\[ \mc H_{G\dashiso}^c(M_{[(H)]}) \simeq \mc H^c(M_{[(H)]}/G) \simeq \mc H^c(M_{[H]}/WH). \]
%The second map is a homeomorphism because we have an identification of the spaces themselves, $M_{[(H)]}/G \cong M_{[H]}/WH$.
So each of these compactified strata of $M$ gives a non-equivariant space of $h$-cobordisms.

Since the submanifolds $M_{[(H_i)]}$, $1 \leq i \leq k$, are disjoint from each other, we have a smooth embedding of $G$-manifolds
\begin{equation*}\label{embed_the_strata}
\coprod_{i=1}^k M_{[(H_i)]} \ra M
\end{equation*}
which induces a map
\begin{equation}\label{left_splitting}
	\prod_{i=1}^k \Hiso^c(M_{[(H_i)]}) \cong \Hiso^c\left( \coprod_{i=1}^k M_{[(H_i)]}\right)  \ra \Hiso^c(M).
\end{equation}
%\[
%	\left[ \prod_{i=j}^k \Hiso(M_{[(H_i)]}) \right] \times \Hiso(M_{[\leq j-1]}) \to \left[ \prod_{i=j+1}^k \Hiso(M_{[(H_i)]}) \right] \times \Hiso(M_{[\leq j]}).
%\]
%Each of these maps is peeling off one copy of $\Hiso(M_{[(H_j)]})$ from the product and combining it into the larger $h$-cobordism space:
%\begin{equation}\label{left_splitting_one_piece}
%	\Hiso(M_{[(H_j)]}) \times \Hiso(M_{[\leq j-1]}) \to \Hiso(M_{[\leq j]}).
%\end{equation}
%This move, in turn, takes a $G$-$h$-cobordism $W_{[(H_j)]}$ over the manifolds $M_{[(H_j)]}$, uses a $U$-shaped stabilization along the normal bundle of $M_{[(H_j)]}$ into $M_{[\leq j]}$, then glues the result to the existing cobordism on $M_{[\leq j-1]}$ to get a cobordism on all of $M_{[\leq j]}$. The smooth structure at the boundary of $M_{[\leq j-1]}$ is not a problem since both cobordisms have germs of lower collars to glue together.
%
%When combined all together, the map \eqref{left_splitting} takes cobordisms on each of the compactified strata and stabilizes them to be cobordisms on codimension 0 submanifolds of $M$. These codimension 0 submanifolds fit together like puzzle pieces to make all of $M$, and the stabilized cobordisms on each one are trivial at the boundaries, so we can glue all of them together to make a single cobordism on all of $M$.

\begin{thm}\label{unstable_splitting_thm}
	Assume that whenever we have a proper containment of components of fixed-point submanifolds
	\[ \emptyset \neq D \subsetneq C, \qquad D \in \pi_0(M^K), \qquad C \in \pi_0(M^H), \]
	then $\dim C \geq 6$. Then the map \eqref{left_splitting} is an equivalence. With no assumptions on dimension, the map is an equivalence when restricted to the connected component of the trivial $h$-cobordism in each space.
\end{thm}

Note that the assumption we make about the dimensions is automatically satisfied when every nonempty fixed point submanifold $M^H$ has dimension $\geq 5$, because then any proper containment would force the dimension of $M^H$ to be $\geq 6$.

\begin{proof}
The map \eqref{left_splitting} is composed of the intermediate maps
\begin{equation}\label{left_splitting_one_piece}
	\Hiso^c(M_{[(H_j)]}) \times \Hiso^c(M_{[\leq j-1]}) \to \Hiso^c(M_{[\leq j]})
\end{equation}
induced by the embeddings $M_{[(H_j)]} \amalg M_{[\leq j-1]} \to M_{[\leq j]}$. It suffices to show that each of these is an equivalence.

This is trivially true when $M_{[(H_j)]} \amalg M_{[\leq j-1]} \to M_{[\leq j]}$ is a diffeomorphism, in other words when $M_{[(H_j)]} = M_{[\leq j]}^{(H_j)}$ is isolated from all points whose stabilizer subgroup is smaller than $(H_j)$ in the total ordering. So without loss of generality we may focus on the components of $M_{[\leq j]}$ in which $M_{[(H_j)]}$ is not isolated, so that by assumption the larger manifold $M_{[\leq j]}$ has dimension at least 6. This implies that the codimension-zero submanifold $M_{[\leq j-1]}$ also has dimension at least 6.

By \autoref{fiber_seq_thm} applied to the $G$-manifold $M_{[\leq j]}$, the subgroup $H_j \leq G$, and the trivial cobordism on $M_{[(H_j)]}$, we get a fiber sequence
\[ \xymatrix{ \Hiso^c(M_{[\leq j-1]}) \ar[r] & \Hiso^c(M_{[\leq j]}) \ar[r] & \Hiso^c(M_{[(H_j)]}). } \]
We can form a map of fiber sequences
\[ \xymatrix{
	\Hiso^c(M_{[\leq j-1]}) \ar@{=}[d] \ar[r] & \Hiso^c(M_{[(H_j)]}) \times \Hiso^c(M_{[\leq j-1]}) \ar[d]^-{\eqref{left_splitting_one_piece}} \ar[r] & \Hiso^c(M_{[(H_j)]}) \ar@{=}[d] \\
	\Hiso^c(M_{[\leq j-1]}) \ar[r] & \Hiso^c(M_{[\leq j]}) \ar[r] & \Hiso^c(M_{[(H_j)]})
} \]
where the middle vertical is \eqref{left_splitting_one_piece}, the map we want to show is an equivalence. The left-hand square commutes up to homotopy since both branches extend a cobordism on $M_{[\leq j-1]}$ by a trivial cobordism to $M_{[\leq j]}$. The right-hand square commutes up to homotopy because stabilizing a cobordism from $M_{[(H_j)]}$ to its neighborhood and then restricting back to the $(H_j)$-fixed points gives the same cobordism back. (Intuitively, it is restricting a ``fan'' to the blade at the center). We conclude that \eqref{left_splitting_one_piece} is an equivalence on those components that map to the component of the trivial cobordism in $\Hiso^c(M_{[(H_j)]})$.

To handle the other components, we form a more elaborate map of fiber sequences where the basepoint on the right-hand side is represented by some cobordism $V$ over $M_{[(H_j)]} = M_{[\leq j]}^{(H_j)}$:
\[ \xymatrix{
	\Hiso^c(M_{[\leq j-1]}) \ar[d]_-{+B(S(U))}^-{(\textup{\autoref{plus_maps}})} \ar[r]^-{(V,\id)} & \Hiso^c(M_{[(H_j)]}) \times \Hiso^c(M_{[\leq j-1]}) \ar[d]^-{\eqref{left_splitting_one_piece}} \ar[r] & \Hiso^c(M_{[(H_j)]}) \ar@{=}[d] \\
	\Hiso^{S(U),c}((S(U) \times I) \cup_{S(U)} M_{[\leq j-1]}) \ar[r]_-{(\textup{\autoref{glue_in_du}})}^-{+D(U)} & \Hiso^c(M_{[\leq j]}) \ar[r] & \Hiso^c(M_{[(H_j)]})
} \]
Here the right-hand square is the same as before, and the map $(V,\id)$ is the constant map to the point $V \in \Hiso^c(M_{[(H_j)]})$, times the identity of $\Hiso^c(M_{[\leq j-1]})$. The bottom row is the fiber sequence of \autoref{fiber_seq_thm}, except that we use $(S(U) \times I) \cup_{S(U)} M_{[\leq j-1]}$ in the place of $M_{[\leq j-1]}$ as the model for $M_{[\leq j]}$ minus a tubular neighborhood of $M_{[(H_j)]}$.

To get the left-hand square to commute up to homotopy, we need to show that for cobordisms $W$ on $M_{[\leq j-1]}$ with trivial sides, there is a natural diffeomorphism
\[ D(U) \cup_{S(U)} B(S(U)) \cup_{S(\nu) \times I} W \cong \St^\nu(V) \cup_{S(\nu) \times I} W. \]
In other words, we want a single diffeomorphism of $h$-cobordisms
\[ D(U) \cup_{S(U)} B(S(U)) \cong \St^\nu(V). \]
%where $D(U) = V \times_{M^{(H)}} D(\nu)$.
We prove this in \autoref{last_diffeomorphism} below. (We also need to know that this diffeomorphism is the identity on the copy of $V$ at the center of $\St^\nu(V)$ and $D(U)$. This guarantees that the homotopies making the above squares commute actually give a map of fiber sequences: they are compatible with the identification between each composite of horizontal maps and the constant map sending everything to $V$. This compatibility, in turn, is what makes the map of fiber sequences commute with the connecting homomorphism in the long exact sequence of homotopy groups.)

This map of fiber sequences is an equivalence in the right-hand column because it is an identity map. In the left-hand column, it is an equivalence by \autoref{fan_away_v_encased}, since the dimension of $M_{[\leq j-1]}$ is at least 6. Therefore \eqref{left_splitting_one_piece} is an equivalence on these components as well.

In conclusion, \eqref{left_splitting_one_piece} is an equivalence when restricted to the preimage of each component of $\Hiso^c(M_{[(H_j)]})$, and is therefore an equivalence.
\end{proof}

\begin{lem}\label{last_diffeomorphism}
	For any equivariant $h$-cobordism $V$ over a manifold $M$ and any vector bundle $\nu \to M$, there is a diffeomorphism
	\[ (V \times_M D(\nu)) \cup_{V \times_M S(\nu)} B(V \times_M S(\nu)) \cong \St^\nu(V) \]
	that is the identity on the copy of $V$ at the center.
\end{lem}

\begin{figure}[h]
	\centerline{\includegraphics[scale=1.0]{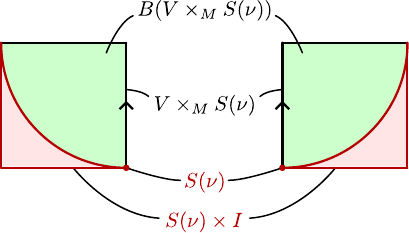}}%\caption{The inverse of a basic $V$-block.}
\end{figure}

\begin{figure}[h]
	\centerline{\includegraphics[scale=1.0]{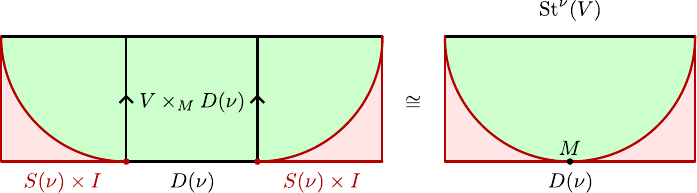}}%\caption{The inverse of a basic $V$-block.}
\end{figure}

\begin{proof}
	The copy of $V$ at the center of $\St^\nu(V)$ has normal bundle that is (equivariantly) isomorphic to the pullback of $\nu$ from $M$ to $V$, since $M \subseteq V$ is an equivariant homotopy equivalence by assumption. Therefore any neat tubular neighborhood of $V$ in $\St^\nu(V)$ is (equivariantly) diffeomorphic to $V \times_M D(\nu)$.
	
	In light of \autoref{blowup_removes_tubular_neighborhood}, it therefore suffices to define a spherical blowdown map $f\colon B(V \times_M S(\nu)) \to \St^\nu(V)$, whose exceptional locus is the face $V \times_M S(\nu)$ on the side of the block $B(V \times_M S(\nu))$, projecting down to the exceptional image $V \subseteq \St^\nu(V)$.
	
	We expand out the definition of $B(V \times_M S(\nu))$ and $\St^\nu(V)$ as
	\begin{align*}
		B(V \times_M S(\nu)) &= \big( (V \times_M S(\nu)) \times_{S(\nu) \times [0,1]} QD(S(\nu) \times \R^2) \big) \cup_{QS(S(\nu) \times \R^2)} QC(S(\nu) \times \R^2) \\
		&= \big( (V \times_M S(\nu)) \times_{[0,1]} QD(\R^2) \big) \cup_{S(\nu) \times QS(\R^2)} (S(\nu) \times QC(\R^2)) \\
		\St^\nu(V) &= \big( V \times_{M \times [0,1]} HD(\nu \times \R) \big) \cup_{HS(\nu \times \R)} HC(\nu \times \R).
	\end{align*}
	Here the space $HD(\nu \times \R)$ refers to the lower half-disc of the bundle $\nu \times \R$, and similarly $HS$ is the lower half-sphere and $HC$ is the lower half-complement, as defined in \autoref{half_disc}.
%	\begin{align*}
%		HD(\nu \times \R) &= \{ (x,u) \in \nu \times \R : \|x\|^2 + u^2 \leq 1, \ u \leq 0 \} \\
%		HS(\nu \times \R) &= \{ (x,u) \in \nu \times \R : \|x\|^2 + u^2 = 1, \ u \leq 0 \} \\
%		HC(\nu \times \R) &= \{ (x,u) \in \nu \times \R : \|x\|^2 + u^2 \geq 1, \ \|x\| \leq 1, \ u \in [-1,0] \}
%	\end{align*}
The spaces $QD$, $QS$, and $QC$ are the corresponding quarter versions of these, as defined in \autoref{quarter_disc}.
%	
%\begin{notn}\label{quarter_disc}
%	We let $QD(\nu \times \R^2)$ refer to the lower quarter-disc of the bundle $\nu \times \R^2$, and similarly $QS$ is the lower quarter-sphere and $QC$ is the lower quarter-complement:
%	\begin{align*}
%		QD(\nu \times \R^2) &= \{ (x,u,v) \in \nu \times \R \times \R : \|x\|^2 + t^2 + u^2 \leq 1, \ t \geq 0, \ u \leq 0 \} \\
%		QS(\nu \times \R^2) &= \{ (x,u,v) \in \nu \times \R \times \R : \|x\|^2 + t^2 + u^2 = 1, \ t \geq 0, \ u \leq 0 \} \\
%		QC(\nu \times \R^2) &= \{ (x,u,v) \in \nu \times \R \times \R : \|x\|^2 + t^2 + u^2 \geq 1, \ \|x\| \leq 1, \ t \in [0,1], \ u \in [-1,0] \}
%	\end{align*}
%	In the case that $\nu$ is the trivial 0-dimensional vector bundle over a point, $QD(\R^2)$ is the lower-right quadrant of the unit disc in $\R^2$, $QS(\R^2)$ is the corresponding quadrant of the unit circle, and $QC(\R^2)$ is the complementary region in the square $[0,1] \times [-1,0]$.
%\end{notn}

%	\begin{align*}
%		QD(\R^2) &= \{ (x,y) \in \R^2 : x^2 + y^2 \leq 1, \ x \geq 0, \ y \leq 0 \} \\
%		QS(\R^2) &= \{ (x,y) \in \R^2 : x^2 + y^2 = 1, \ x \geq 0, \ y \leq 0 \} \\
%		QC(\R^2) &= \{ (x,y) \in \R^2 : x^2 + y^2 \geq 1, \ x \in [0,1], \ y \in [-1,0] \}
%	\end{align*}
	We define a spherical blowdown from $B(V \times_M S(\nu))$ to $\St(V)$ by the formula
	\begin{align*}
		(V \times_M S(\nu)) \times_{[0,1]} QD(\R^2) &\to V \times_{M \times [0,1]} HD(\nu \times \R) \\
		(v,\sigma,(t,u)) &\mapsto (v,(t\sigma,u))
	\end{align*}
	This is smooth and full-rank away from $t = 0$. When $V$ is trivial, it becomes
	\begin{align*}
		S(\nu) \times QD(\R^2) &\to HD(\nu \times \R) \\
		(\sigma,(t,u)) &\mapsto (t\sigma,u)
	\end{align*}
	which extends to a map $S(\nu) \times \R^2 \to \nu \times \R$ with the same formula $(\sigma,(t,u)) \mapsto (t\sigma,u)$. Since the lower collar of $V$ trivializes it near the bottom, we conclude that at the bottom of $V$ (where the boundary of the half-disc meets the half-complement), this formula smoothly extends to the map
	\begin{align*}
		S(\nu) \times QC(\R^2) &\to HC(\nu \times \R) \\
		(\sigma,(t,u)) &\mapsto (t\sigma,u).
	\end{align*}
	The upper collar of $V$ also trivializes it near the top, allowing us to see that this formula defines a smooth function at the singular locus as well.
	
	It remains to show this map is a spherical blowdown. This is true when $V$ is trivial by the above analysis, so we only have to consider those points where $v$ is in the interior of $V$ and where $0<|(t,u)|<1$. However in this region the map is a pullback of a map that we already know to be a blowdown. Specifically, there is a pullback diagram of smooth maps
	\[ \xymatrix{
		(V \times_M S(\nu)) \times_{[0,1]} QD(\R^2) \ar[d] \ar[r] & S(\nu) \times QD(\R^2) \ar[d] \\
		V \times_{M \times [0,1]} HD(\nu \times \R) \ar[r] & HD(\nu \times \R)
	} \]
	in which the right-hand map is a spherical blowdown, and the left-hand map is the one we want. By \autoref{pullback_blowup}, it suffices to check that the bottom horizontal map of this diagram is transverse to the exceptional image, in other words the submanifold
	\[ M \times [0,1] \cong HD(M \times \R) \subseteq HD(\nu \times \R) \]
	corresponding to the zero section of $\nu$, times all of $\R$. Each normal vector to this manifold at the point $(x,u) \in M \times [0,1]$ is represented by a vector in the fiber $\xi \in \nu_x$. Without loss of generality, $|\xi| = |u|$, which is positive since $|(t,u)| > 0$ and we are at a point where $t = 0$.
	
	This normal vector is the tangent vector of the path $s \mapsto (x,(\sin s)\xi, (\cos s)u)$ in the space $HD(\nu \times \R)$. We consider this path in particular because the norm of the last two coordinates stays at the constant value $|\xi| = |u|$ throughout the path. This allows us to easily lift the path to any preimage of the point $(x,0,u)$ in the fiber product $V \times_{M \times [0,1]} HD(\nu \times \R)$, just by taking the $V$ coordinate to be constant. This shows that the projection map from the fiber product to $HD(\nu \times \R)$ is transverse to the exceptional image $M \times [0,1]$, and the proof is complete.
%	By \autoref{luck_442}, it suffices to show the two sides have the same isovariant Whitehead torsion. We do this by following the proof of \autoref{stabiso}, but showing that the inclusion of $D(\nu)$ into both $\St^\nu(V)$ and $B(S(U)) \cup_{S(U)} U$ has zero equivariant torsion on the $H$-stratum, for every subgroup $H \leq G$. We also use \autoref{add_in_iso_wh_side_2} here. This shows that the two sides are diffeomorphic.
\end{proof}

\subsection{Stable splitting}\label{stablesplitting}

We define the stabilized space of $h$-cobordisms and then stabilize the splitting \eqref{left_splitting} to obtain the equivalence asserted in \autoref{mainthm}.

Let $\mc U$ be a fixed complete $G$-universe, in other words an orthogonal $G$-representation of countable dimension that contains every irreducible $G$-representation infinitely many times.

For a compact $G$-manifold $M$, consider the diagram $V\mapsto M \times D(V)$ indexed by finite-dimensional representations $V \subseteq \mc U$, in which each inclusion $V \subseteq W$ induces the inclusion (a smooth embedding)
\[ M \times D(V) \subseteq M \times D(W). \]
The colimit
\begin{equation}\label{eq:stable_h}
	\mc H^{\mc U}_G(M) := \underset{V \subseteq \mc U}\colim\, \mc H_{G\dashiso}^c(M \times D(V)).
\end{equation}
is a functor of $M$ when is regarded as an object of $\Man$, the category of smooth manifolds and embeddings defined in \autoref{m_emb}. Note that this colimit is both the stable space of \emph{isovariant} $h$-cobordisms and the stable space of \emph{equivariant} $h$-cobordisms; these are isomorphic by \autoref{heqishiso}. When $G$ is the trivial group, $\mc H^{\mc U}_G(M)$ is the usual stable $h$-cobordism space
\begin{equation}\label{eq:stable_h_noneq}
	\mc H^\infty(M) := \underset{n \to \infty}\colim\, \mc H^c(M \times D^n).
\end{equation}
of a compact manifold $M$.

By \cite[Theorem 7.14 and Corollary 7.18]{GIMM1}, $\mc H^{\mc U}_G$ has a universal extension to all spaces: there is a functor $\widetilde{\mc H}^{\mc U}_G$ on the larger category of all $G$-CW complexes and continuous maps, whose restriction to $\Man$ is naturally weakly equivalent to $\mc H^{\mc U}_G$. In particular, there is a functor $\widetilde{\mc H}^\infty$ defined on all CW complexes and continuous maps, whose restriction to $\Man$ is naturally weakly equivalent to $\mc H^\infty$. 

Because this extension is defined as a homotopy left Kan extension from the category of finite $G$-CW complexes, we have:
\begin{lem}
	The functors $\widetilde{\mc H}^{\mc U}_G$ and $\widetilde{\mc H}^\infty$ preserve filtered homotopy colimits.
\end{lem}

\begin{thm}\label{stable_splitting_thm}
	For any compact smooth $G$-manifold $M$ with corners, there is an equivalence
	\begin{equation*}\label{eq_stable_splitting}
	\prod_{(H) \leq G} \widetilde{\mc H}^\infty(M^{H}_{hWH}) \simeq \mc H^{\mc U}_G(M),
	\end{equation*}
	where $\mc H_G^{\mc U}(M)$ is the stable space of equivariant/isovariant $h$-cobordisms defined in \eqref{eq:stable_h}.
\end{thm}

To begin the proof, fix a sequence of finite-dimensional representations 
$$
V_1 \subseteq V_2 \subseteq \cdots
$$
in $\mc U$ such that each finite-dimensional subspace $V \subseteq \mc U$ is contained in $V_n$ for some $n$. This ensures that the colimit along all representations is equal to the colimit along the representations $V_n$ as $n \to \infty$.

We construct a homotopy-commuting diagram of embeddings of $G$-manifolds
\[ \xymatrix{
	\coprod_{i=1}^k M_{[(H_i)]} \ar[d] \ar[r] & M \ar[d] \\
	\coprod_{i=1}^k \left(M \times D(V_1)\right)_{[(H_i)]} \ar[d] \ar[r] & M \times D(V_1) \ar[d] \\
	\coprod_{i=1}^k \left(M \times D(V_2)\right)_{[(H_i)]} \ar[d] \ar[r] & M \times D(V_2) \ar[d] \\
	\vdots & \vdots
} \]
as follows. The right-hand column of the diagram is given by the standard embeddings
\[ M \to M \times D(V_1) \to M \times D(V_2) \to \cdots \]
These are morphisms in the category $\Man$ of smooth manifolds and embeddings. The left-hand column is given by the functoriality of compactified strata (as defined by blowing up). By \autoref{manstr_equivalent}, there is a functor $\Manstr \to \Man$ that is surjective on objects and an equivalence on morphism spaces, where $\Manstr$ is the category of manifolds with embedded compactified strata. By the surjectivity, for each $n$ the manifold $M \times D(V_n)$ can be equipped with embeddings of its compactified strata $\left(M \times D(V_n)\right)_{[(H_i)]}$. This choice gives the $n$th horizontal map in the diagram above. Because of the equivalence of morphism spaces, each square in the diagram commutes up to isotopy. %The inclusion $M \times D(V_n) \to M \times D(V_{n+1})$ is isotopic to an embedding that comes from $\Manstr$, in other words one which preserves the compactified strata. Doing this for each $n$ gives us a diagram of the above form that strictly commutes, where the maps in the right-hand column are merely isotopic to the standard embeddings $M \times D(V_n) \to M \times D(V_{n+1})$. Therefore if we let the maps in the right-hand column be the standard embeddings then each square commutes up to an isotopy of embeddings.

Applying the simplicial functor $\mc H_{G\dashiso}$ gives the right half of the following diagram \eqref{big_diagram_1} of simplicial sets. The horizontal map from the middle column to the right-hand column is eventually an equivalence, because as soon as $\dim V_n^G \geq 5$ the assumptions of \autoref{unstable_splitting_thm} hold for the manifold $M \times D(V_n)$.

In the left half of the diagram, the top horizontal arrow is the map that consists of first replacing a cobordism on $M_{[(H_i)]}$ by its $H_i$-fixed points to make a cobordism on $M_{[(H_i)]}^{H_i} =M_{[H_i]}$ and then dividing by the free $WH_i$-action to make a cobordism on $M_{[H_i]}/WH_i$. This map is an equivalence, by \autoref{cobordisms_on_free_manifold}. The $n$th horizontal arrow in the left side is defined in the same way, using $M\times D(V_n)$ in place of $M$.

%We extend the diagram to the left using the equivalences from \autoref{cobordisms_on_free_manifold} that restrict every cobordism on $M_{[(H_i)]}$ to its $H_i$-fixed points, which gives a cobordism on $M_{[(H_i)]}^{H_i} =M_{[H_i]}$, and then takes the $WH_i$-orbits, giving a cobordism on $M_{[H_i]}/WH_i$:
\begin{equation}\label{big_diagram_1}
\xymatrix{
	\prod_{i=1}^k \mc H(M_{[H_i]}/WH_i) \ar[d] \ar@{<-}[r]^-\simeq & \prod_{i=1}^k \mc H_{G\dashiso}(M_{[(H_i)]}) \ar[d] \ar[r] & \mc H_{G\dashiso}(M) \ar[d] \\
	\prod_{i=1}^k \mc H((M \times D(V_1))_{[H_i]}/WH_i) \ar[d] \ar@{<-}[r]^-\simeq & \prod_{i=1}^k \mc H_{G\dashiso}(\left(M \times D(V_1)\right)_{[(H_i)]}) \ar[d] \ar[r] & \mc H_{G\dashiso}(M \times D(V_1)) \ar[d] \\
	\prod_{i=1}^k \mc H((M \times D(V_2))_{[H_i]}/WH_i) \ar[d] \ar@{<-}[r]^-\simeq & \prod_{i=1}^k \mc H_{G\dashiso}(\left(M \times D(V_2)\right)_{[(H_i)]}) \ar[d] \ar[r] & \mc H_{G\dashiso}(M \times D(V_2)) \ar[d] \\
	\vdots & \vdots & \vdots %\\
%	\prod_{i=1}^k \mc H^\infty(M^{H_i}_{hWH_i}) \ar[rr]^-\simeq && \mc H_{G\dashiso}^{\mc U}(M)
}
\end{equation}
Each of the regions on the right commutes up to homotopy since $\mc H_{G\dashiso}$ is a simplicially enriched functor. The regions on the left commute strictly by \autoref{stabilization_commutes_with_fixed_points}.
%Since the splitting is also induced by functoriality, it therefore commutes with this stabilization, giving a splitting of $\Hiso^{\mc U}(M)$.

Now take colimits along the columns. Like all sequential colimits of simplicial sets, these are also homotopy colimits. The diagram gives an equivalence from the colimit along the middle column to that of the left-hand column. It also gives an equivalence from the colimit of the middle column to that of the right-hand column, but in this case the equivalence is only given as a homotopy class. (This unfortunate state of affairs will be remedied in the next section.) 

The colimit along the right-hand column is $\mc H_G^{\mc U}(M)$, by definition. According to the following lemma, the colimit along the left-hand column is the product of the non-equivariant stable $h$-cobordism spaces $\widetilde{\mc H}^\infty(M^{H_i}_{hWH_i})$. This will complete the proof of Theorem \ref{stable_splitting_thm}.

\begin{lem}\label{rewrite_colimit} $\underset{n\to\infty}\colim \mc H((M\times D(V_n))_{[H]}/WH)$ is weakly equivalent to $\widetilde{\mc H}^\infty(M^{H}_{hWH})$.
\end{lem}

\begin{proof}
The natural weak equivalence $\mc H^\infty \simeq \widetilde{\mc H}^\infty$ means that in the homotopy category of functors from $\Man$ to spaces, there is a natural transformation $\mc H \to \widetilde{\mc H}^\infty$. Applying this natural transformation to the colimit system, we get a map of colimits
\begin{equation}\label{h_to_h_twiddle}
\underset{n\to\infty}\colim \mc H((M\times D(V_n))_{[H]}/WH) \ra \underset{n\to\infty}\colim \widetilde{\mc H}^\infty((M\times D(V_n))_{[H]}/WH).
\end{equation}
First we check that \eqref{h_to_h_twiddle} is an equivalence. We expand it out as the map from the colimit along the first column of the following diagram to the colimit of the entire diagram. Here to save space we have written $M_n = (M\times D(V_n))_{[H]}/WH$.
\[ \xymatrix{
	\mc H(M_0) \ar[d] \ar[r] & \mc H(M_0 \times D(\R)) \ar[d] \ar[r] & \mc H(M_0 \times D(\R^2)) \ar[d] \ar[r] & \cdots \ar[r] & \mc H^\infty(M_0) \ar[d] \\
	\mc H(M_1) \ar[d] \ar[r] & \mc H(M_1 \times D(\R)) \ar[d] \ar[r] & \mc H(M_1 \times D(\R^2)) \ar[d] \ar[r] & \cdots \ar[r] & \mc H^\infty(M_1) \ar[d] \\
	\mc H(M_2) \ar[d] \ar[r] & \mc H(M_2 \times D(\R)) \ar[d] \ar[r] & \mc H(M_2 \times D(\R^2)) \ar[d] \ar[r] & \cdots \ar[r] & \mc H^\infty(M_2) \ar[d] \\
	\vdots \ar[d] &\vdots \ar[d] &\vdots \ar[d] & \ddots &\vdots \ar[d] \\
	\hocolim \ar[r]^-\simeq & \hocolim \ar[r]^-\simeq & \hocolim \ar[r]^-\simeq & \cdots \ar[r]^-\simeq & \hocolim
} \]
The stabilized functors along the bottom are homotopy functors of $M$, and the maps $M \times D^m \to M \times D^{m+1}$ are homotopy equivalences, giving the equivalences along the bottom row. Composing them together, we see that \eqref{h_to_h_twiddle} is an equivalence.

Next, observe that we have an equivalence of spaces with free $WH$-action
\[ (M\times D(V))_{[H]}\simeq_{WH} (M\times D(V))_H, \]
and, by \cite[Corollary 4.7]{CaryMona3}, the natural map
\[ \underset{n\to\infty}\colim (M\times D(V_n))_{[H]} \to M^H \]
is a Borel equivalence (i.e. a $WH$-equivariant map that is a weak equivalence of spaces after forgetting the group action). Therefore we have an equivalence between the orbits and the homotopy orbits
\[ \underset{n\to\infty}\colim (M\times D(V_n))_{[H]}/WH \overset\sim\longleftarrow \underset{n\to\infty}\colim ((M\times D(V_n))_{[H]})_{hWH} \simar M^H_{hWH}. \]
%We remark, because we will use it later, that this map would be natural in $M$, if we could take the compactifications $(-)_{[H]}$ in a way that was functorial. At any rate, 
Since the functor $\widetilde{\mc H}^\infty$ preserves filtered homotopy colimits, this leads to a weak equivalence
\begin{equation*}\label{h_twiddle_to_hwh}
\underset{n\to\infty}\colim \widetilde{\mc H}^\infty((M\times D(V_n))_{[H]}/WH) \ra \widetilde{\mc H}^\infty(M^H_{hWH}).
\end{equation*}
With \eqref{h_to_h_twiddle} this completes the proof.
%By \cite[Corollary 4.7]{CaryMona3} we have an equivalence
%$$\mc H^\infty(M^H_{hWH}) \simeq \underset{n\to\infty}\colim \mc H^\infty(((M\times D(V_n))_H)/WH).$$
%Since $\mc H^\infty$ is a homotopy functor and we have an equivalence
%
%we can rewrite this as
%$$\mc H^\infty(M^H_{hWH}) \simeq \underset{n\to\infty}\colim \mc H^\infty(((M\times D(V_n))_{[H]})/WH).$$
%Using the definition of $\mc H^\infty$ for compact smooth manifolds from \eqref{eq:stable_h_noneq}, this expands out to
%\begin{align*}
%	\mc H^\infty(M^H_{hWH}) & \simeq \underset{n\to\infty}\colim \underset{m\to\infty}\colim \mc H^c(((M \times D(V_n))_{[H]})/WH \times D(\R^m)) \\
%	& \simeq \underset{n\to\infty}\colim \underset{m\to\infty}\colim \mc H^c(((M\times D(\R^m) \times D(V_n))_{[H]})/WH). \\
%	& \simeq \underset{n\to\infty}\colim \underset{m\to\infty}\colim \mc H^c(((M\times D(\R^m \oplus V_n))_{[H]})/WH).
%\end{align*}
%However the representations $\R^n \oplus V_n$ can be rewritten as representations in $\mc U$ using any fixed isomorphism $\R^\infty \oplus \mc U \cong \mc U$, and so this colimit is equivalent to just a single colimit
%\[ \underset{n\to\infty}\colim \mc H^c(((M\times D(V_n))_{[H]})/WH). \]
%Since the $WH$-action is free, by \autoref{cobordisms_on_free_manifold} we have an equivalence
%$$\mc H_{WH\dashiso}^c( (M\times D(V))_{[H]} )\simeq \mc H^c(((M\times D(V))_{[H]})/WH).$$
%This commutes with the stabilizations by \autoref{stabilization_commutes_with_orbits}, so the colimit is equivalent to
%\[ \underset{n\to\infty}\colim \mc H^c_{WH\dashiso}((M\times D(V_n))_{[H]}) \]
%as desired.
\end{proof}

%All together, the colimits along the columns of \eqref{big_diagram_1} therefore give a map $\prod_{i=1}^k \widetilde{\mc H}^\infty(M^{H_i}_{hWH_i}) \to \mc H^{\mc U}_G(M)$ in the homotopy category, for each manifold $M$.

\subsection{Naturality}\label{sec:naturality}
To upgrade this equivalence to one that is natural in $M$, we use quasicategories and we use the simplicially enriched category $\Manstr$ in a more serious way. The forgetful functor $\Manstr \to \Man$ that forgets the choices of compactified strata is a Dwyer--Kan equivalence by \autoref{manstr_equivalent}. %, so we do not lose or gain any information if we compose the $h$-cobordism functor with this forgetful functor, and regard everything as a functor on $\Manstr$.

From the definition of $\Manstr$ we see that for each value of $i$ there is a functor $\Manstr \to \Man$ taking $M$ to $M_{[H_i]}$, which has a natural transformation back to the forgetful functor $M \mapsto M$. Therefore the top row of the diagram \eqref{big_diagram_1} is a diagram of natural transformations when we regard everything as a functor on the category $\Manstr$ instead of $\Man$. This is enough to conclude that the splitting of $\mc H_{G\dashiso}(M)$ is functorial on $\Manstr$. In other words, the splitting from \autoref{unstablemainthm} is natural, if we regard both sides as functors on $\Manstr$.

To get the naturality of the stable splitting, we move to the setting of quasicategories. For each simplicial category $\mc C$, let $Q(\mc C) = N^s_\bullet(\Ex^\infty \mc C)$ denote the associated quasicategory, obtained by applying $\Ex^\infty$ to the mapping spaces and then taking the simplicial nerve (or homotopy-coherent nerve). Note that $\Ex^\infty\colon \sSet \to \sSet$ is product-preserving, and in particular it preserves the one-point simplicial set. As a result, $Q$ sends products of simplicial categories to products of simplicial sets
\[ Q(\mc C \times \mc D) \cong Q(\mc C) \times Q(\mc D) \]
and $Q$ sends every discrete category $C$ to its nerve
\[ Q(C) \cong N_\bullet(C). \]
In particular, $Q$ sends the walking arrow category $\{0 \to 1\}$ to the 1-simplex $\Delta[1]$.

Let $\mc S = Q(\sSet)$ denote the quasicategory of spaces. This is equivalent to the quasicategory $Q(\mc{CW})$, where $\mc{CW}$ is the category of CW complexes and the simplicial set of continuous maps between them, $S_\bullet\Map(X,Y)$. We let $G\mc{-CW}$ be the corresponding simplicial category of $G$-CW complexes and continuous maps. Finally, let $\mc H_{G\dashiso}\colon \Man \to \sSet$ denote a model of the isovariant $h$-cobordism functor.

Consider the following diagram of simplicial categories and simplicial functors.
\[ \begin{tikzcd}[column sep = 4em]
	\prod_{i=1}^k WH_i\mc{-CW} \ar[bend left = 15]{rrr}{(-)_{hWH_i}} & \prod_{i=1}^k \Man_{\textup{$WH_i$-free}} \rar{/WH_i} \lar[swap]{\textup{forget}} \ar[equals]{rd} & \prod_{i=1}^k \Man_{1} \ar{rd}{\prod \mc H} \rar{\textup{forget}} & \prod_{i=1}^k \mc{CW} \dar{\prod \widetilde{\mc H}^\infty} & \\
	G\mc{-CW} \uar{(-)^{H_i}} & \Manstr_G \rar{M_{[(H_i)]}} \uar{M_{[(H_i)]}} \ar{rd}{M} & \prod_{i=1}^k \Man_{\textup{$WH_i$-free}} \uar{/WH_i} \dar{\amalg} \rar{\mc H_{G\dashiso}} & \prod_{i=1}^k \sSet \dar{\times} \\
	\Man_G \times \N \uar{(M,n) \mapsto M} \rar{(M,n) \mapsto M \times D(V_n)} & \Man_G \uar \rar[equals] & \Man_G \rar{\mc H_{G\dashiso}} & \sSet
\end{tikzcd} \]

Apply $Q$ to the two outside routes of this diagram to get two functors
\[ Q(\Man_G) \times Q(\N) \cong Q(\Man_G \times \N) \ra Q(\sSet). \]
Taking the adjoint $Q(\Man_G) \to \Fun(Q(\N),Q(\sSet))$ and composing with the colimit functor $\Fun(Q(\N),Q(\sSet)) \to Q(\sSet)$ yields the two functors that we wish to compare, namely
\[ \mc H_{G\dashiso}^{\mc U}(-), \quad \prod_{i=1}^k \widetilde{\mc H}^\infty((-)^{H_i}_{hWH_i}) \quad \colon \quad Q(\Man_G) \ra Q(\sSet). \]
\begin{prop}
	These two functors of quasicategories are naturally isomorphic.
\end{prop}

An isomorphism of functors of quasicategories is a natural transformation that is a weak equivalence at each object, so in practical terms this means that the two spaces produced by these functors are naturally weakly equivalent.

\begin{proof}
	We show that each region in this diagram has a natural transformation with the property that if we pre-compose with \emph{any} sequence of functors in the diagram from $\Man_G \times \N$, post-compose with \emph{any} sequence of functors in the diagram that end in $\sSet$, then apply $Q$, adjoint over the $\times Q(\N)$, and take the colimit, the resulting natural transformation is a natural isomorphism of functors $Q(\Man_G) \to Q(\sSet)$.
	
	Using the well-known fact that colimits in $Q(\sSet)$ correspond to homotopy colimits of enriched diagrams in $\sSet$, this means that for each manifold $M$ separately, the natural transformation in each region induces a weak equivalence on homotopy colimits along $\N$.

% Moreover,  with the following property. If we compose with the other functors and then apply $Q$ to get a natural transformation of functors
%\[ Q(\Man_G) \times Q(\N) \cong Q(\Man_G \times \N) \ra Q(\sSet), \]
%then take the adjoint $Q(\Man_G) \to \Fun(Q(\N),Q(\sSet))$ and compose with the colimit functor $\Fun(Q(\N),Q(\sSet)) \to Q(\sSet)$ to get a natural transformation of functors
%\[ Q(\Man_G) \ra Q(\sSet), \]
%that the result is a natural isomorphism between maps of quasicategories.

Now we proceed one region at a time. We comment only on those regions where there is any difficulty:
\begin{itemize}
	\item For the region at the very top, we use the collapse map $X_{hWH_i} \to X/WH_i$, which is a natural transformation of topologically enriched functors and therefore automatically gives a simplicial natural transformation. It gives a weak equivalence of spaces for every object in the source category $\Man_{\textup{$WH_i$-free}}$, so it is already a natural isomorphism as a functor to $\prod_{i=1}^k Q(\mc{CW})$. Therefore it induces a natural isomorphism in quasicategories (i.e. a weak equivalence of spaces) after applying the remaining functors and taking the homotopy colimit.
	\item For the region at the far left, we use the map $(M \times D(V_n))_{[(H_i)]} \to M^{H_i}$ that forgets the $D(V_n)$ coordinate. This is again automatically a simplicial natural transformation. It is not a weak equivalence, but by the proof of \autoref{rewrite_colimit} it induces a weak equivalence once we apply the remaining functors in the diagram and take the homotopy colimit as $n\to\infty$.
	\item In the middle column, one of the triangles corresponds to the inclusion $\coprod_i M_{[(H_i)]} \to M$. This defines a simplicial natural transformation, as we can check directly by building the simplicial functor
	\[ \{0 \to 1\} \times \Manstr \ra \Man \]
	taking $(0,M)$ to $\coprod_i M_{[(H_i)]}$ and $(1,M)$ to $M$, and the morphism $0 \to 1$ to the embedding $\coprod_i M_{[(H_i)]} \to M$. This extends to a functor on the product $\{0 \to 1\} \times \Manstr$ since $\coprod_i M_{[(H_i)]} \to M$ is a natural transformation of simplicial functors $\Manstr \to \Man$. Again this is not an equivalence, but it becomes an equivalence after evaluating at $M \times D(V_n)$, post-composing with $\mc H_{G\dashiso}$, and taking the homotopy colimit as $n \to \infty$. As discussed above \autoref{stable_splitting_thm}, this holds because for sufficiently large $n$, the hypotheses of \autoref{unstable_splitting_thm} apply to $M \times D(V_n)$.
	\item In the right-hand column, the top triangle is the natural transformation between $\mc H$ and $\widetilde{\mc H}$, which is a weak equivalence of simplicial sets for every manifold $M$.
	\item In the right-hand column, the remaining two regions commute up to a simplicial natural isomorphism, by \autoref{stabilization_commutes_with_orbits} and \autoref{stabilization_commutes_with_coproducts}, respectively.
\end{itemize}
\end{proof}

Now we have proven a natural weak equivalence of functors
\[ \mc H^{\mc U}_G(M) \simeq \prod_{(H) \leq G} \widetilde{\mc H}^\infty(M^{H}_{hWH}) , \]
as maps of quasicategories $Q(\Man_G) \to Q(\sSet)$. We can finally take the left Kan extension of these functors from $Q(\Man_G)$ to $Q(G\mc{-CW})$. On the left-hand side, this gives $\widetilde{\mc H}^{\mc U}_G(X)$ by definition. On the right-hand side this gives $\prod_{(H) \leq G} \widetilde{\mc H}^\infty(X^{H}_{hWH})$, but this uses a few facts:
\begin{itemize}
	\item The left Kan extension is expressed at each object as a filtered homotopy colimit, which commutes with finite products.
	\item The functor $\widetilde{\mc H}^\infty$ is itself defined as a left Kan extension from finite spaces, so it commutes with filtered homotopy colimits, and therefore with this Kan extension.
	\item The homotopy orbits functor $(-)_{hWH}$ is a homotopy colimit so it commutes with all homotopy colimits in $WH\mc{-CW}$.
	\item The fixed points $(-)^H$ as a functor on $G\mc{-CW}$ preserves all homotopy colimits. (Note that the weak equivalences in $G\mc{-CW}$ are the maps that are weak equivalences on the fixed point subspaces $X^H$ for all $H \leq G$. So in the corresponding model category of $G$-spaces, this corresponds to the strict fixed points functor $(-)^H$, not the homotopy fixed points functor $(-)^{hH}$.)
\end{itemize}
In conclusion, we have the following more precise form of \autoref{mainthmspaces}:
\begin{thm}\label{natural_stable_splitting_thm}
	There is a natural isomorphism of functors of quasicategories $Q(G\mc{-CW}) \to Q(\sSet)$
	$$\widetilde{\mc H}_G^\mathcal{U}(X)\simeq \prod_{(H) \leq G} \widetilde{\mc H}^\infty(X^H_{hWH} ).$$
\end{thm}
If desired, this could be turned back into a zig-zag of natural weak equivalences between the corresponding simplicially enriched functors $G\mc{-CW} \to \sSet$.

  \bibliographystyle{alpha}
  \bibliography{references}

\begingroup%
\setlength{\parskip}{\storeparskip}% Restore \parskip within this scope

\end{document}